\documentclass[10pt]{article}%
\usepackage{graphicx}
\usepackage{amsmath}
\usepackage{amsfonts}
\usepackage{amssymb}%
\providecommand{\U}[1]{\protect\rule{.1in}{.1in}}
\begin{document}

\title{Energy decay for Maxwell's equations with Ohm's law \\on partially cubic domains \thanks{This work was supported by the NSF of China
under grants 10771149 and 60974035.}}
\author{Kim Dang Phung\\Yangtze Center of Mathematics, Sichuan University, \\Chengdu 610064, Sichuan Province, China. \\{\small \textit{E-mail:}} {\small \texttt{kim\_dang\_phung@yahoo.fr}}}
\date{\quad}
\maketitle

\bigskip

\bigskip

\textbf{Abstract .-} We prove a polynomial energy decay for the Maxwell's
equations with Ohm's law on partially cubic domains with trapped rays.

\bigskip


\bigskip

\bigskip

\section{Introduction}

\bigskip

The problems dealing with Maxwell's equations with nonzero conductivity are
not only theoretical interesting but also very important in many industrial
applications (see e.g. \cite{B}, \cite{DuL}, \cite{KW}).

\bigskip

Let $\Omega$ be a bounded open connected region in $\mathbb{R}^{3}$, with a
smooth boundary $\partial\Omega$. We suppose that $\Omega$ is simply connected
and $\partial\Omega$ has only one connected component. The domain $\Omega$ is
occupied by an electromagnetic medium of constant electric permittivity
$\varepsilon_{o}$ and constant magnetic permeability $\mu_{o}$. Let $E$ and
$H$ denote the electric and magnetic fields respectively. The Maxwell's
equations with Ohm's law are described by
\begin{equation}
\left\{
\begin{array}
[c]{rl}%
\varepsilon_{o}\partial_{t}E-\operatorname{curl}H+\sigma E=0 & \quad
\text{in}~\Omega\times\left[  0,+\infty\right) \\
\mu_{o}\partial_{t}H+\operatorname{curl}E=0 & \quad\text{in}~\Omega
\times\left[  0,+\infty\right) \\
\operatorname{div}\left(  \mu_{o}H\right)  =0 & \quad\text{in}~\Omega
\times\left[  0,+\infty\right) \\
E\times\nu=H\cdot\nu=0 & \quad\text{on}~\partial\Omega\times\left[
0,+\infty\right) \\
\left(  E,H\right)  \left(  \cdot,0\right)  =\left(  E_{o},H_{o}\right)  &
\quad\text{in}~\Omega\text{ .}%
\end{array}
\right.  \tag{1.1}\label{1.1}%
\end{equation}
Here, $\left(  E_{o},H_{o}\right)  $ are the initial data in the energy space
$L^{2}\left(  \Omega\right)  ^{6}$ and $\nu$ denotes the outward unit normal
vector to $\partial\Omega$. The conductivity is such that $\sigma\in
L^{\infty}\left(  \Omega\right)  $ and $\sigma\geq0$. It is well-known that
when the conductivity is identically null, then the above system is
conservative and when $\sigma$ is bounded from below by a positive constant,
then an exponential energy decay rate holds for the Maxwell's equations with
Ohm's law in the energy space. The situation becomes more delicate when we
only assume that
\[
\sigma\left(  x\right)  \geq constant>0\quad\forall x\in\omega
\]
for some non-empty connected open subset $\omega$ of $\Omega$. Observe that
the condition $\operatorname{div}\left(  \varepsilon_{o}E\right)  =0$ in
$\Omega\times\left[  0,+\infty\right)  $ does not appear because the free
divergence is not preserved by the Maxwell's equations with Ohm's law. Here we
know that the above system is dissipative and its energy tends to zero in
large time. However, we would like to establish the energy decay rate as well.
In the field of control theory, the exponential energy decay rate of a linear
dissipative system is deduced from an observability estimate. Precisely, in
order to get an exponential decay rate in the energy space we should have the
following observability inequality%
\[
\exists C,T_{c}>0\quad\forall\zeta\geq0\quad\int_{\Omega}\left\vert \left(
E,H\right)  \left(  \cdot,\zeta\right)  \right\vert ^{2}dx\leq C\int_{\zeta
}^{\zeta+T_{c}}\int_{\Omega}\sigma\left\vert E\right\vert ^{2}dxdt
\]
or simply, in virtue of a semigroup property,
\[
\exists C,T_{c}>0\quad\int_{\Omega}\left\vert \left(  E_{o},H_{o}\right)
\right\vert ^{2}dx\leq C\int_{0}^{T_{c}}\int_{\Omega}\sigma\left\vert
E\right\vert ^{2}dxdt
\]
for any initial data $\left(  E_{o},H_{o}\right)  $ in the energy space
$L^{2}\left(  \Omega\right)  ^{6}$. We can also look for establishing the
above observability inequality for any initial data in the energy space
intersecting suitable invariant subspaces but not with the condition
$\operatorname{div}E_{o}=0$ in $\Omega$. Such estimate is established in
\cite{P1}\ under the geometric control condition of Bardos, Lebeau and Rauch
\cite{BLR} for the scalar wave operator and when the conductivity has the
property that $\sigma\left(  x\right)  \geq constant>0$ for all $x\in\omega$
and $\sigma\left(  x\right)  =0$ for all $x\in\Omega\left\backslash
\overline{\omega}\right.  $. From now, we consider a subset $\omega$ such that
the geometric control condition for the scalar wave operator or other
assumptions based on the multiplier method fail. In such geometry, we do not
hope an exponential energy decay rate in the energy space. Our geometry
(described precisely in Section 3) presents parallel trapped rays and can be
compared to the one in \cite{P2} or in \cite{BH},\cite{N} for the two
dimensional case. It generalises the cube (see \cite{KW}) and therefore
explicit and analytical results are harder to obtain. Our main result gives a
polynomial energy decay with regular initial data. Our proof is based on a new
kind of observation inequality (see (\ref{4.33}) below) which can also be seen
as an interpolation estimate. It relies with the construction of a particular
solution for the operator $i\partial_{s}+h\left(  \Delta-\partial_{t}%
^{2}\right)  $ inspired by the gaussian beam techniques. Also the dispersion
property for the one dimensional Schr\"{o}dinger operator will play a key role.

\bigskip

The plan of the paper is as follows. In the next section, we recall the known
results about the Maxwell's equations with Ohm's law that will be used in the
following. Section 3 contains the statement of our main result, while Section
4 is concerned with its proof. In Section 5, we present the interpolation
estimate, while Section 6 includes its proof. Finally, two appendix are added
dealing with inequalities involving Fourier analysis.

\bigskip

\bigskip

\section{The Maxwell's equations with Ohm's law}

\bigskip

We begin to recall some well-known results concerning the Maxwell's equations
with Ohm's law: well-posedness, energy identity, standard orthogonal
decomposition and asymptotic behaviour in time of the energy of the
electromagnetic field.

\bigskip

\bigskip

\subsection{Well-posedness of the problem}

\bigskip

Let us introduce the spaces%
\begin{equation}
\mathcal{V}=L^{2}\left(  \Omega\right)  ^{3}\times\left\{  G\in L^{2}\left(
\Omega\right)  ^{3};~\operatorname{div}G=0,~G\cdot\nu_{\left\vert
\partial\Omega\right.  }=0\right\}  \text{ ,} \tag{2.1.1}\label{2.1.1}%
\end{equation}%
\begin{equation}%
\begin{array}
[c]{rr}
& \mathcal{W}=\left\{  \left(  F,G\right)  \in L^{2}\left(  \Omega\right)
^{6};~\operatorname{curl}F\in L^{2}\left(  \Omega\right)  ^{3},~F\times
\nu_{\left\vert \partial\Omega\right.  }=0,\right.  \quad\quad\quad\quad
\quad\quad\\
& \left.  \operatorname{div}G=0,~G\cdot\nu_{\left\vert \partial\Omega\right.
}=0,~\operatorname{curl}G\in L^{2}\left(  \Omega\right)  ^{3}\right\}  \text{
.}%
\end{array}
\tag{2.1.2}\label{2.1.2}%
\end{equation}
It is well-known that if $\left(  E_{o},H_{o}\right)  \in\mathcal{V}$, there
is a unique weak solution $\left(  E,H\right)  \in C^{0}\left(  \left[
0,+\infty\right)  ,\mathcal{V}\right)  $. Further, if $\left(  E_{o}%
,H_{o}\right)  \in\mathcal{W}$, there is a unique strong solution $\left(
E,H\right)  \in C^{0}\left(  \left[  0,+\infty\right)  ,\mathcal{W}\right)
\cap C^{1}\left(  \left[  0,+\infty\right)  ,\mathcal{V}\right)  $. Let us
define the functionals of energy%
\begin{equation}
\mathcal{E}\left(  t\right)  =\frac{1}{2}\int_{\Omega}\left(  \varepsilon
_{o}\left\vert E\left(  x,t\right)  \right\vert ^{2}+\mu_{o}\left\vert
H\left(  x,t\right)  \right\vert ^{2}\right)  dx\text{ ,} \tag{2.1.3}%
\label{2.1.3}%
\end{equation}%
\begin{equation}
\mathcal{E}_{1}\left(  t\right)  =\frac{1}{2}\int_{\Omega}\left(
\varepsilon_{o}\left\vert \partial_{t}E\left(  x,t\right)  \right\vert
^{2}+\mu_{o}\left\vert \partial_{t}H\left(  x,t\right)  \right\vert
^{2}\right)  dx\text{ .} \tag{2.1.4}\label{2.1.4}%
\end{equation}
We can easily check that the energy $\mathcal{E}$ is a continuous positive
non-increasing real function on $\left[  0,+\infty\right)  $ and further for
any initial data $\left(  E_{o},H_{o}\right)  \in\mathcal{W}$,
\begin{equation}
\frac{d}{dt}\mathcal{E}\left(  t\right)  +\int_{\Omega}\sigma\left(  x\right)
\left\vert E\left(  x,t\right)  \right\vert ^{2}dx=0\text{ ,} \tag{2.1.5}%
\label{2.1.5}%
\end{equation}
and for any $t_{2}>t_{1}\geq0$,
\begin{equation}
\mathcal{E}\left(  t_{2}\right)  -\mathcal{E}\left(  t_{1}\right)
+\int_{t_{1}}^{t_{2}}\int_{\Omega}\sigma\left(  x\right)  \left\vert E\left(
x,t\right)  \right\vert ^{2}dxdt=0\text{ ,} \tag{2.1.6}\label{2.1.6}%
\end{equation}%
\begin{equation}
\mathcal{E}_{1}\left(  t_{2}\right)  -\mathcal{E}_{1}\left(  t_{1}\right)
+\int_{t_{1}}^{t_{2}}\int_{\Omega}\sigma\left(  x\right)  \left\vert
\partial_{t}E\left(  x,t\right)  \right\vert ^{2}dxdt=0\text{ .}
\tag{2.1.7}\label{2.1.7}%
\end{equation}

\bigskip

\bigskip

\subsection{Orthogonal decomposition}

\bigskip

Both $E$ and $\mu_{o}H$ can be described, by means of the scalar and vector
potentials $p$ and $A$ with the Coulomb gauge, in an unique way as follows.

\bigskip

\textbf{Proposition 2.1}\ -. \textit{For any initial data }$\left(
E_{o},H_{o}\right)  \in\mathcal{W}$\textit{, there is a unique} $\left(
p,A\right)  \in C^{1}\left(  \left[  0,+\infty\right)  ,H_{0}^{1}\left(
\Omega\right)  \right)  \times C^{2}\left(  \left[  0,+\infty\right)
,H^{1}\left(  \Omega\right)  ^{3}\right)  $\textit{ such that }$\left(
E,H\right)  $\textit{ the solution of (\ref{1.1}) the Maxwell's equations with
Ohm's law satisfies}%
\begin{equation}
\left\{
\begin{array}
[c]{ll}%
E=-\nabla p-\partial_{t}A & \\
\mu_{o}H=\operatorname{curl}A &
\end{array}
\right.  \tag{2.2.1}\label{2.2.1}%
\end{equation}%
\begin{equation}
\left\{
\begin{array}
[c]{rl}%
\varepsilon_{o}\mu_{o}\partial_{t}^{2}A+\operatorname{curl}\operatorname{curl}%
A=\mu_{o}\left(  -\varepsilon_{o}\partial_{t}\nabla p+\sigma E\right)  &
\quad\text{\textit{in}}~\Omega\times\left[  0,+\infty\right) \\
\operatorname{div}A=0 & \quad\text{\textit{in}}~\Omega\times\left[
0,+\infty\right) \\
A\times\nu=0 & \quad\text{\textit{on}}~\partial\Omega\times\left[
0,+\infty\right)
\end{array}
\right.  \tag{2.2.2}\label{2.2.2}%
\end{equation}
\textit{and we have the following relations}%
\begin{equation}
\left\Vert E\right\Vert _{L^{2}\left(  \Omega\right)  ^{3}}^{2}=\left\Vert
\nabla p\right\Vert _{L^{2}\left(  \Omega\right)  ^{3}}^{2}+\left\Vert
\partial_{t}A\right\Vert _{L^{2}\left(  \Omega\right)  ^{3}}^{2}\text{
\textit{,}} \tag{2.2.3}\label{2.2.3}%
\end{equation}%
\begin{equation}
\left\Vert \varepsilon_{o}\partial_{t}\nabla p\right\Vert _{L^{2}\left(
\Omega\right)  ^{3}}\leq\left\Vert \sigma E\right\Vert _{L^{2}\left(
\Omega\right)  ^{3}}\text{ \textit{,}} \tag{2.2.4}\label{2.2.4}%
\end{equation}%
\begin{equation}
\exists c>0\qquad\left\Vert A\right\Vert _{L^{2}\left(  \Omega\right)  ^{3}%
}^{2}\leq c\left\Vert \operatorname{curl}A\right\Vert _{L^{2}\left(
\Omega\right)  ^{3}}^{2}\text{ \textit{.}} \tag{2.2.5}\label{2.2.5}%
\end{equation}
\textit{Further, since }$\operatorname{curl}H\in L^{2}\left(  \Omega\right)
^{3}$\textit{,} $\operatorname{curl}\operatorname{curl}A\in L^{2}\left(
\Omega\right)  ^{3}$ \textit{and} $\operatorname{div}A\in H_{0}^{1}\left(
\Omega\right)  $\textit{.}

\bigskip

The proof is essentially given in \cite[page 121]{P1} from a Hodge
decomposition and is omitted here. Now, the vector field $A$ has the nice
property of free divergence and satisfies a second order vector wave equation
with homogeneous boundary condition $A\times\nu=\operatorname{div}A=0$ and
with a second member in $C^{1}\left(  \left[  0,+\infty\right)  ,L^{2}\left(
\Omega\right)  ^{3}\right)  $ bounded by $2\mu_{o}\left\Vert \sigma
E\right\Vert _{L^{2}\left(  \Omega\right)  ^{3}}$. For the sake of simplicity,
we assume from now that $\varepsilon_{o}\mu_{o}=1$.

\bigskip

\bigskip

\subsection{Invariant subspaces, asymptotic behavior and exponential energy
decay}

\bigskip

Let $\omega_{+}$ be a non-empty connected open subset of $\Omega$ with a
Lipschitz boundary $\partial\omega_{+}$. Suppose that $\sigma\in L^{\infty
}\left(  \Omega\right)  $ with the property that $\sigma\left(  x\right)  \geq
constant>0$ for all $x\in\omega_{+}$ and $\sigma\left(  x\right)  =0$ for all
$x\in\Omega\left\backslash \overline{\omega_{+}}\right.  $. Define $\omega
_{-}=\Omega\left\backslash \text{supp}\sigma\right.  $and suppose that its
boundary $\partial\omega_{-}$ is Lipschitz and has no more than two connected
components $\gamma_{1},\gamma_{2}$.

\bigskip

We recall that the range of the curl, $\operatorname{curl}H^{1}\left(
\omega_{-}\right)  ^{3}$, is closed in $L^{2}\left(  \omega_{-}\right)  ^{3}$
(see \cite[page 257]{DaL} or \cite[page 54]{C}) and
\begin{equation}
\operatorname{curl}H^{1}\left(  \omega_{-}\right)  ^{3}=\left\{  U\in
L^{2}\left(  \omega_{-}\right)  ^{3};~\operatorname{div}U=0\text{ in }%
\omega_{-},~\int_{\gamma_{i}}U\cdot\nu=0\quad\text{for }i\in\left\{
1,2\right\}  \right\}  \text{ .} \tag{2.3.1}\label{2.3.1}%
\end{equation}
Its orthogonal space for the $\left(  L^{2}\left(  \omega_{-}\right)  \right)
^{3}$ norm is
\begin{equation}
\left(  \operatorname{curl}H^{1}\left(  \omega_{-}\right)  ^{3}\right)
^{\bot}=\left\{  V\in L^{2}\left(  \omega_{-}\right)  ^{3}%
;~\operatorname{curl}V=0\text{ in }\omega_{-},~V\times\nu=0\text{ on }%
\partial\omega_{-}\right\}  \text{ .} \tag{2.3.2}\label{2.3.2}%
\end{equation}

\bigskip

Let us introduce $\mathcal{S}_{\sigma}=\left(  \operatorname{curl}H^{1}\left(
\omega_{-}\right)  ^{3}\cap L^{2}\left(  \Omega\right)  ^{3}\right)  \times
L^{2}\left(  \Omega\right)  ^{3}$. The space $\mathcal{W}\cap\mathcal{S}%
_{\sigma}$ is stable for the system of Maxwell's equations with Ohm's law,
which can be seen by multiplying by $V\in\left(  \operatorname{curl}%
H^{1}\left(  \omega_{-}\right)  ^{3}\right)  ^{\bot}$ the equation
$\varepsilon_{o}\partial_{t}E-\operatorname{curl}H+\sigma E=0$. Then, we can
add the following well-posedness result. If $\left(  E_{o},H_{o}\right)
\in\mathcal{W}\cap\mathcal{S}_{\sigma}$, there is a unique solution $\left(
E,H\right)  \in C^{0}\left(  \left[  0,+\infty\right)  ,\mathcal{W}%
\cap\mathcal{S}_{\sigma}\right)  \cap C^{1}\left(  \left[  0,+\infty\right)
,\mathcal{V}\cap\mathcal{S}_{\sigma}\right)  $.

\bigskip

It has been proved (see \cite[page 124]{P1}) that if $\omega_{-}$ is a
non-empty connected open set then $\underset{t\rightarrow+\infty}{\lim
}\mathcal{E}\left(  t\right)  =0$ for any initial data $\left(  E_{o}%
,H_{o}\right)  \in\mathcal{V}\cap\mathcal{S}_{\sigma}$. Further, the following
result (see \cite[page 124]{P1}) plays a key role.

\bigskip

\textbf{Proposition 2.2}\ -. \textit{If }$\omega_{-}$\textit{ is a non-empty
connected open set, then there exists }$c>0$\textit{ such that for all initial
data} $\left(  E_{o},H_{o}\right)  \in\mathcal{W}\cap\mathcal{S}_{\sigma}%
$\textit{ of the system (\ref{1.1}) of Maxwell's equations with Ohm's law, we
have}%
\begin{equation}
\forall t\geq0\quad\mathcal{E}\left(  t\right)  \leq c~\mathcal{E}_{1}\left(
t\right)  \text{ .} \tag{2.3.3}\label{2.3.3}%
\end{equation}

\bigskip

\textbf{Remark 2.3}\ -. The estimate (\ref{2.3.3}) is still true if
$\partial\omega_{+}\cap\partial\Omega\neq\emptyset$. Indeed, the proof given
in \cite[page 127]{P1} can be divided into two steps. In the first step, we
begin to establish the existence of $c>0$ such that
\begin{equation}
\left\Vert \nabla p\right\Vert _{L^{2}\left(  \omega_{+}\right)  ^{3}}%
^{2}+\left\Vert \partial_{t}A\right\Vert _{L^{2}\left(  \Omega\right)  ^{3}%
}^{2}+\left\Vert H\right\Vert _{L^{2}\left(  \Omega\right)  ^{3}}^{2}\leq
c\left(  \mathcal{E}_{1}\left(  t\right)  +\sqrt{\mathcal{E}\left(  t\right)
}\sqrt{\mathcal{E}_{1}\left(  t\right)  }\right)  \tag{2.3.4}\label{2.3.4}%
\end{equation}
for any $\left(  E_{o},H_{o}\right)  \in\mathcal{W}\cap\mathcal{S}_{\sigma}$.
Here, we used a standard compactness-uniqueness argument for $H$,
(\ref{2.2.5}) of Proposition 2.1 for $\partial_{t}A$, and for $\nabla p$ from
the fact that $\sigma\left(  x\right)  \geq constant>0$ for all $x\in
\omega_{+}$ and (\ref{2.1.5}). Till now, we did not need that $\omega_{-}$ is
a connected set. The second step (see \cite[page 128]{P1}) did consist to
prove that
\begin{equation}
\left\Vert \nabla p\right\Vert _{L^{2}\left(  \omega_{-}\right)  ^{3}}^{2}\leq
c\left(  \left\Vert \nabla p\right\Vert _{L^{2}\left(  \omega_{+}\right)
^{3}}^{2}+\left\Vert \partial_{t}A\right\Vert _{L^{2}\left(  \Omega\right)
^{3}}^{2}\right)  \text{ .} \tag{2.3.5}\label{2.3.5}%
\end{equation}
Finally, we concluded by virtue of (\ref{2.2.3}) of Proposition 2.1. This last
estimate becomes easier to obtain under the assumption $\partial\omega_{+}%
\cap\partial\Omega\neq\emptyset$ and without adding the hypothesis saying that
$\omega_{-}$ is a connected set. Indeed, since $\left(  E_{o},H_{o}\right)
\in\mathcal{W}\cap\mathcal{S}_{\sigma}$ and $-\Delta p=divE$, $p\in H_{0}%
^{1}\left(  \Omega\right)  $ solves the following elliptic system%
\begin{equation}
\left\{
\begin{array}
[c]{c}%
\Delta p=0\quad\text{in}~\omega_{-}\\
p\in H^{1/2}\left(  \partial\omega_{-}\right) \\
p=0\quad\text{on}~\partial\omega_{+}\cap\partial\Omega\neq\emptyset\text{ .}%
\end{array}
\right.  \tag{2.3.6}\label{2.3.6}%
\end{equation}
Thus, by the elliptic regularity, the trace theorem and the Poincar\'{e}
inequality, we have the following estimate%
\begin{equation}
\left\Vert \nabla p\right\Vert _{L^{2}\left(  \omega_{-}\right)  ^{3}}\leq
c_{1}\left\Vert p\right\Vert _{H^{1/2}\left(  \partial\omega_{+}\right)  }\leq
c_{2}\left\Vert p\right\Vert _{H^{1/2}\left(  \partial\omega_{+}\right)  }\leq
c_{3}\left\Vert \nabla p\right\Vert _{L^{2}\left(  \omega_{+}\right)  ^{3}}
\tag{2.3.7}\label{2.3.7}%
\end{equation}
for suitable constants $c_{1},c_{2},c_{3}>0$. Hence, combining (\ref{2.3.4})
and (\ref{2.3.7}) with (\ref{2.2.3}), (\ref{2.3.3}) follows if $\partial
\omega_{+}\cap\partial\Omega\neq\emptyset$.

\bigskip

The exponential energy decay rate for the Maxwell's equations with Ohm's law
in the energy space is as follows.

\bigskip

\textbf{Proposition 2.4}\ -. \textit{Let }$\vartheta$\textit{ be a subset of
}$\Omega$\textit{ such that any generalized ray of the scalar wave operator
}$\partial_{t}^{2}-\Delta$\textit{ meets }$\overline{\vartheta}$\textit{.
Suppose that }$\overline{\vartheta}\cap\Omega\subset\omega_{+}$\textit{.
Further if }$\partial\omega_{+}\cap\partial\Omega\neq\emptyset$\textit{ or
}$\omega_{-}$\textit{ is a non-empty connected open set, then there exist
}$c>0$\textit{ and }$\beta>0$\textit{ such that for all initial data }$\left(
E_{o},H_{o}\right)  \in\mathcal{V}\cap\mathcal{S}_{\sigma}$\textit{ of the
system (\ref{1.1}) of Maxwell's equations with Ohm's law, we have}%
\begin{equation}
\forall t\geq0\quad\mathcal{E}\left(  t\right)  \leq ce^{-\beta t}%
~\mathcal{E}\left(  0\right)  \text{ \textit{.}} \tag{2.3.8}\label{2.3.8}%
\end{equation}

\bigskip

The proof of Proposition 2.4 is done in \cite[page 129]{P1} when $\omega_{-}$
is a non-empty connected open set. Here, we simply recall the key points of
the proof. From the geometric control condition, the following estimate holds
without using the fact that $\omega_{-}$ is a non-empty connected open set.%
\begin{equation}
\exists C,T_{c}>0\quad\forall\zeta\geq0\quad\mathcal{E}_{1}\left(
\zeta\right)  \leq C\int_{\zeta}^{T_{c}+\zeta}\int_{\Omega}\left(
\sigma\left\vert \partial_{t}E\right\vert ^{2}+\sigma\left\vert E\right\vert
^{2}\right)  dxdt \tag{2.3.9}\label{2.3.9}%
\end{equation}
for any initial data $\left(  E_{o},H_{o}\right)  \in\mathcal{W}%
\cap\mathcal{S}_{\sigma}$. Next by (\ref{2.3.3}), we deduced that
\begin{equation}
\exists C,T_{c}>0\quad\forall\zeta\geq0\quad\mathcal{E}\left(  \zeta\right)
+\mathcal{E}_{1}\left(  \zeta\right)  \leq C\int_{\zeta}^{T_{c}+\zeta}%
\int_{\Omega}\left(  \sigma\left\vert \partial_{t}E\right\vert ^{2}%
+\sigma\left\vert E\right\vert ^{2}\right)  dxdt\text{ .} \tag{2.3.10}%
\label{2.3.10}%
\end{equation}
Finally, we concluded by virtue of a semigroup property. The proof works as
well when $\partial\omega_{+}\cap\partial\Omega\neq\emptyset$ thanks to Remark 2.3.

\bigskip

\bigskip

\section{Geometric setting and main result}

\bigskip

Let us introduce the geometry on which we work in this paper.

\bigskip

We set $D\left(  r_{1},r_{2}\right)  =\left\{  \left(  x_{1},x_{2}\right)
\in\mathbb{R}^{2};\left\vert x_{1}\right\vert <r_{1},\left\vert x_{2}%
\right\vert <r_{2}\right\}  $ where $r_{1},r_{2}>0$. Let $m_{1},m_{2},\rho>0$.
We choose $\Omega$ a connected open set in $\mathbb{R}^{3}$ bounded by
$\Gamma_{1}$, $\Gamma_{2}$, $\Upsilon$ where

\begin{description}
\item $\quad\Gamma_{1}=\overline{D\left(  m_{1},m_{2}\right)  }\times\left\{
\rho\right\}  $, with boundary $\partial\Gamma_{1}$,

\item $\quad\Gamma_{2}=\overline{D\left(  m_{1},m_{2}\right)  }\times\left\{
-\rho\right\}  $, with boundary $\partial\Gamma_{2}$,

\item $\quad\Upsilon$ is a surface with boundary $\partial\Upsilon
=\partial\Gamma_{1}\cup\partial\Gamma_{2}$.
\end{description}

\noindent Therefore, the boundary of $\Omega$ is $\partial\Omega=\Gamma
_{1}\cup\Gamma_{2}\cup\Upsilon$. Further, we suppose that $\partial\Omega$ is
$C^{\infty}$ with $\Upsilon\subset\left(  \mathbb{R}^{2}\left\backslash
D\left(  m_{1},m_{2}\right)  \right.  \right)  \times\mathbb{R}$. In
particular, $\Omega$ is simply connected and $\partial\Omega$ has only one
connected component.

\bigskip

\noindent Let $\Theta$ be a small neighborhood of $\Upsilon$ in $\mathbb{R}%
^{3}$ such that $\Theta\cap D\left(  M_{1},M_{2}\right)  \times\left[
-\rho,\rho\right]  =\emptyset$ for some $M_{1}\in\left(  0,m_{1}\right)  $ and
$M_{2}\in\left(  0,m_{2}\right)  $. Further, we suppose that the boundaries
$\partial\left(  \Theta\cap\Omega\right)  $ and $\partial\left(
\Omega\left\backslash \overline{\Theta\cap\Omega}\right.  \right)  $ are at
least Lipschitz.

\bigskip

After these preparations, we are now able to state our main result.

\bigskip

\textbf{Theorem 3.1} .- \textit{Let }$\omega=\Theta\cap\Omega$\textit{. If
}$\sigma\in L^{\infty}\left(  \Omega\right)  $\textit{ is such that }%
$\sigma\left(  x\right)  \geq constant>0$\textit{ for all }$x\in\omega
$\textit{ and }$\sigma\left(  x\right)  =0$\textit{ for all }$x\in
\Omega\left\backslash \overline{\omega}\right.  $\textit{, then there exist
}$c>0$\textit{ and }$\gamma>0$\textit{ such that for any }$t\geq0$%
\begin{equation}
\mathcal{E}\left(  t\right)  \leq\frac{c}{t^{\gamma}}\left(  \mathcal{E}%
\left(  0\right)  +\mathcal{E}_{1}\left(  0\right)  \right)  \tag{3.1}%
\label{3.1}%
\end{equation}
\textit{for every solution of the system (\ref{1.1}) of Maxwell's equations
with Ohm's law with initial data }$\left(  E_{o},H_{o}\right)  $\textit{ in
}$\mathcal{W}\cap\mathcal{S}_{\sigma}$.

\bigskip

\textbf{Remark 3.2}\ -. From now $\omega=\Theta\cap\Omega$ and $\sigma\in
L^{\infty}\left(  \Omega\right)  $ is such that $\sigma\left(  x\right)  \geq
constant>0$ for all $x\in\omega$ and $\sigma\left(  x\right)  =0$ for all
$x\in\Omega\left\backslash \overline{\omega}\right.  $. Notice that $\omega$
and $\Omega\left\backslash \overline{\omega}\right.  $ are two non-empty
connected open sets with Lipschitz boundaries. Therefore by Proposition 2.2,
there exists $c>0$ such that for all initial data $\left(  E_{o},H_{o}\right)
\in\mathcal{W}\cap\mathcal{S}_{\sigma}$ of the system (\ref{1.1}) of Maxwell's
equations with Ohm's law, $\mathcal{E}\left(  t\right)  \leq c~\mathcal{E}%
_{1}\left(  t\right)  $ for any $t\geq0$.

\bigskip

\textbf{Remark 3.3}\ -. Notice the existence of trapped rays bouncing up and
down from $\Gamma_{1}$ to $\Gamma_{2}$.

\bigskip

\bigskip

\section{Proof of the main result}

\bigskip

Let us consider the solution $U$ of the following system%
\begin{equation}
\left\{
\begin{array}
[c]{rl}%
\partial_{t}^{2}U+\operatorname{curl}\operatorname{curl}U=0 & \quad
\text{in}~\Omega\times\mathbb{R}\\
\operatorname{div}U=0 & \quad\text{in}~\Omega\times\mathbb{R}\\
U\times\nu=0 & \quad\text{on}~\partial\Omega\times\mathbb{R}\\
\left(  U\left(  \cdot,0\right)  ,\partial_{t}U\left(  \cdot,0\right)
\right)  =\left(  U^{0},U^{1}\right)  & \quad\text{in}~\Omega\text{ ,}\\
\left(  U^{0},U^{1}\right)  \in\mathcal{X} & \text{,}\\
\left(  U^{1},\operatorname{curl}\operatorname{curl}U^{0}\right)
\in\mathcal{X} & \text{,}%
\end{array}
\right.  \tag{4.1}\label{4.1}%
\end{equation}
where%
\begin{equation}
\mathcal{X}=\left\{  \left(  F,G\right)  \in L^{2}\left(  \Omega\right)
^{6};~\operatorname{curl}F\in L^{2}\left(  \Omega\right)  ^{3},F\times
\nu_{\left\vert \partial\Omega\right.  }=0,\operatorname{div}%
F=0,\operatorname{div}G=0\right\}  \text{ .} \tag{4.2}\label{4.2}%
\end{equation}
It is well-known that the above system is well-posed with a unique solution
$U$ such that $\left(  U\left(  \cdot,t\right)  ,\partial_{t}U\left(
\cdot,t\right)  \right)  $ and $\left(  \partial_{t}U\left(  \cdot,t\right)
,\partial_{t}^{2}U\left(  \cdot,t\right)  \right)  $ belong to $\mathcal{X}$
for any $t\in\mathbb{R}$. Let us define the following two conservations of
energies.%
\begin{equation}
\mathcal{G}\left(  U,0\right)  =\mathcal{G}\left(  U,t\right)  \equiv
\int_{\Omega}\left(  \left\vert \partial_{t}U\left(  x,t\right)  \right\vert
^{2}+\left\vert \operatorname{curl}U\left(  x,t\right)  \right\vert
^{2}\right)  dx\text{ ,} \tag{4.3}\label{4.3}%
\end{equation}%
\begin{equation}
\mathcal{G}\left(  \partial_{t}U,0\right)  =\mathcal{G}\left(  \partial
_{t}U,t\right)  \equiv\int_{\Omega}\left(  \left\vert \operatorname{curl}%
\operatorname{curl}U\left(  x,t\right)  \right\vert ^{2}+\left\vert
\operatorname{curl}\partial_{t}U\left(  x,t\right)  \right\vert ^{2}\right)
dx\text{ .} \tag{4.4}\label{4.4}%
\end{equation}
Further, for such solution $U$, the following two inequalities hold by
standard compactness-uniqueness argument and classical embedding (see
\cite{ABDG} and \cite[page 50]{C}).%
\begin{equation}
\mathcal{G}\left(  U,t\right)  \leq c\mathcal{G}\left(  \partial
_{t}U,t\right)  \text{ ,} \tag{4.5}\label{4.5}%
\end{equation}%
\begin{equation}
\left\Vert U\left(  \cdot,t\right)  \right\Vert _{H^{1}\left(  \Omega\right)
^{3}}^{2}\leq c\left\Vert \operatorname{curl}U\left(  \cdot,t\right)
\right\Vert _{L^{2}\left(  \Omega\right)  ^{3}}^{2}\text{ ,} \tag{4.6}%
\label{4.6}%
\end{equation}
for some $c>0$ and any $t\in\mathbb{R}$.

\bigskip

Since $\Theta$ is a small neighborhood of $\Upsilon$ in $\mathbb{R}^{3}$ such
that $\Theta\cap D\left(  M_{1},M_{2}\right)  \times\left[  -\rho,\rho\right]
=\emptyset$ for some $M_{1}\in\left(  0,m_{1}\right)  $ and $M_{2}\in\left(
0,m_{2}\right)  $, there exists a positive real number $r_{o}<\min\left(
m_{1}-M_{1},m_{2}-M_{2},\rho\right)  /2$ such that $D\left(  m_{1}%
,m_{2}\right)  \left\backslash D\left(  m_{1}-2r_{o},m_{2}-2r_{o}\right)
\right.  \times\left(  \rho-2r_{o},\rho+2r_{o}\right)  \cup\left(
-\rho-2r_{o},-\rho+2r_{o}\right)  \subset\Theta$. Now, we define
\begin{equation}
\omega_{o}=D\left(  m_{1}-r_{o},m_{2}-r_{o}\right)  \times\left(  \rho
-2r_{o},\rho-r_{o}\right)  \text{ .} \tag{4.7}\label{4.7}%
\end{equation}

\bigskip

\textbf{Proposition 4.1}\ -. \textit{There exist }$h_{o},c,\gamma>0$\textit{
such that for any }$T_{o}>0$\textit{, }$\zeta\geq0$\textit{ and }$h\in\left(
0,h_{o}\right]  $\textit{, the solution }$U$\textit{ of (\ref{4.1}) satisfies}%
\begin{equation}
\int_{\zeta+c\frac{1}{h^{\gamma}}}^{T_{o}+\zeta+c\frac{1}{h^{\gamma}}}%
\int_{\omega_{o}}\left\vert \partial_{t}U\right\vert ^{2}dxdt\leq c\frac
{1}{h^{\gamma}}\int_{\zeta}^{\zeta+2c\frac{1}{h^{\gamma}}}\int_{\omega}\left(
\left\vert \partial_{t}U\right\vert ^{2}+\left\vert U\right\vert ^{2}\right)
dxdt+ch\mathcal{G}\left(  \partial_{t}U,\zeta\right)  \text{ \textit{.}}
\tag{4.8}\label{4.8}%
\end{equation}

\bigskip

We shall leave the proof of Proposition 4.1 till later (see Section 5). Now we
turn to prove Theorem 3.1.

\bigskip

We start by choosing $\widetilde{\omega_{o}}\subset\omega_{o}$ such that
$\omega\cup\widetilde{\omega_{o}}$ is a non-empty connected open set and such
that the boundaries $\partial\left(  \omega\cup\widetilde{\omega_{o}}\right)
$ and $\partial\left(  \Omega\left\backslash \overline{\omega\cup
\widetilde{\omega_{o}}}\right.  \right)  $ are Lipschitz. Notice that
$\partial\left(  \omega\cup\widetilde{\omega_{o}}\right)  \cap\partial
\Omega\neq\emptyset$ and there exists $\vartheta$ a subset of $\Omega$ such
that $\overline{\vartheta}\cap\Omega\subset\left(  \omega\cup\widetilde
{\omega_{o}}\right)  $ and such that any generalized ray of the scalar wave
operator $\partial_{t}^{2}-\Delta$ meets $\overline{\vartheta}$.

\bigskip

Let $\zeta,T_{h}\geq0$. Let $\left(  \widetilde{E},\widetilde{H}\right)  $
denote the electromagnetic field of the following Maxwell's equations with
Ohm's law
\begin{equation}
\left\{
\begin{array}
[c]{rl}%
\varepsilon_{o}\partial_{t}\widetilde{E}-\operatorname{curl}\widetilde
{H}+\left(  \sigma+1_{\left\vert \widetilde{\omega_{o}}\right.  }\right)
\widetilde{E}=0 & \quad\text{in}~\Omega\times\left[  0,+\infty\right) \\
\mu_{o}\partial_{t}\widetilde{H}+\operatorname{curl}\widetilde{E}=0 &
\quad\text{in}~\Omega\times\left[  0,+\infty\right) \\
\operatorname{div}\left(  \mu_{o}\widetilde{H}\right)  =0 & \quad
\text{in}~\Omega\times\left[  0,+\infty\right) \\
\widetilde{E}\times\nu=\widetilde{H}\cdot\nu=0 & \quad\text{on}~\partial
\Omega\times\left[  0,+\infty\right) \\
\left(  \widetilde{E},\widetilde{H}\right)  \left(  \cdot,\zeta+T_{h}\right)
=\left(  E,H\right)  \left(  \cdot,\zeta+T_{h}\right)  & \quad\text{in}%
~\Omega\text{ .}%
\end{array}
\right.  \tag{4.9}\label{4.9}%
\end{equation}
The conductivity $\left(  \sigma+1_{\left\vert \widetilde{\omega_{o}}\right.
}\right)  $ is such that $\left(  \sigma+1_{\left\vert \widetilde{\omega_{o}%
}\right.  }\right)  \geq constant>0$ in $\omega\cup\widetilde{\omega_{o}}$ and
$\left(  \sigma+1_{\left\vert \widetilde{\omega_{o}}\right.  }\right)  =0$ in
$\Omega\left\backslash \left(  \overline{\omega\cup\widetilde{\omega_{o}}%
}\right)  \right.  $. Also notice that $\left(  E,H\right)  \in\mathcal{W}%
\cap\mathcal{S}_{\sigma}\subset\left(  \mathcal{V}\cap\mathcal{S}%
_{\sigma+1_{\left\vert \widetilde{\omega_{o}}\right.  }}\right)  $. Therefore
by Proposition 2.4, there exist $c,\beta>0$ (independent of $\zeta,T_{h}$)
such that for any $t\geq0$ we have
\begin{equation}
\int_{\Omega}\left(  \varepsilon_{o}\left\vert \widetilde{E}\left(
x,t+\zeta+T_{h}\right)  \right\vert ^{2}+\mu_{o}\left\vert \widetilde
{H}\left(  x,t+\zeta+T_{h}\right)  \right\vert ^{2}\right)  dx\leq ce^{-\beta
t}~\mathcal{E}\left(  \zeta+T_{h}\right)  \text{ .} \tag{4.10}\label{4.10}%
\end{equation}
On the other hand, let $\left(  \overline{E},\overline{H}\right)  =\left(
\widetilde{E},\widetilde{H}\right)  -\left(  E,H\right)  $. Then it solves
\begin{equation}
\left\{
\begin{array}
[c]{rl}%
\varepsilon_{o}\partial_{t}\overline{E}-\operatorname{curl}\overline
{H}+\left(  \sigma+1_{\left\vert \widetilde{\omega_{o}}\right.  }\right)
\overline{E}=-1_{\left\vert \widetilde{\omega_{o}}\right.  }E & \quad
\text{in}~\Omega\times\left[  0,+\infty\right) \\
\mu_{o}\partial_{t}\overline{H}+\operatorname{curl}\overline{E}=0 &
\quad\text{in}~\Omega\times\left[  0,+\infty\right) \\
\operatorname{div}\left(  \mu_{o}\overline{H}\right)  =0 & \quad
\text{in}~\Omega\times\left[  0,+\infty\right) \\
\overline{E}\times\nu=\overline{H}\cdot\nu=0 & \quad\text{on}~\partial
\Omega\times\left[  0,+\infty\right) \\
\left(  \overline{E}\left(  \cdot,\zeta+T_{h}\right)  ,\overline{H}\left(
\cdot,\zeta+T_{h}\right)  \right)  =\left(  0,0\right)  & \quad\text{in}%
~\Omega\text{ ,}%
\end{array}
\right.  \tag{4.11}\label{4.11}%
\end{equation}
and by a standard energy method and the fact that $\widetilde{\omega_{o}%
}\subset\omega_{o}$, we get that for any $t\geq0$
\begin{equation}
\int_{\Omega}\left(  \varepsilon_{o}\left\vert \overline{E}\left(
x,t+\zeta+T_{h}\right)  \right\vert ^{2}+\mu_{o}\left\vert \overline{H}\left(
x,t+\zeta+T_{h}\right)  \right\vert ^{2}\right)  dx\leq\frac{t}{\varepsilon
_{o}}\int_{\zeta+T_{h}}^{t+\zeta+T_{h}}\int_{\omega_{o}}\left\vert E\left(
x.s\right)  \right\vert ^{2}dxds\text{ .} \tag{4.12}\label{4.12}%
\end{equation}
Now we are able to bound the quantity $\mathcal{E}\left(  \zeta+T_{h}\right)
=\mathcal{E}\left(  t+\zeta+T_{h}\right)  +\int_{\zeta+T_{h}}^{t+\zeta+T_{h}%
}\int_{\Omega}\sigma\left(  x\right)  \left\vert E\left(  x,s\right)
\right\vert ^{2}dxds$ as follows. By using (\ref{4.10}) and (\ref{4.12}), we
deduce that
\begin{equation}
\mathcal{E}\left(  \zeta+T_{h}\right)  \leq2ce^{-\beta t}~\mathcal{E}\left(
\zeta+T_{h}\right)  +\frac{2t}{\varepsilon_{o}}\int_{\zeta+T_{h}}%
^{t+\zeta+T_{h}}\int_{\omega_{o}}\left\vert E\left(  x.s\right)  \right\vert
^{2}dxds+\int_{\zeta+T_{h}}^{t+\zeta+T_{h}}\int_{\Omega}\sigma\left(
x\right)  \left\vert E\left(  x,s\right)  \right\vert ^{2}dxds \tag{4.13}%
\label{4.13}%
\end{equation}
which implies by taking $t$ large enough, the existence of constants
$C,T_{c}>1$ such that
\begin{equation}
\mathcal{E}\left(  \zeta+T_{h}\right)  \leq C\int_{\zeta+T_{h}}^{T_{c}%
+\zeta+T_{h}}\left(  \int_{\Omega}\sigma\left\vert E\right\vert ^{2}%
+\int_{\omega_{o}}\left\vert E\right\vert ^{2}dx\right)  dxdt\text{ .}
\tag{4.14}\label{4.14}%
\end{equation}

\bigskip

Recall the existence of the vector potential $A$ from Proposition 2.1 and let
$U$ be the solution of%
\begin{equation}
\left\{
\begin{array}
[c]{rl}%
\partial_{t}^{2}U+\operatorname{curl}\operatorname{curl}U=0 & \quad
\text{in}~\Omega\times\mathbb{R}\\
\operatorname{div}U=0 & \quad\text{in}~\Omega\times\mathbb{R}\\
U\times\nu=0 & \quad\text{on}~\partial\Omega\times\mathbb{R}\\
\left(  U,\partial_{t}U\right)  \left(  \cdot,\zeta\right)  =\left(
A,\partial_{t}A\right)  \left(  \cdot,\zeta\right)  & \quad\text{in}%
~\Omega\text{ ,}%
\end{array}
\right.  \tag{4.15}\label{4.15}%
\end{equation}
then by a standard energy method, for any $T_{1}>0$%
\begin{equation}
\int_{\zeta}^{\zeta+T_{1}}\int_{\Omega}\left(  \left\vert \partial_{t}\left(
U-A\right)  \right\vert ^{2}+\left\vert \operatorname{curl}\left(  U-A\right)
\right\vert ^{2}\right)  dxdt\leq T_{1}^{2}\int_{\zeta}^{\zeta+T_{1}}%
\int_{\Omega}\mu_{o}\left\vert -\varepsilon_{o}\partial_{t}\nabla p+\sigma
E\right\vert ^{2}dxdt \tag{4.16}\label{4.16}%
\end{equation}
which implies from (\ref{2.2.4}) of Proposition 2.1\ that%
\begin{equation}
\int_{\zeta}^{\zeta+T_{1}}\int_{\Omega}\left(  \left\vert \partial_{t}\left(
U-A\right)  \right\vert ^{2}+\left\vert \operatorname{curl}\left(  U-A\right)
\right\vert ^{2}\right)  dxdt\leq4\mu_{o}T_{1}^{2}\int_{\zeta}^{\zeta+T_{1}%
}\int_{\Omega}\left\vert \sigma E\right\vert ^{2}dxdt\text{ .} \tag{4.17}%
\label{4.17}%
\end{equation}

\bigskip

Now we are able to bound the quantity $\mathcal{E}\left(  \zeta\right)
=\mathcal{E}\left(  \zeta+T_{h}\right)  +\int_{\zeta}^{\zeta+T_{h}}%
\int_{\Omega}\sigma\left(  x\right)  \left\vert E\left(  x,s\right)
\right\vert ^{2}dxds$ as follows. Since $E=-\nabla p+\partial_{t}\left(
U-A\right)  -\partial_{t}U$, we deduce by (\ref{4.14}) and (\ref{4.17}) that%
\begin{equation}%
\begin{array}
[c]{ll}%
\mathcal{E}\left(  \zeta\right)  & \leq C%
{\displaystyle\int_{\zeta+T_{h}}^{T_{c}+\zeta+T_{h}}}
\left(
{\displaystyle\int_{\Omega}}
\sigma\left\vert E\right\vert ^{2}dxdt+%
{\displaystyle\int_{\omega_{o}}}
\left\vert -\nabla p+\partial_{t}\left(  U-A\right)  -\partial_{t}U\right\vert
^{2}dx\right)  dt+%
{\displaystyle\int_{\zeta}^{\zeta+T_{h}}}
{\displaystyle\int_{\Omega}}
\sigma\left\vert E\right\vert ^{2}dxdt\\
& \leq C\left(  1+T_{h}^{2}\right)
{\displaystyle\int_{\zeta}^{T_{c}+\zeta+T_{h}}}
{\displaystyle\int_{\Omega}}
\left(  \sigma\left\vert E\right\vert ^{2}+\left\vert \nabla p\right\vert
^{2}\right)  dxdt+C%
{\displaystyle\int_{\zeta+T_{h}}^{T_{c}+\zeta+T_{h}}}
{\displaystyle\int_{\omega_{o}}}
\left\vert \partial_{t}U\right\vert ^{2}dxdt\text{ .}%
\end{array}
\tag{4.18}\label{4.18}%
\end{equation}
Here and hereafter, $C$ will be used to denote a generic constant, not
necessarily the same in any two places.

\bigskip

Now we fix $T_{h}=c\frac{1}{h^{\gamma}}$ where $c$ and $\gamma$ are given by
Proposition 4.1. Taking $T_{o}=T_{c}$ in Proposition 4.1, we obtain that for
any $\zeta\geq0$ and $h\in\left(  0,h_{o}\right]  $,
\begin{equation}
\int_{\zeta+T_{h}}^{T_{c}+\zeta+T_{h}}\int_{\omega_{o}}\left\vert \partial
_{t}U\right\vert ^{2}dxdt\leq CT_{h}\int_{\zeta}^{\zeta+2T_{h}}\int_{\omega
}\left(  \left\vert \partial_{t}U\right\vert ^{2}+\left\vert U\right\vert
^{2}\right)  dxdt+ch\mathcal{G}\left(  \partial_{t}U,\zeta\right)  \text{ .}
\tag{4.19}\label{4.19}%
\end{equation}
But%
\begin{equation}%
\begin{array}
[c]{ll}%
{\displaystyle\int_{\zeta}^{\zeta+2T_{h}}}
{\displaystyle\int_{\omega}}
\left\vert \partial_{t}U\right\vert ^{2}dxdt & =%
{\displaystyle\int_{\zeta}^{\zeta+2T_{h}}}
{\displaystyle\int_{\omega}}
\left\vert -E-\nabla p+\partial_{t}\left(  U-A\right)  \right\vert ^{2}dxdt\\
& \leq CT_{h}^{2}%
{\displaystyle\int_{\zeta}^{\zeta+2T_{h}}}
{\displaystyle\int_{\Omega}}
\left(  \sigma\left\vert E\right\vert ^{2}+\left\vert \nabla p\right\vert
^{2}\right)  dxdt
\end{array}
\tag{4.20}\label{4.20}%
\end{equation}
and%
\begin{equation}%
\begin{array}
[c]{ll}%
{\displaystyle\int_{\zeta}^{\zeta+2T_{h}}}
{\displaystyle\int_{\omega}}
\left\vert U\right\vert ^{2}dxdt & \leq2%
{\displaystyle\int_{\zeta}^{\zeta+2T_{h}}}
{\displaystyle\int_{\Omega}}
\left\vert U-A\right\vert ^{2}dxdt+2%
{\displaystyle\int_{\zeta}^{\zeta+2T_{h}}}
{\displaystyle\int_{\omega}}
\left\vert A\right\vert ^{2}dxdt\\
& \leq C%
{\displaystyle\int_{\zeta}^{\zeta+2T_{h}}}
{\displaystyle\int_{\Omega}}
\left\vert \operatorname{curl}\left(  U-A\right)  \right\vert ^{2}dxdt+2%
{\displaystyle\int_{\zeta}^{\zeta+2T_{h}}}
{\displaystyle\int_{\omega}}
\left\vert A\right\vert ^{2}dxdt\\
& \leq CT_{h}^{2}%
{\displaystyle\int_{\zeta}^{\zeta+2T_{h}}}
{\displaystyle\int_{\Omega}}
\sigma\left\vert E\right\vert ^{2}dxdt+2%
{\displaystyle\int_{\zeta}^{\zeta+2T_{h}}}
{\displaystyle\int_{\omega}}
\left\vert A\right\vert ^{2}dxdt
\end{array}
\tag{4.21}\label{4.21}%
\end{equation}
therefore (\ref{4.19}) becomes
\begin{equation}
\int_{\zeta+T_{h}}^{T_{c}+\zeta+T_{h}}\int_{\omega_{o}}\left\vert \partial
_{t}U\right\vert ^{2}dxdt\leq CT_{h}^{3}\int_{\zeta}^{\zeta+2T_{h}}\left(
\int_{\Omega}\left(  \sigma\left\vert E\right\vert ^{2}+\left\vert \nabla
p\right\vert ^{2}\right)  dx+\int_{\omega}\left\vert A\right\vert
^{2}dx\right)  dt+ch\mathcal{G}\left(  \partial_{t}U,\zeta\right)
\tag{4.22}\label{4.22}%
\end{equation}
and finally, combining (\ref{4.18}) and (\ref{4.22}), we get%
\begin{equation}
\mathcal{E}\left(  \zeta\right)  \leq CT_{h}^{3}\int_{\zeta}^{\zeta+2T_{h}%
}\left(  \int_{\Omega}\left(  \sigma\left\vert E\right\vert ^{2}+\left\vert
\nabla p\right\vert ^{2}\right)  dx+\int_{\omega}\left\vert A\right\vert
^{2}dx\right)  dt+ch\mathcal{G}\left(  \partial_{t}U,\zeta\right)  \text{ .}
\tag{4.23}\label{4.23}%
\end{equation}

\bigskip

We have proved that there exist $h_{o},c,\gamma>0$ such that for any
$\zeta\geq0$ and $h\in\left(  0,h_{o}\right]  $, the solution $\left(
E,H\right)  $ of (\ref{1.1}) satisfies%
\begin{equation}
\mathcal{E}\left(  \zeta\right)  \leq c\frac{1}{h^{\gamma}}\int_{\zeta}%
^{\zeta+c\frac{1}{h^{\gamma}}}\left(  \int_{\Omega}\left(  \sigma\left\vert
E\right\vert ^{2}+\left\vert \nabla p\right\vert ^{2}\right)  dx+\int_{\omega
}\left\vert A\right\vert ^{2}dx\right)  dt+ch\mathcal{G}\left(  \partial
_{t}U,\zeta\right)  \text{ .} \tag{4.24}\label{4.24}%
\end{equation}
By formula (\ref{2.1.6}) and since $\mathcal{G}\left(  \partial_{t}%
U,\zeta+mc\frac{1}{h^{\gamma}}\right)  =\mathcal{G}\left(  \partial_{t}%
U,\zeta\right)  $ for any $m$, this last inequality becomes
\begin{equation}%
\begin{array}
[c]{ll}
& \quad N\mathcal{E}\left(  \zeta\right)  -\sum\limits_{m=0,..,N-1}%
{\displaystyle\int_{\zeta}^{\zeta+mc\frac{1}{h^{\gamma}}}}
{\displaystyle\int_{\Omega}}
\sigma\left\vert E\right\vert ^{2}dxdt\\
& =\sum\limits_{m=0,..,N-1}\mathcal{E}\left(  \zeta+mc\frac{1}{h^{\gamma}%
}\right) \\
& \leq c\frac{1}{h^{\gamma}}\sum\limits_{m=0,..,N-1}%
{\displaystyle\int_{\zeta+mc\frac{1}{h^{\gamma}}}^{\zeta+\left(  m+1\right)
c\frac{1}{h^{\gamma}}}}
\left(
{\displaystyle\int_{\Omega}}
\left(  \sigma\left\vert E\right\vert ^{2}+\left\vert \nabla p\right\vert
^{2}\right)  dx+%
{\displaystyle\int_{\omega}}
\left\vert A\right\vert ^{2}dx\right)  dt+Nch\mathcal{G}\left(  \partial
_{t}U,\zeta\right)  \text{ ,}%
\end{array}
\tag{4.25}\label{4.25}%
\end{equation}
for any $N>1$. We choose $N\in\left(  c\frac{1}{h^{\gamma}},1+c\frac
{1}{h^{\gamma}}\right]  $. Therefore, there exist $c,\gamma>0$ such that for
any $h\in\left(  0,h_{o}\right]  $,
\begin{equation}
\mathcal{E}\left(  \zeta\right)  \leq c\int_{\zeta}^{\zeta+c\left(  \frac
{1}{h}\right)  ^{\gamma}}\left(  \int_{\Omega}\left(  \sigma\left\vert
E\right\vert ^{2}+\left\vert \nabla p\right\vert ^{2}\right)  dx+\int_{\omega
}\left\vert A\right\vert ^{2}dx\right)  dt+ch\mathcal{G}\left(  \partial
_{t}U,\zeta\right)  \text{ .} \tag{4.26}\label{4.26}%
\end{equation}
On the other hand, since $\left(  U,\partial_{t}U\right)  \left(  \cdot
,\zeta\right)  =\left(  A,\partial_{t}A\right)  \left(  \cdot,\zeta\right)  $,%
\begin{equation}%
\begin{array}
[c]{ll}%
\mathcal{G}\left(  \partial_{t}U,\zeta\right)  & =\left\Vert
\operatorname{curl}E\left(  \cdot,\zeta\right)  \right\Vert _{L^{2}\left(
\Omega\right)  ^{3}}^{2}+\left\Vert \mu_{o}\operatorname{curl}H\left(
\cdot,\zeta\right)  \right\Vert _{L^{2}\left(  \Omega\right)  ^{3}}^{2}\\
& =\mu_{o}^{2}\left\Vert \partial_{t}H\left(  \cdot,\zeta\right)  \right\Vert
_{L^{2}\left(  \Omega\right)  ^{3}}^{2}+\left\Vert \left(  \partial_{t}%
E+\mu_{o}\sigma E\right)  \left(  \cdot,\zeta\right)  \right\Vert
_{L^{2}\left(  \Omega\right)  ^{3}}^{2}\\
& \leq c\left(  \mathcal{E}_{1}\left(  0\right)  +\mathcal{E}\left(  0\right)
\right)  \leq c\left\Vert \left(  E_{o},H_{o}\right)  \right\Vert _{D\left(
\mathcal{M}\right)  }^{2}%
\end{array}
\tag{4.27}\label{4.27}%
\end{equation}
where $\mathcal{M}$ is the m-accretive operator in $\mathcal{V}$ with domain
$D\left(  \mathcal{M}\right)  =\mathcal{W}$, defined as follows.%
\begin{equation}%
\begin{array}
[c]{l}%
\left\Vert \left(  F,G\right)  \right\Vert _{\mathcal{V}}^{2}=\varepsilon
_{o}\left\Vert F\right\Vert _{L^{2}\left(  \Omega\right)  ^{3}}^{2}+\mu
_{o}\left\Vert G\right\Vert _{L^{2}\left(  \Omega\right)  ^{3}}^{2}\text{ ,}\\
\mathcal{M}=\left(
\begin{array}
[c]{cc}%
\frac{1}{\varepsilon_{o}}\sigma & -\frac{1}{\varepsilon_{o}}%
\operatorname{curl}\\
\frac{1}{\mu_{o}}\operatorname{curl} & 0
\end{array}
\right)  \text{ .}%
\end{array}
\tag{4.28}\label{4.28}%
\end{equation}
Therefore, combining (\ref{4.26}) and (\ref{4.27}), we get the existence of
constants $c,\gamma>0$ such that for any $\zeta\geq0$ and $h\in\left(
0,h_{o}\right]  $,
\begin{equation}
\mathcal{E}\left(  \zeta\right)  \leq c\int_{\zeta}^{\zeta+c\left(  \frac
{1}{h}\right)  ^{\gamma}}\left(  \int_{\Omega}\left(  \sigma\left\vert
E\right\vert ^{2}+\left\vert \nabla p\right\vert ^{2}\right)  dx+\int_{\omega
}\left\vert A\right\vert ^{2}dx\right)  dt+ch\left\Vert \left(  E_{o}%
,H_{o}\right)  \right\Vert _{D\left(  \mathcal{M}\right)  }^{2}\text{ .}
\tag{4.29}\label{4.29}%
\end{equation}

\bigskip

Denote $\left(  \mathcal{T}\left(  t\right)  \right)  _{t\geq0}$\ the unique
semigroup of contractions generated by $-\mathcal{M}$. First, suppose that
$\left(  E_{o},H_{o}\right)  \in D\left(  \mathcal{M}^{3}\right)  $ and let us
define the functional of energy%
\begin{equation}
\mathcal{E}_{2}\left(  t\right)  =\frac{1}{2}\int_{\Omega}\left(
\varepsilon_{o}\left\vert \partial_{t}^{2}E\left(  x,t\right)  \right\vert
^{2}+\mu_{o}\left\vert \partial_{t}^{2}H\left(  x,t\right)  \right\vert
^{2}\right)  dx \tag{4.30}\label{4.30}%
\end{equation}
which satisfies%
\begin{equation}
\mathcal{E}_{2}\left(  t_{2}\right)  -\mathcal{E}_{2}\left(  t_{1}\right)
+\int_{t_{1}}^{t_{2}}\int_{\Omega}\sigma\left(  x\right)  \left\vert
\partial_{t}^{2}E\left(  x,t\right)  \right\vert ^{2}dxdt=0\text{ .}
\tag{4.31}\label{4.31}%
\end{equation}
Let $X_{o}=-\mathcal{M}^{2}\left(  E_{o},H_{o}\right)  $, then $\left(
\mathcal{T}\left(  t\right)  \right)  _{t\geq0}X_{o}=\left(  \partial_{t}%
^{2}E,\partial_{t}^{2}H\right)  $, $\left\Vert \mathcal{T}\left(
\zeta\right)  X_{o}\right\Vert _{\mathcal{V}}^{2}=2\mathcal{E}_{2}\left(
\zeta\right)  $ and $\left\Vert X_{o}\right\Vert _{D\left(  \mathcal{M}%
\right)  }^{2}\leq c\left\Vert \left(  E_{o},H_{o}\right)  \right\Vert
_{D\left(  \mathcal{M}^{3}\right)  }^{2}$. Further, by uniqueness of the
orthogonal decomposition in (\ref{2.1.1}) of Proposition 2.1, (\ref{4.29})
implies that for any $\left(  E_{o},H_{o}\right)  \in D\left(  \mathcal{M}%
^{3}\right)  $%
\begin{equation}
\mathcal{E}_{2}\left(  \zeta\right)  \leq c\int_{\zeta}^{\zeta+c\left(
\frac{1}{h}\right)  ^{\gamma}}\left(  \int_{\Omega}\left(  \sigma\left\vert
\partial_{t}^{2}E\right\vert ^{2}+\left\vert \partial_{t}^{2}\nabla
p\right\vert ^{2}\right)  dx+\int_{\omega}\left\vert \partial_{t}%
^{2}A\right\vert ^{2}dx\right)  dt+ch\left\Vert \left(  E_{o},H_{o}\right)
\right\Vert _{D\left(  \mathcal{M}^{3}\right)  }^{2}\text{ .} \tag{4.32}%
\label{4.32}%
\end{equation}
Since by Proposition 2.2, $\mathcal{E}\left(  \zeta\right)  \leq
c\mathcal{E}_{1}\left(  \zeta\right)  $ and in a similar way $\mathcal{E}%
_{1}\left(  \zeta\right)  \leq c\mathcal{E}_{2}\left(  \zeta\right)  $ for
some $c>0$, taking account of the first line of (\ref{2.2.1}) and
(\ref{2.2.4}), (\ref{4.32}) becomes
\begin{equation}
\mathcal{E}\left(  \zeta\right)  +\mathcal{E}_{1}\left(  \zeta\right)
+\mathcal{E}_{2}\left(  \zeta\right)  \leq c\int_{\zeta}^{\zeta+c\left(
\frac{1}{h}\right)  ^{\gamma}}\int_{\Omega}\left(  \sigma\left\vert
E\right\vert ^{2}+\sigma\left\vert \partial_{t}E\right\vert ^{2}%
+\sigma\left\vert \partial_{t}^{2}E\right\vert ^{2}\right)  dxdt+ch\left\Vert
\left(  E_{o},H_{o}\right)  \right\Vert _{D\left(  \mathcal{M}^{3}\right)
}^{2}\text{ .} \tag{4.33}\label{4.33}%
\end{equation}
Denote
\begin{equation}
\mathcal{H}\left(  \zeta\right)  =\frac{\mathcal{E}_{2}\left(  \zeta\right)
+\mathcal{E}_{1}\left(  \zeta\right)  +\mathcal{E}\left(  \zeta\right)
}{\left\Vert \left(  E_{o},H_{o}\right)  \right\Vert _{D\left(  \mathcal{M}%
^{3}\right)  }^{2}}\text{ .} \tag{4.34}\label{4.34}%
\end{equation}
Since $\mathcal{H}\leq1$, the inequality (\ref{4.33}) holds for any $h>0$.
Taking $h=c_{0}\mathcal{H}\left(  \zeta\right)  $ with some suitable small
constant $c_{0}$, we get the existence of constants $c,\gamma>0$ such that for
any $\zeta\geq0$,
\begin{equation}
\mathcal{E}\left(  \zeta\right)  +\mathcal{E}_{1}\left(  \zeta\right)
+\mathcal{E}_{2}\left(  \zeta\right)  \leq c\int_{\zeta}^{\zeta+c\left(
\frac{1}{\mathcal{H}\left(  \zeta\right)  }\right)  ^{\gamma}}\int_{\Omega
}\left(  \sigma\left\vert E\right\vert ^{2}+\sigma\left\vert \partial
_{t}E\right\vert ^{2}+\sigma\left\vert \partial_{t}^{2}E\right\vert
^{2}\right)  dxdt\text{ .} \tag{4.35}\label{4.35}%
\end{equation}
The function $\mathcal{H}$ is a continuous positive decreasing real function
on $\left[  0,+\infty\right)  $, bounded by one and satisfying from
(\ref{2.1.6}), (\ref{2.1.7}), (\ref{4.31}) and (\ref{4.35}),
\begin{equation}
\mathcal{H}\left(  \zeta\right)  \leq c\left(  \mathcal{H}\left(
\zeta\right)  -\mathcal{H}\left(  \left(  \frac{c}{\mathcal{H}\left(
\zeta\right)  }\right)  ^{\gamma}+\zeta\right)  \right)  \quad\forall\zeta
\geq0\text{ .} \tag{4.36}\label{4.36}%
\end{equation}
From \cite[p.122, Lemma B]{P2}, we deduce that there exist $C,\gamma>0$ such
that for any $t>0$%
\begin{equation}
\mathcal{E}\left(  t\right)  +\mathcal{E}_{1}\left(  t\right)  +\mathcal{E}%
_{2}\left(  \zeta\right)  \leq\frac{C}{t^{\gamma}}\left\Vert \left(
E_{o},H_{o}\right)  \right\Vert _{D\left(  \mathcal{M}^{3}\right)  }^{2}
\tag{4.37}\label{4.37}%
\end{equation}
that is
\begin{equation}
\left\Vert \mathcal{T}\left(  t\right)  Y_{o}\right\Vert _{D\left(
\mathcal{M}^{2}\right)  }^{2}\leq\frac{C}{t^{\gamma}}\left\Vert Y_{o}%
\right\Vert _{D\left(  \mathcal{M}^{3}\right)  }^{2}\quad\forall Y_{o}\in
D\left(  \mathcal{M}^{3}\right)  \text{ .} \tag{4.38}\label{4.38}%
\end{equation}

\bigskip

Since $\mathcal{M}$ is an m-accretive operator in $\mathcal{V}$ with dense
domain, one can restrict it to $D\left(  \mathcal{M}^{2}\right)  $ in a way
that its restriction operator is m-accretive. Thus the following two
properties holds.%
\begin{equation}
\forall Z_{o}\in D\left(  \mathcal{M}^{2}\right)  \quad\exists!Y_{o}\in
D\left(  \mathcal{M}^{3}\right)  \quad Y_{o}+\mathcal{M}Y_{o}=Z_{o}\text{ ,}
\tag{4.39}\label{4.39}%
\end{equation}%
\begin{equation}
\left\Vert Y_{o}\right\Vert _{D\left(  \mathcal{M}^{2}\right)  }\leq\left\Vert
Y_{o}+\mathcal{M}Y_{o}\right\Vert _{D\left(  \mathcal{M}^{2}\right)  }%
\quad\forall Y_{o}\in D\left(  \mathcal{M}^{3}\right)  \text{ .}
\tag{4.40}\label{4.40}%
\end{equation}
Consequently,
\begin{equation}%
\begin{array}
[c]{ll}%
\left\Vert \mathcal{T}\left(  t\right)  Z_{o}\right\Vert _{D\left(
\mathcal{M}\right)  }^{2} & =\left\Vert \mathcal{T}\left(  t\right)  \left(
Y_{o}+\mathcal{M}Y_{o}\right)  \right\Vert _{D\left(  \mathcal{M}\right)
}^{2}\quad\text{by (\ref{4.39})}\\
& \leq C_{1}\left\Vert \mathcal{T}\left(  t\right)  Y_{o}\right\Vert
_{D\left(  \mathcal{M}^{2}\right)  }^{2}\\
& \leq\frac{C_{2}}{t^{\gamma}}\left\Vert Y_{o}\right\Vert _{D\left(
\mathcal{M}^{3}\right)  }^{2}\quad\text{by (\ref{4.38})}\\
& \leq\frac{C_{3}}{t^{\gamma}}\left(  \left\Vert Y_{o}\right\Vert _{D\left(
\mathcal{M}^{2}\right)  }^{2}+\left\Vert Y_{o}+\mathcal{M}Y_{o}\right\Vert
_{D\left(  \mathcal{M}^{2}\right)  }^{2}\right) \\
& \leq\frac{C_{4}}{t^{\gamma}}\left\Vert Y_{o}+\mathcal{M}Y_{o}\right\Vert
_{D\left(  \mathcal{M}^{2}\right)  }^{2}\quad\text{by (\ref{4.40})}\\
& \leq\frac{C_{5}}{t^{\gamma}}\left\Vert Z_{o}\right\Vert _{D\left(
\mathcal{M}^{2}\right)  }^{2}\quad\text{by (\ref{4.39})}%
\end{array}
\tag{4.41}\label{4.41}%
\end{equation}
for suitable positive constants $C_{1},C_{2},C_{3},C_{4},C_{5}>0$.

\bigskip

Now, suppose that $\left(  E_{o},H_{o}\right)  \in D\left(  \mathcal{M}%
\right)  $. Since $\mathcal{M}$ is an m-accretive operator in $\mathcal{V}$
with dense domain, one can restrict it to $D\left(  \mathcal{M}\right)  $ in a
way that its restriction operator is m-accretive. Thus the following two
properties holds.%
\begin{equation}
\forall\left(  E_{o},H_{o}\right)  \in D\left(  \mathcal{M}\right)
\quad\exists!Z_{o}\in D\left(  \mathcal{M}^{2}\right)  \quad Z_{o}%
+\mathcal{M}Z_{o}=\left(  E_{o},H_{o}\right)  \text{ ,} \tag{4.42}\label{4.42}%
\end{equation}%
\begin{equation}
\left\Vert Z_{o}\right\Vert _{D\left(  \mathcal{M}\right)  }\leq\left\Vert
Z_{o}+\mathcal{M}Z_{o}\right\Vert _{D\left(  \mathcal{M}\right)  }\quad\forall
Z_{o}\in D\left(  \mathcal{M}^{2}\right)  \text{ .} \tag{4.43}\label{4.43}%
\end{equation}
We conclude that%
\begin{equation}%
\begin{array}
[c]{ll}%
\mathcal{E}\left(  t\right)  =\left\Vert \mathcal{T}\left(  t\right)  \left(
E_{o},H_{o}\right)  \right\Vert _{\mathcal{V}}^{2} & =\left\Vert
\mathcal{T}\left(  t\right)  \left(  Z_{o}+\mathcal{M}Z_{o}\right)
\right\Vert _{\mathcal{V}}^{2}\quad\text{by (\ref{4.42})}\\
& \leq C_{1}\left\Vert \mathcal{T}\left(  t\right)  Z_{o}\right\Vert
_{D\left(  \mathcal{M}\right)  }^{2}\\
& \leq\frac{C_{2}}{t^{\gamma}}\left\Vert Z_{o}\right\Vert _{D\left(
\mathcal{M}^{2}\right)  }^{2}\quad\text{by (\ref{4.41})}\\
& \leq\frac{C_{3}}{t^{\gamma}}\left(  \left\Vert Z_{o}\right\Vert _{D\left(
\mathcal{M}\right)  }^{2}+\left\Vert Z_{o}+\mathcal{M}Z_{o}\right\Vert
_{D\left(  \mathcal{M}\right)  }^{2}\right) \\
& \leq\frac{C_{4}}{t^{\gamma}}\left\Vert Z_{o}+\mathcal{M}Z_{o}\right\Vert
_{D\left(  \mathcal{M}\right)  }^{2}\quad\text{by (\ref{4.43})}\\
& \leq\frac{C_{5}}{t^{\gamma}}\left\Vert \left(  E_{o},H_{o}\right)
\right\Vert _{D\left(  \mathcal{M}\right)  }^{2}\quad\text{by (\ref{4.42})}\\
& \leq\frac{C_{6}}{t^{\gamma}}\left(  \mathcal{E}\left(  0\right)
+\mathcal{E}_{1}\left(  0\right)  \right)
\end{array}
\tag{4.44}\label{4.44}%
\end{equation}
for suitable positive constants $C_{1},C_{2},C_{3},C_{4},C_{5},C_{6}>0$.

\bigskip

\bigskip

\section{Proof of Proposition 4.1}

\bigskip

Recall that the definition of $\omega_{o}$ and the solution $U$ were given in
Section 4.

\bigskip

Notice that the hypothesis saying that $\Upsilon\subset\left(  \mathbb{R}%
^{2}\left\backslash D\left(  m_{1},m_{2}\right)  \right.  \right)
\times\mathbb{R}$ implies that $C_{0}^{\infty}\left(  B\left(  x_{o}%
,r_{o}/2\right)  \right)  \subset C_{0}^{\infty}\left(  \Omega\right)  $ for
any $x_{o}\in\overline{\omega_{o}}$, where $B\left(  x_{o},r\right)  $ denotes
the ball of center $x_{o}$ and radius $r$.

\bigskip

Let $\ell\in C^{\infty}\left(  \mathbb{R}^{3}\right)  $ be such that
$0\leq\ell\left(  x\right)  \leq1$, $\ell=1$ in $\mathbb{R}^{3}\left\backslash
\Theta\right.  $, $\ell\left(  x\right)  \geq\ell_{o}>0$ for any
$x\in\overline{\omega_{o}}$, $\ell=\partial_{\nu}\ell=0$ on $\Upsilon$ and
both $\nabla\ell$ and $\Delta\ell$ have support in $\Theta$.

\bigskip

The proof of Proposition 4.1 comes from the following result.

\bigskip

\textbf{Proposition 5.1}\ -. \textit{There exists }$c>0$\textit{ such that for
any }$x_{o}\in\overline{\omega_{o}}$\textit{ and any }$h\in\left(  0,1\right]
$\textit{, }$L\geq1$\textit{, }$\lambda\geq1$\textit{, the solution }%
$U$\textit{ of (\ref{4.1}) satisfies}%
\begin{equation}%
\begin{array}
[c]{ll}
& \quad%
{\displaystyle\int_{\Omega\times\mathbb{R}}}
\chi\left(  x\right)  e^{-\frac{1}{2}\left(  \frac{1}{h}\left\vert
x-x_{o}\right\vert ^{2}+t^{2}\right)  }\ell\left(  x\right)  \left\vert
\partial_{t}U\left(  x,t\right)  \right\vert ^{2}dxdt\\
& \leq c\left[  \sqrt{\frac{1}{\lambda}}+\left[  \frac{1}{\sqrt{L}}+\left(
1+hL\lambda\right)  e^{-\frac{1}{ch}}\right]  \sqrt{\lambda}\left(
\frac{\lambda^{2}}{\sqrt{h}}+\frac{1}{h}\right)  \right]  \mathcal{G}\left(
\partial_{t}U,0\right) \\
& \quad+c\left[  h\left(  1+\sqrt{hL}\right)  \sqrt{\lambda}\left(
\frac{\lambda^{2}}{\sqrt{h}}+\frac{1}{h}\right)  \right]  ^{2}\sqrt{\lambda
}\left\Vert \left(  U,\partial_{t}U\right)  \right\Vert _{L^{2}\left(
\omega\times\left(  -1-T,T+1\right)  \right)  ^{6}}^{2}%
\end{array}
\tag{5.1}\label{5.1}%
\end{equation}
\textit{where }$\chi\in C_{0}^{\infty}\left(  B\left(  x_{o},r_{o}/2\right)
\right)  $\textit{, }$0\leq\chi\leq1$\textit{ and}
\begin{equation}
T=4\left[  \frac{\lambda hL}{\sqrt{2}}+\sqrt{h}L+\frac{\sqrt{2}}{\sqrt{h}%
}\right]  \text{ \textit{.}} \tag{5.2}\label{5.2}%
\end{equation}

\bigskip

\bigskip

We shall leave the proof of Proposition 5.1 till later (see Section 6). Now we
turn to prove Proposition 4.1.

\bigskip

Let $h\in\left(  0,h_{o}\right]  $ where $h_{o}=\min\left(  1,\left(
r_{o}/8\right)  ^{2}\right)  $. We begin by covering $\overline{\omega_{o}}$
with a finite collection of balls $B\left(  x_{o}^{i},2\sqrt{h}\right)  $ for
$i\in I$ with $x_{o}^{i}\in\overline{\omega_{o}}$ and where $I$ is a countable
set such that the number of elements of $I$ is $\frac{c_{o}}{h\sqrt{h}}$ for
some constant $c_{o}>0$ independent of $h$. Then, for each $x_{o}^{i}$, we
introduce $\chi_{x_{o}^{i}}\in C_{0}^{\infty}\left(  B\left(  x_{o}^{i}%
,r_{o}/2\right)  \right)  \subset C_{0}^{\infty}\left(  \Omega\right)  $ be
such that $0\leq\chi_{x_{o}^{i}}\leq1$ and $\chi_{x_{o}^{i}}=1$ on $B\left(
x_{o}^{i},r_{o}/4\right)  \supset B\left(  x_{o}^{i},2\sqrt{h}\right)  $.
Consequently, for any $T_{o}>0$,
\begin{equation}%
\begin{array}
[c]{ll}%
{\displaystyle\int_{0}^{T_{o}}}
{\displaystyle\int_{\omega_{o}}}
\left\vert \partial_{t}U\right\vert ^{2}dxdt & \leq\frac{1}{\ell_{o}}%
e^{\frac{1}{2}T_{o}^{2}}%
{\displaystyle\int_{0}^{T_{o}}}
{\displaystyle\int_{\omega_{o}}}
e^{-\frac{1}{2}t^{2}}\ell\left(  x\right)  \left\vert \partial_{t}U\left(
x,t\right)  \right\vert ^{2}dxdt\\
& \leq\frac{1}{\ell_{o}}e^{\frac{1}{2}T_{o}^{2}+2}\sum\limits_{i\in I}%
{\displaystyle\int_{0}^{T_{o}}}
{\displaystyle\int_{B\left(  x_{o}^{i},2\sqrt{h}\right)  }}
\chi_{x_{o}^{i}}\left(  x\right)  e^{-\frac{1}{2}\left(  \frac{1}{h}\left\vert
x-x_{o}^{i}\right\vert ^{2}+t^{2}\right)  }\ell\left(  x\right)  \left\vert
\partial_{t}U\left(  x,t\right)  \right\vert ^{2}dxdt\\
& \leq\frac{1}{\ell_{o}}e^{\frac{1}{2}T_{o}^{2}+2}\sum\limits_{i\in I}%
{\displaystyle\int_{\Omega\times\mathbb{R}}}
\chi_{x_{o}^{i}}\left(  x\right)  e^{-\frac{1}{2}\left(  \frac{1}{h}\left\vert
x-x_{o}^{i}\right\vert ^{2}+t^{2}\right)  }\ell\left(  x\right)  \left\vert
\partial_{t}U\left(  x,t\right)  \right\vert ^{2}dxdt\text{ .}%
\end{array}
\tag{5.3}\label{5.3}%
\end{equation}
By virtue of Proposition 5.1, $\int_{\Omega\times\mathbb{R}}\chi_{x_{o}^{i}%
}\left(  x\right)  e^{-\frac{1}{2}\left(  \frac{1}{h}\left\vert x-x_{o}%
^{i}\right\vert ^{2}+t^{2}\right)  }\ell\left(  x\right)  \left\vert
\partial_{t}U\left(  x,t\right)  \right\vert ^{2}dxdt$ is bounded
independently of $x_{o}^{i}$ and it implies that for some constant $c>0$,
\begin{equation}%
\begin{array}
[c]{ll}%
{\displaystyle\int_{0}^{T_{o}}}
{\displaystyle\int_{\omega_{o}}}
\left\vert \partial_{t}U\right\vert ^{2}dxdt & \leq c\frac{1}{h\sqrt{h}%
}\left[  \sqrt{\frac{1}{\lambda}}+\left[  \frac{1}{\sqrt{L}}+\left(
1+hL\lambda\right)  e^{-\frac{1}{ch}}\right]  \sqrt{\lambda}\left(
\frac{\lambda^{2}}{\sqrt{h}}+\frac{1}{h}\right)  \right]  \mathcal{G}\left(
\partial_{t}U,0\right) \\
& \quad+c\frac{1}{h\sqrt{h}}\left[  h\left(  1+\sqrt{hL}\right)  \sqrt
{\lambda}\left(  \frac{\lambda^{2}}{\sqrt{h}}+\frac{1}{h}\right)  \right]
^{2}\sqrt{\lambda}\left\Vert \left(  U,\partial_{t}U\right)  \right\Vert
_{L^{2}\left(  \omega\times\left(  -1-T,T+1\right)  \right)  ^{6}}^{2}\text{
.}%
\end{array}
\tag{5.4}\label{5.4}%
\end{equation}
First, we choose $\lambda\geq1$ be such that $\lambda=\left(  \frac{h_{o}}%
{h}\right)  ^{5}$ in order that $\frac{1}{h\sqrt{h}}\frac{1}{\sqrt{\lambda}%
}\leq Ch$, then there exist $c,\delta>0$ such that for any $h\in\left(
0,h_{o}\right]  $,
\begin{equation}%
\begin{array}
[c]{ll}%
{\displaystyle\int_{0}^{T_{o}}}
{\displaystyle\int_{\omega_{o}}}
\left\vert \partial_{t}U\right\vert ^{2}dxdt & \leq ch\mathcal{G}\left(
\partial_{t}U,0\right)  +c\frac{1}{h^{\delta}}\left[  \frac{1}{\sqrt{L}%
}+Le^{-\frac{1}{ch}}\right]  \mathcal{G}\left(  \partial_{t}U,0\right) \\
& \quad+c\frac{L}{h^{\delta}}\left\Vert \left(  U,\partial_{t}U\right)
\right\Vert _{L^{2}\left(  \omega\times\left(  -1-T,T+1\right)  \right)  ^{6}%
}^{2}\text{ .}%
\end{array}
\tag{5.5}\label{5.5}%
\end{equation}
Next, we choose $L\geq1$ be such that $L=\left(  \frac{h_{o}}{h}\right)
^{2\left(  \delta+1\right)  }$ in order that $\frac{1}{h^{\delta}}\frac
{1}{\sqrt{L}}\leq Ch$, then there exist $c,c^{\prime},\gamma>0$ such that for
any $h\in\left(  0,h_{o}\right]  $,%
\begin{equation}
T=4\left[  \frac{\lambda hL}{\sqrt{2}}+\sqrt{h}L+\frac{\sqrt{2}}{\sqrt{h}%
}\right]  \leq c^{\prime}\frac{1}{h^{\gamma}}\text{ ,} \tag{5.6}\label{5.6}%
\end{equation}
and further,
\begin{equation}
\int_{0}^{T_{o}}\int_{\omega_{o}}\left\vert \partial_{t}U\right\vert
^{2}dxdt\leq c\frac{1}{h^{\gamma}}\int_{-c\frac{1}{h^{\gamma}}}^{c\frac
{1}{h^{\gamma}}}\int_{\omega}\left(  \left\vert \partial_{t}U\right\vert
^{2}+\left\vert U\right\vert ^{2}\right)  dxdt+ch\mathcal{G}\left(
\partial_{t}U,0\right)  \text{ .} \tag{5.7}\label{5.7}%
\end{equation}
By a translation in time, Proposition 4.1 follows.

\bigskip

\bigskip

\section{Proof of Proposition 5.1}

\bigskip

Let $h\in\left(  0,1\right]  $, $x_{o}\in\overline{\omega_{o}}$, $\chi\in
C_{0}^{\infty}\left(  B\left(  x_{o},r_{o}/2\right)  \right)  $ be such that
$0\leq\chi\leq1$ and $U$ be a solution of (\ref{4.1}).

\bigskip

By integrations by parts on the time variable, we can check that
\begin{equation}%
\begin{array}
[c]{ll}
& \quad%
{\displaystyle\int_{\Omega\times\mathbb{R}}}
\chi\left(  x\right)  e^{-\frac{1}{2}\left(  \frac{1}{h}\left\vert
x-x_{o}\right\vert ^{2}+t^{2}\right)  }\ell\left(  x\right)  \left\vert
\partial_{t}U\left(  x,t\right)  \right\vert ^{2}dxdt\\
& \leq2\left\vert
{\displaystyle\int_{\Omega\times\mathbb{R}}}
\chi\left(  x\right)  e^{-\frac{1}{2}\left(  \frac{1}{h}\left\vert
x-x_{o}\right\vert ^{2}+\frac{1}{2}t^{2}\right)  }\ell\left(  x\right)
\left\vert U\left(  x,t\right)  \right\vert ^{2}dxdt\right\vert \\
& \quad+\left\vert
{\displaystyle\int_{\Omega\times\mathbb{R}}}
\chi\left(  x\right)  e^{-\frac{1}{2}\left(  \frac{1}{h}\left\vert
x-x_{o}\right\vert ^{2}+t^{2}\right)  }\partial_{t}^{2}U\left(  x,t\right)
\cdot\ell\left(  x\right)  U\left(  x,t\right)  dxdt\right\vert \text{ .}%
\end{array}
\tag{6.1}\label{6.1}%
\end{equation}

\bigskip

Let us introduce for any $\theta\in\left\{  1,2\right\}  $,%
\begin{equation}
a_{o,\theta}\left(  x,t\right)  =e^{-\frac{1}{4}\left(  \frac{1}{h}\left\vert
x-x_{o}\right\vert ^{2}+\frac{1}{\theta}t^{2}\right)  }\text{ and }%
\varphi_{\theta}\left(  x,t\right)  =\chi\left(  x\right)  a_{o,\theta}\left(
x,t\right)  \text{ .} \tag{6.2}\label{6.2}%
\end{equation}
By the Fourier inversion formula,
\begin{equation}%
\begin{array}
[c]{ll}
& \quad2\left\vert
{\displaystyle\int_{\Omega\times\mathbb{R}}}
\chi\left(  x\right)  e^{-\frac{1}{2}\left(  \frac{1}{h}\left\vert
x-x_{o}\right\vert ^{2}+\frac{1}{2}t^{2}\right)  }\ell\left(  x\right)
\left\vert U\left(  x,t\right)  \right\vert ^{2}dxdt\right\vert \\
& \quad+\left\vert
{\displaystyle\int_{\Omega\times\mathbb{R}}}
\chi\left(  x\right)  e^{-\frac{1}{2}\left(  \frac{1}{h}\left\vert
x-x_{o}\right\vert ^{2}+t^{2}\right)  }\partial_{t}^{2}U\left(  x,t\right)
\cdot\ell\left(  x\right)  U\left(  x,t\right)  dxdt\right\vert \\
& =2\left\vert
{\displaystyle\int_{\Omega\times\mathbb{R}}}
\varphi_{2}\left(  x,t\right)  a_{o,2}\left(  x,t\right)  \ell\left(
x\right)  \left\vert U\left(  x,t\right)  \right\vert ^{2}dxdt\right\vert \\
& \quad+\left\vert
{\displaystyle\int_{\Omega\times\mathbb{R}}}
\varphi_{1}\left(  x,t\right)  a_{o,1}\left(  x,t\right)  \partial_{t}%
^{2}U\left(  x,t\right)  \cdot\ell\left(  x\right)  U\left(  x,t\right)
dxdt\right\vert \\
& =2\left\vert
{\displaystyle\int_{\Omega\times\mathbb{R}}}
\left(  \frac{1}{\left(  2\pi\right)  ^{4}}\int_{\mathbb{R}^{4}}e^{i\left(
x\xi+t\tau\right)  }~\widehat{\varphi_{2}U}\left(  \xi,\tau\right)  d\xi
d\tau\right)  \cdot a_{o,2}\left(  x,t\right)  \ell\left(  x\right)  U\left(
x,t\right)  dxdt\right\vert \\
& \quad+\left\vert
{\displaystyle\int_{\Omega\times\mathbb{R}}}
\left(  \frac{1}{\left(  2\pi\right)  ^{4}}\int_{\mathbb{R}^{4}}e^{i\left(
x\xi+t\tau\right)  }~\widehat{\varphi_{1}\partial_{t}^{2}U}\left(  \xi
,\tau\right)  d\xi d\tau\right)  \cdot a_{o,1}\left(  x,t\right)  \ell\left(
x\right)  U\left(  x,t\right)  dxdt\right\vert \\
& \leq2\left\vert
{\displaystyle\int_{\Omega\times\mathbb{R}}}
\left(  \frac{1}{\left(  2\pi\right)  ^{4}}\int_{\mathbb{R}^{3}}%
\int_{\left\vert \tau\right\vert <\lambda}e^{i\left(  x\xi+t\tau\right)
}~\widehat{\varphi_{2}U}\left(  \xi,\tau\right)  d\xi d\tau\right)  \cdot
a_{o,2}\left(  x,t\right)  \ell\left(  x\right)  U\left(  x,t\right)
dxdt\right\vert \\
& \quad+\left\vert
{\displaystyle\int_{\Omega\times\mathbb{R}}}
\left(  \frac{1}{\left(  2\pi\right)  ^{4}}\int_{\mathbb{R}^{3}}%
\int_{\left\vert \tau\right\vert <\lambda}e^{i\left(  x\xi+t\tau\right)
}~\widehat{\varphi_{1}\partial_{t}^{2}U}\left(  \xi,\tau\right)  d\xi
d\tau\right)  \cdot a_{o,1}\left(  x,t\right)  \ell\left(  x\right)  U\left(
x,t\right)  dxdt\right\vert \\
& \quad+2\left\vert
{\displaystyle\int_{\Omega\times\mathbb{R}}}
\left(  \frac{1}{\left(  2\pi\right)  ^{4}}\int_{\mathbb{R}^{3}}%
\int_{\left\vert \tau\right\vert \geq\lambda}e^{i\left(  x\xi+t\tau\right)
}~\widehat{\varphi_{2}U}\left(  \xi,\tau\right)  d\xi d\tau\right)  \cdot
a_{o,2}\left(  x,t\right)  \ell\left(  x\right)  U\left(  x,t\right)
dxdt\right\vert \\
& \quad+\left\vert
{\displaystyle\int_{\Omega\times\mathbb{R}}}
\left(  \frac{1}{\left(  2\pi\right)  ^{4}}\int_{\mathbb{R}^{3}}%
\int_{\left\vert \tau\right\vert \geq\lambda}e^{i\left(  x\xi+t\tau\right)
}~\widehat{\varphi_{1}\partial_{t}^{2}U}\left(  \xi,\tau\right)  d\xi
d\tau\right)  \cdot a_{o,1}\left(  x,t\right)  \ell\left(  x\right)  U\left(
x,t\right)  dxdt\right\vert
\end{array}
\tag{6.3}\label{6.3}%
\end{equation}
for any $\lambda\geq1$. Here we recall that
\begin{equation}
\widehat{F}\left(  \xi,\tau\right)  =\int_{\mathbb{R}^{4}}e^{-i\left(
x\xi+t\tau\right)  }~F\left(  x,t\right)  dxdt\text{\quad and\quad}F\left(
x,t\right)  =\frac{1}{\left(  2\pi\right)  ^{4}}\int_{\mathbb{R}^{4}%
}e^{i\left(  x\xi+t\tau\right)  }~\widehat{F}\left(  \xi,\tau\right)  d\xi
d\tau\tag{6.4}\label{6.4}%
\end{equation}
when $F$ and $\widehat{F}$ belong to $L^{1}\left(  \mathbb{R}^{4}\right)
^{3}$. On the other hand, from (A1) of Appendix A, we have that%
\begin{equation}%
\begin{array}
[c]{ll}
& \quad2\left\vert
{\displaystyle\int_{\Omega\times\mathbb{R}}}
\left(  \frac{1}{\left(  2\pi\right)  ^{4}}%
{\displaystyle\int_{\mathbb{R}^{3}}}
{\displaystyle\int_{\left\vert \tau\right\vert \geq\lambda}}
e^{i\left(  x\xi+t\tau\right)  }~\widehat{\varphi_{2}U}\left(  \xi
,\tau\right)  d\xi d\tau\right)  \cdot a_{o,2}\left(  x,t\right)  \ell\left(
x\right)  U\left(  x,t\right)  dxdt\right\vert \\
& \quad+\left\vert
{\displaystyle\int_{\Omega\times\mathbb{R}}}
\left(  \frac{1}{\left(  2\pi\right)  ^{4}}%
{\displaystyle\int_{\mathbb{R}^{3}}}
{\displaystyle\int_{\left\vert \tau\right\vert \geq\lambda}}
e^{i\left(  x\xi+t\tau\right)  }~\widehat{\varphi_{1}\partial_{t}^{2}U}\left(
\xi,\tau\right)  d\xi d\tau\right)  \cdot a_{o,1}\left(  x,t\right)
\ell\left(  x\right)  U\left(  x,t\right)  dxdt\right\vert \\
& \leq c\sqrt{\frac{1}{\lambda}}\mathcal{G}\left(  \partial_{t}U,0\right)
\text{ .}%
\end{array}
\tag{6.5}\label{6.5}%
\end{equation}
It remains to study the following two quantities
\begin{equation}%
{\displaystyle\int_{\Omega\times\mathbb{R}}}
\left(  \frac{1}{\left(  2\pi\right)  ^{4}}\int_{\mathbb{R}^{3}}%
\int_{\left\vert \tau\right\vert <\lambda}e^{i\left(  x\xi+t\tau\right)
}~\widehat{\varphi_{2}U}\left(  \xi,\tau\right)  d\xi d\tau\right)  \cdot
a_{o,2}\left(  x,t\right)  \ell\left(  x\right)  U\left(  x,t\right)  dxdt
\tag{6.6}\label{6.6}%
\end{equation}
and
\begin{equation}
\int_{\Omega\times\mathbb{R}}\left(  \frac{1}{\left(  2\pi\right)  ^{4}}%
\int_{\mathbb{R}^{3}}\int_{\left\vert \tau\right\vert <\lambda}e^{i\left(
x\xi+t\tau\right)  }~\widehat{\varphi_{1}\partial_{t}^{2}U}\left(  \xi
,\tau\right)  d\xi d\tau\right)  \cdot a_{o,1}\left(  x,t\right)  \ell\left(
x\right)  U\left(  x,t\right)  dxdt\text{ .} \tag{6.7}\label{6.7}%
\end{equation}
We claim that
\begin{equation}%
\begin{array}
[c]{ll}
& \quad2\left\vert
{\displaystyle\int_{\Omega\times\mathbb{R}}}
\left(  \frac{1}{\left(  2\pi\right)  ^{4}}%
{\displaystyle\int_{\mathbb{R}^{3}}}
{\displaystyle\int_{\left\vert \tau\right\vert <\lambda}}
e^{i\left(  x\xi+t\tau\right)  }~\widehat{\varphi_{2}U}\left(  \xi
,\tau\right)  d\xi d\tau\right)  \cdot a_{o,2}\left(  x,t\right)  \ell\left(
x\right)  U\left(  x,t\right)  dxdt\right\vert \\
& \quad+\left\vert
{\displaystyle\int_{\Omega\times\mathbb{R}}}
\left(  \frac{1}{\left(  2\pi\right)  ^{4}}%
{\displaystyle\int_{\mathbb{R}^{3}}}
{\displaystyle\int_{\left\vert \tau\right\vert <\lambda}}
e^{i\left(  x\xi+t\tau\right)  }~\widehat{\varphi_{1}\partial_{t}^{2}U}\left(
\xi,\tau\right)  d\xi d\tau\right)  \cdot a_{o,1}\left(  x,t\right)
\ell\left(  x\right)  U\left(  x,t\right)  dxdt\right\vert \\
& \leq c\left[  \left(  1+hL\lambda\right)  e^{-\frac{1}{ch}}+\frac{1}%
{\sqrt{L}}\right]  \sqrt{\lambda}\left(  \frac{\lambda^{2}}{\sqrt{h}}+\frac
{1}{h}\right)  \mathcal{G}\left(  \partial_{t}U,0\right) \\
& \quad+ch\left(  1+\sqrt{hL}\right)  \left(  \left\Vert \left(
U,\partial_{t}U\right)  \right\Vert _{L^{2}\left(  \omega\times\left(
-1-T,T+1\right)  \right)  ^{6}}\right)  \sqrt{\lambda}\left(  \frac
{\lambda^{2}}{\sqrt{h}}+\frac{1}{h}\right)  \sqrt{\mathcal{G}\left(
\partial_{t}U,0\right)  }%
\end{array}
\tag{6.8}\label{6.8}%
\end{equation}
with $T$ given by (\ref{5.2}) which implies Proposition 5.1 using (\ref{6.1}),
(\ref{6.5}) and Cauchy-Schwarz inequality.

\bigskip

The proof of our claim is divided into nine subsections. In the next
subsection, we introduce suitable sequences of Fourier integral operators.
First, we add a new variable $s\in\left[  0,L\right]  $. Next, we construct a
particular solution of the equation (6.1.10) below for $\left(  x,t,s\right)
\in\mathbb{R}^{4}\times\left[  0,L\right]  $ with good properties on
$\Gamma_{1}\cup\Gamma_{2}$.

\bigskip

\bigskip

\subsection{Fourier integral operators}

\bigskip

Let $\varphi\in C_{0}^{\infty}\left(  \Omega\right)  $ and $f=f\left(
x,t\right)  \in L^{\infty}\left(  \mathbb{R};L^{2}\left(  \mathbb{R}%
^{3}\right)  \right)  $ be such that $\widehat{\varphi f}\in L^{1}\left(
\mathbb{R}^{4}\right)  $. Let $h\in\left(  0,1\right]  $, $L\geq1$,
$\lambda\geq1$ and $\left(  x_{o},\xi_{o3}\right)  \in\overline{\omega_{o}%
}\times\left(  2\mathbb{Z}+1\right)  $. Denote $x=\left(  x_{1},x_{2}%
,x_{3}\right)  $ and $x_{o}=\left(  x_{o1},x_{o2},x_{o3}\right)  $. First, let
us introduce for any $\left(  x,t,s\right)  \in\mathbb{R}^{4}\times\left[
0,L\right]  $ and $n\in\mathbb{Z}$,%
\begin{equation}%
\begin{array}
[c]{ll}
& \quad\left(  \mathcal{A}\left(  x_{o},\xi_{o3},n\right)  f\right)  \left(
x,t,s\right) \\
& =\frac{\left(  -1\right)  ^{n}}{\left(  2\pi\right)  ^{4}}%
{\displaystyle\int_{\mathbb{R}^{2}}}
{\displaystyle\int_{\xi_{o3}-1}^{\xi_{o3}+1}}
{\displaystyle\int_{\left\vert \tau\right\vert <\lambda}}
e^{i\left(  x_{1}\xi_{1}+x_{2}\xi_{2}+t\tau\right)  }~e^{i\left[  \left(
-1\right)  ^{n}x_{3}+2n\frac{\xi_{o3}}{\left\vert \xi_{o3}\right\vert }%
\rho\right]  \xi_{3}}~e^{-i\left(  \left\vert \xi\right\vert ^{2}-\tau
^{2}\right)  hs}~\widehat{\varphi f}\left(  \xi,\tau\right) \\
& \qquad a_{\theta}\left(  x_{1}-x_{o1}-2\xi_{1}hs,x_{2}-x_{o2}-2\xi
_{2}hs,\left(  -1\right)  ^{n}x_{3}+2n\frac{\xi_{o3}}{\left\vert \xi
_{o3}\right\vert }\rho-x_{o3}-2\xi_{3}hs,t+2\tau hs,s\right)  d\xi d\tau
\end{array}
\tag{6.1.1}\label{6.1.1}%
\end{equation}
where $\xi=\left(  \xi_{1},\xi_{2},\xi_{3}\right)  \in\mathbb{R}^{2}%
\times\left[  \xi_{o3}-1,\xi_{o3}+1\right]  $,
\begin{equation}
a_{\theta}\left(  x,t,s\right)  =\left(  \frac{1}{\left(  is+1\right)  ^{3/2}%
}~e^{-\frac{1}{4h}\frac{\left\vert x\right\vert ^{2}}{is+1}}\right)  \left(
\frac{\sqrt{\theta}}{\sqrt{-ihs+\theta}}~e^{-\frac{1}{4}\frac{t^{2}%
}{-ihs+\theta}}\right)  \quad\forall\theta\in\left\{  1,2\right\}  \text{ .}
\tag{6.1.2}\label{6.1.2}%
\end{equation}
Next, let us introduce for any $\left(  x,t,s\right)  \in\mathbb{R}^{4}%
\times\left[  0,L\right]  $,%
\begin{equation}
\left(  \mathbb{A}\left(  x_{o},\xi_{o3}\right)  f\right)  \left(
x,t,s\right)  =\sum\limits_{n=-2Q}^{2P+1}\left(  \mathcal{A}\left(  x_{o}%
,\xi_{o3},n\right)  f\right)  \left(  x,t,s\right)  \text{ ,} \tag{6.1.3}%
\label{6.1.3}%
\end{equation}%
\begin{equation}
\left(  \mathbb{B}\left(  x_{o},\xi_{o3}\right)  f\right)  \left(
x,t,s\right)  =\sum\limits_{n=-2Q}^{2P+1}\left(  -1\right)  ^{n}\left(
\mathcal{A}\left(  x_{o},\xi_{o3},n\right)  f\right)  \left(  x,t,s\right)
\text{ ,} \tag{6.1.4}\label{6.1.4}%
\end{equation}
where $\left(  P,Q\right)  \in\mathbb{N}^{2}$ is the first couple of integer
numbers satisfying
\begin{equation}
\left\{
\begin{array}
[c]{ll}%
P\geq\frac{1}{4\rho}\left(  \sqrt{\left(  \left\vert \xi_{o3}\right\vert
+2\right)  \left(  L^{2}+1\right)  }+2\left(  \left\vert \xi_{o3}\right\vert
+1\right)  L\right)  \text{ ,} & \\
Q\geq\frac{1}{4\rho}\left(  \sqrt{\left(  \left\vert \xi_{o3}\right\vert
+2\right)  \left(  L^{2}+1\right)  }+2\left(  \rho-r_{o}\right)  \right)
\text{ .} &
\end{array}
\right.  \tag{6.1.5}\label{6.1.5}%
\end{equation}

\bigskip

We check after a lengthy but straightforward calculation that for any $\left(
x,t,s\right)  \in\mathbb{R}^{4}\times\left[  0,L\right]  $,
\begin{equation}
\left\{
\begin{array}
[c]{cc}%
\left(  i\partial_{s}+h\left(  \Delta-\partial_{t}^{2}\right)  \right)
\left(  \mathbb{A}\left(  x_{o},\xi_{o3}\right)  f\right)  \left(
x,t,s\right)  =0\text{ ,} & \\
\left(  i\partial_{s}+h\left(  \Delta-\partial_{t}^{2}\right)  \right)
\left(  \mathbb{B}\left(  x_{o},\xi_{o3}\right)  f\right)  \left(
x,t,s\right)  =0\text{ ,} &
\end{array}
\right.  \tag{6.1.6}\label{6.1.6}%
\end{equation}
and that for any $\left(  x_{1},x_{2},t,s\right)  \in\mathbb{R}^{3}%
\times\left[  0,L\right]  $,
\begin{equation}
\left\{
\begin{array}
[c]{rl}%
\left(  \mathbb{A}\left(  x_{o},\xi_{o3}\right)  f\right)  \left(  x_{1}%
,x_{2},\frac{\xi_{o3}}{\left\vert \xi_{o3}\right\vert }\rho,t,s\right)  &
=0\text{ ,}\\
\left(  \mathbb{A}\left(  x_{o},\xi_{o3}\right)  f\right)  \left(  x_{1}%
,x_{2},-\frac{\xi_{o3}}{\left\vert \xi_{o3}\right\vert }\rho,t,s\right)  &
=\left(  \mathcal{A}\left(  x_{o},\xi_{o3},-2Q\right)  f\right)  \left(
x_{1},x_{2},-\frac{\xi_{o3}}{\left\vert \xi_{o3}\right\vert }\rho,t,s\right)
\\
& \quad+\left(  \mathcal{A}\left(  x_{o},\xi_{o3},2P+1\right)  f\right)
\left(  x_{1},x_{2},-\frac{\xi_{o3}}{\left\vert \xi_{o3}\right\vert }%
\rho,t,s\right)  \text{ ,}%
\end{array}
\right.  \tag{6.1.7}\label{6.1.7}%
\end{equation}%
\begin{equation}
\left\{
\begin{array}
[c]{rl}%
\partial_{x_{3}}\left(  \mathbb{B}\left(  x_{o},\xi_{o3}\right)  f\right)
\left(  x_{1},x_{2},\frac{\xi_{o3}}{\left\vert \xi_{o3}\right\vert }%
\rho,t,s\right)  & =0\text{ ,}\\
\partial_{x_{3}}\left(  \mathbb{B}\left(  x_{o},\xi_{o3}\right)  f\right)
\left(  x_{1},x_{2},-\frac{\xi_{o3}}{\left\vert \xi_{o3}\right\vert }%
\rho,t,s\right)  & =\partial_{x_{3}}\left(  \mathcal{A}\left(  x_{o},\xi
_{o3},-2Q\right)  f\right)  \left(  x_{1},x_{2},-\frac{\xi_{o3}}{\left\vert
\xi_{o3}\right\vert }\rho,t,s\right) \\
& \quad-\partial_{x_{3}}\left(  \mathcal{A}\left(  x_{o},\xi_{o3},2P+1\right)
f\right)  \left(  x_{1},x_{2},-\frac{\xi_{o3}}{\left\vert \xi_{o3}\right\vert
}\rho,t,s\right)  \text{ .}%
\end{array}
\right.  \tag{6.1.8}\label{6.1.8}%
\end{equation}

\bigskip

Let $f_{j}=f_{j}\left(  x,t\right)  \in L^{\infty}\left(  \mathbb{R}%
;L^{2}\left(  \mathbb{R}^{3}\right)  \right)  $ be such that $\widehat{\varphi
f_{j}}\in L^{1}\left(  \mathbb{R}^{4}\right)  $ for any $j\in\left\{
1,2,3\right\}  $. Let us introduce
\begin{equation}
F=\left(
\begin{array}
[c]{c}%
f_{1}\\
f_{2}\\
f_{3}%
\end{array}
\right)  \quad\text{and\quad}V\left(  x_{o},\xi_{o3}\right)  F=\left(
\begin{array}
[c]{c}%
\mathbb{A}\left(  x_{o},\xi_{o3}\right)  f_{1}\\
\mathbb{A}\left(  x_{o},\xi_{o3}\right)  f_{2}\\
\mathbb{B}\left(  x_{o},\xi_{o3}\right)  f_{3}%
\end{array}
\right)  \tag{6.1.9}\label{6.1.9}%
\end{equation}
then
\begin{equation}
\left(  i\partial_{s}+h\left(  \Delta-\partial_{t}^{2}\right)  \right)
\left(  V\left(  x_{o},\xi_{o3}\right)  F\right)  \left(  x,t,s\right)
=0\text{\quad}\forall\left(  x,t,s\right)  \in\mathbb{R}^{4}\times\left[
0,L\right]  \text{ .} \tag{6.1.10}\label{6.1.10}%
\end{equation}
On another hand, let $U$ be the solution of (\ref{4.1}). Denote%
\begin{equation}
U=\left(
\begin{array}
[c]{c}%
u_{1}\\
u_{2}\\
u_{3}%
\end{array}
\right)  \quad\text{then\quad}\left\{
\begin{array}
[c]{rl}%
\forall j\in\left\{  1,2,3\right\}  \quad\partial_{t}^{2}u_{j}-\Delta
u_{j}=0 & \quad\text{in}~\Omega\times\mathbb{R}\\
u_{1}=u_{2}=0 & \quad\text{on}~\left(  \Gamma_{1}\cup\Gamma_{2}\right)
\times\mathbb{R}\\
\partial_{x_{3}}u_{3}=0 & \quad\text{on}~\left(  \Gamma_{1}\cup\Gamma
_{2}\right)  \times\mathbb{R}%
\end{array}
\right.  \tag{6.1.11}\label{6.1.11}%
\end{equation}
because $\operatorname{div}U=0$ and $U\times\nu=0$. Further, by (\ref{4.3})
and (\ref{4.6}),
\begin{equation}
\exists c>0\qquad\left\Vert u_{j}\left(  \cdot,t\right)  \right\Vert
_{H^{1}\left(  \Omega\right)  }^{2}+\left\Vert \partial_{t}u_{j}\left(
\cdot,t\right)  \right\Vert _{L^{2}\left(  \Omega\right)  }^{2}\leq
c\mathcal{G}\left(  U,0\right)  \quad\forall j\in\left\{  1,2,3\right\}
\text{ .} \tag{6.1.12}\label{6.1.12}%
\end{equation}

\bigskip

By multiplying the equation (\ref{6.1.10}) by $\ell\left(  x\right)  U\left(
x,t\right)  $ and integrating by parts over $\Omega\times\left[  -T,T\right]
\times\left[  0,L\right]  $, we have that for all $\left(  x_{o},\xi
_{o3}\right)  \in\overline{\omega_{o}}\times\left(  2\mathbb{Z}+1\right)  $
and all $h\in\left(  0,1\right]  $, $L\geq1$, $T>0$,%
\begin{equation}%
\begin{array}
[c]{rl}%
0= & -i%
{\displaystyle\int_{\Omega}}
{\displaystyle\int_{-T}^{T}}
\left(  V\left(  x_{o},\xi_{o3}\right)  F\right)  \left(  \cdot,\cdot
,0\right)  \cdot\ell Udxdt\\
& +i%
{\displaystyle\int_{\Omega}}
{\displaystyle\int_{-T}^{T}}
\left(  V\left(  x_{o},\xi_{o3}\right)  F\right)  \left(  \cdot,\cdot
,L\right)  \cdot\ell Udxdt\\
& -h%
{\displaystyle\int_{\Gamma_{1}\cup\Gamma_{2}}}
{\displaystyle\int_{-T}^{T}}
\left\{  \left(
{\displaystyle\int_{0}^{L}}
\mathbb{A}\left(  x_{o},\xi_{o3}\right)  f_{1}ds\right)  \ell\partial_{\nu
}u_{1}+\left(
{\displaystyle\int_{0}^{L}}
\mathbb{A}\left(  x_{o},\xi_{o3}\right)  f_{2}ds\right)  \ell\partial_{\nu
}u_{2}\right\}  d\sigma dt\\
& +h%
{\displaystyle\int_{\Gamma_{1}\cup\Gamma_{2}}}
{\displaystyle\int_{-T}^{T}}
\left(
{\displaystyle\int_{0}^{L}}
\partial_{\nu}\left(  \mathbb{B}\left(  x_{o},\xi_{o3}\right)  f_{3}\right)
ds\right)  \ell u_{3}d\sigma dt\\
& -h%
{\displaystyle\int}
{\displaystyle\int_{\left(  \Gamma_{1}\cup\Gamma_{2}\right)  \cap\Theta}}
{\displaystyle\int_{-T}^{T}}
\left(
{\displaystyle\int_{0}^{L}}
\mathbb{B}\left(  x_{o},\xi_{o3}\right)  f_{3}ds\right)  \partial_{\nu}\ell
u_{3}d\sigma dt\\
& -h\int_{\Omega}\left[  \left(
{\displaystyle\int_{0}^{L}}
\partial_{t}\left(  V\left(  x_{o},\xi_{o3}\right)  F\right)  \left(
\cdot,t,\cdot\right)  ds\right)  \cdot\ell U\left(  \cdot,t\right)  -\left(
{\displaystyle\int_{0}^{L}}
\left(  V\left(  x_{o},\xi_{o3}\right)  F\right)  \left(  \cdot,t,\cdot
\right)  ds\right)  \cdot\ell\partial_{t}U\left(  \cdot,t\right)  \right]
_{-T}^{T}dx\\
& +h%
{\displaystyle\int_{\omega}}
{\displaystyle\int_{-T}^{T}}
\left(
{\displaystyle\int_{0}^{L}}
V\left(  x_{o},\xi_{o3}\right)  Fds\right)  \cdot\left[  2\left(  \nabla
\ell\cdot\nabla\right)  U+\Delta\ell U\right]  dxdt\\
\equiv & \mathcal{I}_{1}+\mathcal{I}_{2}+\mathcal{I}_{3}+\mathcal{I}%
_{4}+\mathcal{I}_{5}+\mathcal{I}_{6}+\mathcal{I}_{7}\text{ .}%
\end{array}
\tag{6.1.13}\label{6.1.13}%
\end{equation}
The different terms of the last equality will be estimated separately. The
quantity $\mathcal{I}_{1}$ will allow us to recover (\ref{6.6}) (resp.
(\ref{6.7})) when $\theta=2$ and $\varphi F=\varphi_{2}U$ (resp. when
$\theta=1$ and $\varphi F=\varphi_{1}\partial_{t}^{2}U$). The dispersion
property for the one dimensional Schr\"{o}dinger operator will be used for
making $\mathcal{I}_{2}$ small for large $L$. We treat $\mathcal{I}_{3}$
(resp. $\mathcal{I}_{4}$) by applying the formula (\ref{6.1.7}) (resp.
(\ref{6.1.8})). The quantity $\mathcal{I}_{5}$ and $\mathcal{I}_{7}$ will
correspond to a term localized in $\omega$. Finally, an appropriate choice of
$T$ will bound $\mathcal{I}_{6}$ and give the desired inequality (\ref{6.9.2}) below.

\bigskip

\bigskip

\subsection{Estimate for $\mathcal{I}_{1}$ (the term at $s=0$)}

\bigskip

We estimate $\mathcal{I}_{1}=-i\int_{\Omega}\int_{-T}^{T}\left(  V\left(
x_{o},\xi_{o3}\right)  F\right)  \left(  x,t,0\right)  \cdot\ell\left(
x\right)  U\left(  x,t\right)  dxdt$ as follows.

\bigskip

\textbf{Lemma 6.1 .-} \textit{There exists }$c>0$\textit{ such that for any
}$\left(  x_{o},\xi_{o3}\right)  \in\overline{\omega_{o}}\times\left(
2\mathbb{Z}+1\right)  $\textit{ and }$h\in\left(  0,1\right]  $\textit{,
}$\lambda\geq1$\textit{, }$T>0$\textit{, we have}%
\begin{equation}%
\begin{array}
[c]{ll}
& \quad\left\vert \mathcal{I}_{1}+i%
{\displaystyle\int_{\Omega\times\mathbb{R}}}
\left(  \frac{1}{\left(  2\pi\right)  ^{4}}%
{\displaystyle\int_{\mathbb{R}^{2}}}
{\displaystyle\int_{\xi_{o3}-1}^{\xi_{o3}+1}}
{\displaystyle\int_{\left\vert \tau\right\vert <\lambda}}
e^{i\left(  x\xi+t\tau\right)  }~\widehat{\varphi F}\left(  \xi,\tau\right)
d\xi d\tau\right)  \cdot a_{o,\theta}\left(  x,t\right)  \ell\left(  x\right)
U\left(  x,t\right)  dxdt\right\vert \\
& \leq c\left(  e^{-\frac{1}{ch}}+e^{-\frac{1}{8}T^{2}}\right)  \sqrt
{\mathcal{G}\left(  U,0\right)  }\left(
{\displaystyle\int_{\mathbb{R}^{2}}}
{\displaystyle\int_{\xi_{o3}-1}^{\xi_{o3}+1}}
{\displaystyle\int_{\left\vert \tau\right\vert <\lambda}}
\left\vert \widehat{\varphi F}\left(  \xi,\tau\right)  \right\vert d\xi
d\tau\right)  \text{ .}%
\end{array}
\tag{6.2.1}\label{6.2.1}%
\end{equation}

\bigskip

Proof .- We start with the third component of $V\left(  x_{o},\xi_{o3}\right)
F$. First,
\begin{equation}%
\begin{array}
[c]{ll}
& \quad\left(  \mathbb{B}\left(  x_{o},\xi_{o3}\right)  f\right)  \left(
x,t,0\right) \\
& =\sum\limits_{n=-2Q}^{2P+1}\left[  a_{\theta}\left(  x_{1}-x_{o1}%
,x_{2}-x_{o2},\left(  -1\right)  ^{n}x_{3}+2n\frac{\xi_{o3}}{\left\vert
\xi_{o3}\right\vert }\rho-x_{o3},t,0\right)  \right. \\
& \qquad\qquad\left.  \frac{1}{\left(  2\pi\right)  ^{4}}%
{\displaystyle\int_{\mathbb{R}^{2}}}
{\displaystyle\int_{\xi_{o3}-1}^{\xi_{o3}+1}}
{\displaystyle\int_{\left\vert \tau\right\vert <\lambda}}
e^{i\left(  x_{1}\xi_{1}+x_{2}\xi_{2}+t\tau\right)  }~e^{i\left[  \left(
-1\right)  ^{n}x_{3}+2n\frac{\xi_{o3}}{\left\vert \xi_{o3}\right\vert }%
\rho\right]  \xi_{3}}~\widehat{\varphi f}\left(  \xi,\tau\right)  d\xi
d\tau\right] \\
& =a_{o,\theta}\left(  x,t\right)  \frac{1}{\left(  2\pi\right)  ^{4}}%
{\displaystyle\int_{\mathbb{R}^{2}}}
{\displaystyle\int_{\xi_{o3}-1}^{\xi_{o3}+1}}
{\displaystyle\int_{\left\vert \tau\right\vert <\lambda}}
e^{i\left(  x\xi+t\tau\right)  }~\widehat{\varphi f}\left(  \xi,\tau\right)
d\xi d\tau\\
& \quad+\sum\limits_{n\in\left\{  -2Q,\cdot\cdot\cdot,2P+1\right\}
\left\backslash \left\{  0\right\}  \right.  }\left[  a_{\theta}\left(
x_{1}-x_{o1},x_{2}-x_{o2},\left(  -1\right)  ^{n}x_{3}+2n\frac{\xi_{o3}%
}{\left\vert \xi_{o3}\right\vert }\rho-x_{o3},t,0\right)  \right. \\
& \qquad\qquad\qquad\qquad\qquad\left.  \frac{1}{\left(  2\pi\right)  ^{4}}%
{\displaystyle\int_{\mathbb{R}^{2}}}
{\displaystyle\int_{\xi_{o3}-1}^{\xi_{o3}+1}}
{\displaystyle\int_{\left\vert \tau\right\vert <\lambda}}
e^{i\left(  x_{1}\xi_{1}+x_{2}\xi_{2}+t\tau\right)  }~e^{i\left[  \left(
-1\right)  ^{n}x_{3}+2n\frac{\xi_{o3}}{\left\vert \xi_{o3}\right\vert }%
\rho\right]  \xi_{3}}~\widehat{\varphi f}\left(  \xi,\tau\right)  d\xi
d\tau\right]  \text{ .}%
\end{array}
\tag{6.2.2}\label{6.2.2}%
\end{equation}
Next, we estimate the discrete sum over $\left\{  -2Q,\cdot\cdot
\cdot,2P+1\right\}  \left\backslash \left\{  0\right\}  \right.  $.
\begin{equation}%
\begin{array}
[c]{ll}
& \left\vert \sum\limits_{n\in\left\{  -2Q,\cdot\cdot\cdot,2P+1\right\}
\left\backslash \left\{  0\right\}  \right.  }\left[  a_{\theta}\left(
x_{1}-x_{o1},x_{2}-x_{o2},\left(  -1\right)  ^{n}x_{3}+2n\frac{\xi_{o3}%
}{\left\vert \xi_{o3}\right\vert }\rho-x_{o3},t,0\right)  \right.  \right. \\
& \qquad\qquad\qquad\qquad\left.  \left.  \frac{1}{\left(  2\pi\right)  ^{4}}%
{\displaystyle\int_{\mathbb{R}^{2}}}
{\displaystyle\int_{\xi_{o3}-1}^{\xi_{o3}+1}}
{\displaystyle\int_{\left\vert \tau\right\vert <\lambda}}
e^{i\left(  x_{1}\xi_{1}+x_{2}\xi_{2}+t\tau\right)  }~e^{i\left[  \left(
-1\right)  ^{n}x_{3}+2n\frac{\xi_{o3}}{\left\vert \xi_{o3}\right\vert }%
\rho\right]  \xi_{3}}~\widehat{\varphi f}\left(  \xi,\tau\right)  d\xi
d\tau\right]  \right\vert \\
& \leq\frac{1}{\left(  2\pi\right)  ^{4}}%
{\displaystyle\int_{\mathbb{R}^{2}}}
{\displaystyle\int_{\xi_{o3}-1}^{\xi_{o3}+1}}
{\displaystyle\int_{\left\vert \tau\right\vert <\lambda}}
\left\vert \widehat{\varphi f}\left(  \xi,\tau\right)  \right\vert d\xi
d\tau~e^{-\frac{t^{2}}{8}}~e^{-\frac{\left\vert \left(  x_{1}-x_{o1}%
,x_{2}-x_{o2}\right)  \right\vert ^{2}}{4h}}\sum\limits_{n\in\mathbb{Z}%
\left\backslash \left\{  0\right\}  \right.  }e^{-\frac{\left(  \left(
-1\right)  ^{n}x_{3}+2n\frac{\xi_{o3}}{\left\vert \xi_{o3}\right\vert }%
\rho-x_{o3}\right)  ^{2}}{4h}}\text{ .}%
\end{array}
\tag{6.2.3}\label{6.2.3}%
\end{equation}
Remark that for any $\left(  x_{o1},x_{o2},x_{o3}\right)  \in\overline
{\omega_{o}}$ and $\left(  x_{1},x_{2},x_{3}\right)  \in\Omega$, the following
two cases appears. If $\left(  x_{1},x_{2}\right)  \notin\overline{D\left(
m_{1}-r_{o}/2,m_{2}-r_{o}/2\right)  }$ then we have $\left(  x_{1}%
-x_{o1}\right)  ^{2}+\left(  x_{2}-x_{o2}\right)  ^{2}\geq\left(
r_{o}/2\right)  $. If $\left(  x_{1},x_{2}\right)  \in\overline{D\left(
m_{1}-r_{o}/2,m_{2}-r_{o}/2\right)  }$, then $x_{3}\in\left[  -\rho
,\rho\right]  $ and we get $2r_{o}\leq\pm\left(  \left(  -1\right)  ^{n}%
x_{3}-x_{o3}\right)  \frac{\xi_{o3}}{\left\vert \xi_{o3}\right\vert }+2\rho$.
Therefore, for some $c>0$,%
\begin{equation}%
\begin{array}
[c]{ll}
& \quad e^{-\frac{\left\vert \left(  x_{1}-x_{o1},x_{2}-x_{o2}\right)
\right\vert ^{2}}{4h}}\sum\limits_{n\in\mathbb{Z}\left\backslash \left\{
0\right\}  \right.  }e^{-\frac{\left(  \left(  -1\right)  ^{n}x_{3}%
+2n\frac{\xi_{o3}}{\left\vert \xi_{o3}\right\vert }\rho-x_{o3}\right)  ^{2}%
}{4h}}\\
& \leq e^{-\frac{r_{o}}{8h}}\sum\limits_{n\in\mathbb{Z}}e^{-\frac{\left(
\left(  -1\right)  ^{n}x_{3}+2n\frac{\xi_{o3}}{\left\vert \xi_{o3}\right\vert
}\rho-x_{o3}\right)  ^{2}}{4h}}+2\sum\limits_{n\geq0}e^{-\frac{n^{2}\rho^{2}%
}{h}}e^{-\frac{r_{o}^{2}}{h}}\\
& \leq ce^{-\frac{1}{ch}}\text{ .}%
\end{array}
\tag{6.2.4}\label{6.2.4}%
\end{equation}
Now, we deduce from (\ref{6.2.3}) and (\ref{6.2.4}) that
\begin{equation}%
\begin{array}
[c]{ll}
& \quad\left\vert -i%
{\displaystyle\int_{\Omega}}
{\displaystyle\int_{-T}^{T}}
\left(  \mathbb{B}\left(  x_{o},\xi_{o3}\right)  f\right)  \left(
x,t,0\right)  \ell\left(  x\right)  u_{3}\left(  x,t\right)  dxdt\right. \\
& \qquad\left.  +i%
{\displaystyle\int_{\Omega}}
{\displaystyle\int_{-T}^{T}}
\left(  \frac{1}{\left(  2\pi\right)  ^{4}}%
{\displaystyle\int_{\mathbb{R}^{2}}}
{\displaystyle\int_{\xi_{o3}-1}^{\xi_{o3}+1}}
{\displaystyle\int_{\left\vert \tau\right\vert <\lambda}}
e^{i\left(  x\xi+t\tau\right)  }~\widehat{\varphi f}\left(  \xi,\tau\right)
d\xi d\tau\right)  a_{o,\theta}\left(  x,t\right)  \ell\left(  x\right)
u_{3}\left(  x,t\right)  dxdt\right\vert \\
& \leq\frac{c}{\left(  2\pi\right)  ^{4}}e^{-\frac{1}{ch}}%
{\displaystyle\int_{\Omega}}
{\displaystyle\int_{-T}^{T}}
e^{-\frac{t^{2}}{8}}\left\vert \ell\left(  x\right)  u_{3}\left(  x,t\right)
\right\vert dxdt\left(
{\displaystyle\int_{\mathbb{R}^{2}}}
{\displaystyle\int_{\xi_{o3}-1}^{\xi_{o3}+1}}
{\displaystyle\int_{\left\vert \tau\right\vert <\lambda}}
\left\vert \widehat{\varphi f}\left(  \xi,\tau\right)  \right\vert d\xi
d\tau\right) \\
& \leq ce^{-\frac{1}{ch}}\sqrt{\mathcal{G}\left(  U,0\right)  }\left(
{\displaystyle\int_{\mathbb{R}^{2}}}
{\displaystyle\int_{\xi_{o3}-1}^{\xi_{o3}+1}}
{\displaystyle\int_{\left\vert \tau\right\vert <\lambda}}
\left\vert \widehat{\varphi f}\left(  \xi,\tau\right)  \right\vert d\xi
d\tau\right)
\end{array}
\tag{6.2.5}\label{6.2.5}%
\end{equation}
where in the last line we have used the fact that the solution $U$ has the
following property, from Cauchy-Schwarz inequality and (\ref{6.1.12}),%
\begin{equation}%
\begin{array}
[c]{ll}%
{\displaystyle\int_{-T}^{T}}
e^{-\frac{t^{2}}{8}}%
{\displaystyle\int_{\Omega}}
\left\vert \ell\left(  x\right)  u_{3}\left(  x,t\right)  \right\vert dxdt &
\leq c\sqrt{\left\vert \Omega\right\vert }\left(
{\displaystyle\int_{-\infty}^{\infty}}
e^{-\frac{t^{2}}{8}}dt\right)  \sqrt{\mathcal{G}\left(  U,0\right)  }\\
& \leq c\sqrt{\left\vert \Omega\right\vert }\left(  2\sqrt{2\pi}\right)
\sqrt{\mathcal{G}\left(  U,0\right)  }\text{ .}%
\end{array}
\tag{6.2.6}\label{6.2.6}%
\end{equation}
Here and hereafter, $c$ will be used to denote a generic constant, not
necessarily the same in any two places. On the other hand, we cut the integral
on time into two parts to obtain
\begin{equation}%
\begin{array}
[c]{ll}
& \quad%
{\displaystyle\int_{\Omega\times\mathbb{R}}}
\left(  \frac{1}{\left(  2\pi\right)  ^{4}}%
{\displaystyle\int_{\mathbb{R}^{2}}}
{\displaystyle\int_{\xi_{o3}-1}^{\xi_{o3}+1}}
{\displaystyle\int_{\left\vert \tau\right\vert <\lambda}}
e^{i\left(  x\xi+t\tau\right)  }~\widehat{\varphi f}\left(  \xi,\tau\right)
d\xi d\tau\right)  a_{o,\theta}\left(  x,t\right)  \ell\left(  x\right)
u_{3}\left(  x,t\right)  dxdt\\
& =%
{\displaystyle\int_{\Omega}}
{\displaystyle\int_{-T}^{T}}
\left(  \frac{1}{\left(  2\pi\right)  ^{4}}%
{\displaystyle\int_{\mathbb{R}^{2}}}
{\displaystyle\int_{\xi_{o3}-1}^{\xi_{o3}+1}}
{\displaystyle\int_{\left\vert \tau\right\vert <\lambda}}
e^{i\left(  x\xi+t\tau\right)  }~\widehat{\varphi f}\left(  \xi,\tau\right)
d\xi d\tau\right)  a_{o,\theta}\left(  x,t\right)  \ell\left(  x\right)
u_{3}\left(  x,t\right)  dxdt\\
& +%
{\displaystyle\int_{\Omega}}
{\displaystyle\int_{\mathbb{R}\left\backslash \left(  -T,T\right)  \right.  }}
\left(  \frac{1}{\left(  2\pi\right)  ^{4}}%
{\displaystyle\int_{\mathbb{R}^{2}}}
{\displaystyle\int_{\xi_{o3}-1}^{\xi_{o3}+1}}
{\displaystyle\int_{\left\vert \tau\right\vert <\lambda}}
e^{i\left(  x\xi+t\tau\right)  }~\widehat{\varphi f}\left(  \xi,\tau\right)
d\xi d\tau\right)  e^{-\frac{1}{4}\left(  \frac{1}{h}\left\vert x-x_{o}%
\right\vert ^{2}+\frac{1}{\theta}t^{2}\right)  }\ell\left(  x\right)
u_{3}\left(  x,t\right)  dxdt
\end{array}
\tag{6.2.7}\label{6.2.7}%
\end{equation}
and, by using (\ref{6.1.12}) and Cauchy-Schwarz inequality, we have%
\begin{equation}%
\begin{array}
[c]{ll}
& \quad\left\vert i%
{\displaystyle\int_{\Omega}}
{\displaystyle\int_{\mathbb{R}\left\backslash \left(  -T,T\right)  \right.  }}
\left(  \frac{1}{\left(  2\pi\right)  ^{4}}%
{\displaystyle\int_{\mathbb{R}^{2}}}
{\displaystyle\int_{\xi_{o3}-1}^{\xi_{o3}+1}}
{\displaystyle\int_{\left\vert \tau\right\vert <\lambda}}
e^{i\left(  x\xi+t\tau\right)  }~\widehat{\varphi f}\left(  \xi,\tau\right)
d\xi d\tau\right)  e^{-\frac{1}{4}\left(  \frac{1}{h}\left\vert x-x_{o}%
\right\vert ^{2}+\frac{1}{\theta}t^{2}\right)  }\ell\left(  x\right)
u_{3}\left(  x,t\right)  dxdt\right\vert \\
& \leq ce^{-\frac{1}{8}T^{2}}\sqrt{\mathcal{G}\left(  U,0\right)  }\left(
{\displaystyle\int_{\mathbb{R}^{2}}}
{\displaystyle\int_{\xi_{o3}-1}^{\xi_{o3}+1}}
{\displaystyle\int_{\left\vert \tau\right\vert <\lambda}}
\left\vert \widehat{\varphi f}\left(  \xi,\tau\right)  \right\vert d\xi
d\tau\right)  \text{ .}%
\end{array}
\tag{6.2.8}\label{6.2.8}%
\end{equation}
We conclude from (\ref{6.2.5}), (\ref{6.2.7}) and (\ref{6.2.8}) that
\begin{equation}%
\begin{array}
[c]{ll}
& \quad\left\vert -i%
{\displaystyle\int_{\Omega}}
{\displaystyle\int_{-T}^{T}}
\left(  \mathbb{B}\left(  x_{o},\xi_{o3}\right)  f\right)  \left(
x,t,0\right)  \ell\left(  x\right)  u_{3}\left(  x,t\right)  dxdt\right. \\
& \qquad\left.  +i%
{\displaystyle\int_{\Omega\times\mathbb{R}}}
\left(  \frac{1}{\left(  2\pi\right)  ^{4}}%
{\displaystyle\int_{\mathbb{R}^{2}}}
{\displaystyle\int_{\xi_{o3}-1}^{\xi_{o3}+1}}
{\displaystyle\int_{\left\vert \tau\right\vert <\lambda}}
e^{i\left(  x\xi+t\tau\right)  }~\widehat{\varphi f}\left(  \xi,\tau\right)
d\xi d\tau\right)  a_{o,\theta}\left(  x,t\right)  \ell\left(  x\right)
u_{3}\left(  x,t\right)  dxdt\right\vert \\
& \leq c\left(  e^{-\frac{1}{ch}}+e^{-\frac{1}{8}T^{2}}\right)  \sqrt
{\mathcal{G}\left(  U,0\right)  }\left(
{\displaystyle\int_{\mathbb{R}^{2}}}
{\displaystyle\int_{\xi_{o3}-1}^{\xi_{o3}+1}}
{\displaystyle\int_{\left\vert \tau\right\vert <\lambda}}
\left\vert \widehat{\varphi f}\left(  \xi,\tau\right)  \right\vert d\xi
d\tau\right)  \text{ .}%
\end{array}
\tag{6.2.9}\label{6.2.9}%
\end{equation}
Similarly,
\begin{equation}%
\begin{array}
[c]{ll}
& \quad\left\vert -i%
{\displaystyle\int_{\Omega}}
{\displaystyle\int_{-T}^{T}}
\left(  \mathbb{A}\left(  x_{o},\xi_{o3}\right)  f\right)  \left(
x,t,0\right)  \ell\left(  x\right)  u_{j}\left(  x,t\right)  dxdt\right. \\
& \qquad\left.  +i\int_{\Omega\times\mathbb{R}}\left(  \frac{1}{\left(
2\pi\right)  ^{4}}%
{\displaystyle\int_{\mathbb{R}^{2}}}
{\displaystyle\int_{\xi_{o3}-1}^{\xi_{o3}+1}}
{\displaystyle\int_{\left\vert \tau\right\vert <\lambda}}
e^{i\left(  x\xi+t\tau\right)  }~\widehat{\varphi f}\left(  \xi,\tau\right)
d\xi d\tau\right)  a_{o,\theta}\left(  x,t\right)  \ell\left(  x\right)
u_{j}\left(  x,t\right)  dxdt\right\vert \\
& \leq c\left(  e^{-\frac{1}{ch}}+e^{-\frac{1}{8}T^{2}}\right)  \sqrt
{\mathcal{G}\left(  U,0\right)  }\left(
{\displaystyle\int_{\mathbb{R}^{2}}}
{\displaystyle\int_{\xi_{o3}-1}^{\xi_{o3}+1}}
{\displaystyle\int_{\left\vert \tau\right\vert <\lambda}}
\left\vert \widehat{\varphi f}\left(  \xi,\tau\right)  \right\vert d\xi
d\tau\right)  \text{ .}%
\end{array}
\tag{6.2.10}\label{6.2.10}%
\end{equation}
This completes the proof.

\bigskip

\bigskip

\subsection{Estimate for $\mathcal{I}_{2}$ (the term at $s=L$)}

\bigskip

We estimate $\mathcal{I}_{2}=i\int_{\Omega}\int_{-T}^{T}\left(  V\left(
x_{o},\xi_{o3}\right)  F\right)  \left(  x,t,L\right)  \cdot\ell\left(
x\right)  U\left(  x,t\right)  dxdt$ as follows.

\bigskip

\textbf{Lemma 6.2 .-} \textit{There exists }$c>0$\textit{ such that for any
}$\left(  x_{o},\xi_{o3}\right)  \in\overline{\omega_{o}}\times\left(
2\mathbb{Z}+1\right)  $\textit{ and }$h\in\left(  0,1\right]  $\textit{,
}$L\geq1$\textit{, }$\lambda\geq1$\textit{, }$T>0$\textit{,} \textit{we have}%
\begin{equation}
\left\vert \mathcal{I}_{2}\right\vert \leq c\left(  \frac{1}{\sqrt{L}%
}+e^{-\frac{1}{4h}}\right)  \sqrt{\mathcal{G}\left(  U,0\right)  }\left(
\int_{\mathbb{R}^{2}}\int_{\xi_{o3}-1}^{\xi_{o3}+1}\int_{\left\vert
\tau\right\vert <\lambda}\left\vert \widehat{\varphi F}\left(  \xi
,\tau\right)  \right\vert d\xi d\tau\right)  \text{ .} \tag{6.3.1}%
\label{6.3.1}%
\end{equation}

\bigskip

Proof .- We start with the third component of $V\left(  x_{o},\xi_{o3}\right)
F$. First,%
\begin{equation}%
\begin{array}
[c]{ll}
& \quad\left\vert \left(  \mathbb{B}\left(  x_{o},\xi_{o3}\right)  f\right)
\left(  x,t,L\right)  \right\vert \\
& \leq\sum\limits_{n\in\mathbb{Z}\left\backslash \left\{  -2Q,\cdot\cdot
\cdot,2P+1\right\}  \right.  }\left\vert \left(  \mathcal{A}\left(  x_{o}%
,\xi_{o3},n\right)  f\right)  \left(  x,t,L\right)  \right\vert +\left\vert
\sum\limits_{n\in\mathbb{Z}}\left(  -1\right)  ^{n}\left(  \mathcal{A}\left(
x_{o},\xi_{o3},n\right)  f\right)  \left(  x,t,L\right)  \right\vert \text{ .}%
\end{array}
\tag{6.3.2}\label{6.3.2}%
\end{equation}
Next,
\begin{equation}%
\begin{array}
[c]{ll}
& \quad\sum\limits_{n\in\mathbb{Z}\left\backslash \left\{  -2Q,\cdot\cdot
\cdot,2P+1\right\}  \right.  }\left\vert \left(  \mathcal{A}\left(  x_{o}%
,\xi_{o3},n\right)  f\right)  \left(  x,t,L\right)  \right\vert \\
& \leq\frac{1}{\left(  2\pi\right)  ^{4}}%
{\displaystyle\int_{\mathbb{R}^{2}}}
{\displaystyle\int_{\xi_{o3}-1}^{\xi_{o3}+1}}
{\displaystyle\int_{\left\vert \tau\right\vert <\lambda}}
\left\vert \widehat{\varphi f}\left(  \xi,\tau\right)  \right\vert \left(
\frac{\sqrt{\theta}}{\left(  \sqrt{\left(  hL\right)  ^{2}+\theta^{2}}\right)
^{1/2}}~e^{-\frac{\theta}{4}\frac{\left(  t+2\tau hL\right)  ^{2}}{\left(
hL\right)  ^{2}+\theta^{2}}}\right) \\
& \qquad\left(  \frac{1}{\sqrt{L^{2}+1}}~e^{-\frac{1}{4h}\frac{\left\vert
\left(  x_{1}-x_{o1}-2\xi_{1}hL,x_{2}-x_{o2}-2\xi_{2}hL\right)  \right\vert
^{2}}{L^{2}+1}}\right) \\
& \qquad\frac{1}{\left(  \sqrt{L^{2}+1}\right)  ^{1/2}}\left(  \sum
\limits_{n\in\mathbb{Z}\left\backslash \left\{  -2Q,\cdot\cdot\cdot
,2P+1\right\}  \right.  }~e^{-\frac{1}{4h}\frac{\left(  \left(  -1\right)
^{n}x_{3}+2n\frac{\xi_{o3}}{\left\vert \xi_{o3}\right\vert }\rho-x_{o3}%
-2\xi_{3}hL\right)  ^{2}}{L^{2}+1}}\right)  d\xi d\tau\text{ .}%
\end{array}
\tag{6.3.3}\label{6.3.3}%
\end{equation}
When $\xi_{3}\in\left(  \xi_{o3}-1,\xi_{o3}+1\right)  $ with $\xi_{o3}%
\in\left(  2\mathbb{Z}+1\right)  $,
\begin{equation}%
\begin{array}
[c]{ll}%
\sqrt{L^{2}+1} & \leq4P\rho-2\left(  \left\vert \xi_{o3}\right\vert +1\right)
L\quad\text{from our choice of }P\\
& \leq4P\rho-2\left\vert \xi_{3}\right\vert hL+2\rho-\left\vert x_{3}%
\right\vert -\left\vert x_{o3}\right\vert \quad\text{because }\left\vert
x_{3}\right\vert +\left\vert x_{o3}\right\vert \leq2\left(  \rho-r_{o}\right)
\\
& \leq4P\rho-2\left\vert \xi_{3}\right\vert hL+2\rho+\frac{\xi_{o3}%
}{\left\vert \xi_{o3}\right\vert }\left[  -\left(  -1\right)  ^{n}x_{3}%
-x_{o3}\right]  \quad\forall n\in\mathbb{Z}%
\end{array}
\tag{6.3.4}\label{6.3.4}%
\end{equation}
thus%
\begin{equation}%
\begin{array}
[c]{ll}
& \sum\limits_{n\geq2P+2}e^{-\frac{1}{4h}\frac{\left(  \left(  -1\right)
^{n}x_{3}+2n\frac{\xi_{o3}}{\left\vert \xi_{o3}\right\vert }\rho-x_{o3}%
-2\xi_{3}hL\right)  ^{2}}{L^{2}+1}}=\sum\limits_{n\geq2P+2}e^{-\frac{1}%
{4h}\frac{\left(  2n\rho-2\left\vert \xi_{3}\right\vert hL+\frac{\xi_{o3}%
}{\left\vert \xi_{o3}\right\vert }\left[  \left(  -1\right)  ^{n}x_{3}%
-x_{o3}\right]  \right)  ^{2}}{L^{2}+1}}\\
& =\sum\limits_{n\geq1}e^{-\frac{1}{4h}\frac{\left(  2n\rho+4P\rho-2\left\vert
\xi_{3}\right\vert hL+2\rho+\frac{\xi_{o3}}{\left\vert \xi_{o3}\right\vert
}\left[  -\left(  -1\right)  ^{n}x_{3}-x_{o3}\right]  \right)  ^{2}}{L^{2}+1}%
}\\
& \leq\sum\limits_{n\geq1}e^{-\frac{1}{4h}\frac{\left(  2n\rho\right)  ^{2}%
}{L^{2}+1}}e^{-\frac{1}{4h}\frac{\left(  4P\rho-2\left\vert \xi_{3}\right\vert
hL+2\rho+\frac{\xi_{o3}}{\left\vert \xi_{o3}\right\vert }\left[  -\left(
-1\right)  ^{n}x_{3}-x_{o3}\right]  \right)  ^{2}}{L^{2}+1}}\\
& \leq e^{-\frac{1}{4h}}\sum\limits_{n\geq1}e^{-\frac{1}{h}\frac{\left(
n\rho\right)  ^{2}}{L^{2}+1}}\leq e^{-\frac{1}{4h}}\left(  \frac{\sqrt{\pi}%
}{2}\frac{\sqrt{h}\sqrt{L^{2}+1}}{\rho}\right)  \text{ .}%
\end{array}
\tag{6.3.5}\label{6.3.5}%
\end{equation}
Also,
\begin{equation}%
\begin{array}
[c]{ll}%
\sqrt{L^{2}+1} & \leq4Q\rho-2\left(  \rho-r_{o}\right)  \quad\text{from our
choice of }Q\\
& \leq4Q\rho+2\left\vert \xi_{3}\right\vert hL-\left\vert x_{3}\right\vert
-\left\vert x_{o3}\right\vert \quad\text{because }\left\vert x_{3}\right\vert
+\left\vert x_{o3}\right\vert \leq2\left(  \rho-r_{o}\right) \\
& \leq4Q\rho+2\left\vert \xi_{3}\right\vert hL-\frac{\xi_{o3}}{\left\vert
\xi_{o3}\right\vert }\left[  \left(  -1\right)  ^{n}x_{3}-x_{o3}\right]
\quad\forall n\in\mathbb{Z}%
\end{array}
\tag{6.3.6}\label{6.3.6}%
\end{equation}
thus
\begin{equation}%
\begin{array}
[c]{ll}
& \sum\limits_{n\leq-2Q-1}e^{-\frac{1}{4h}\frac{\left(  \left(  -1\right)
^{n}x_{3}+2n\frac{\xi_{o3}}{\left\vert \xi_{o3}\right\vert }\rho-x_{o3}%
-2\xi_{3}hL\right)  ^{2}}{L^{2}+1}}=\sum\limits_{n\geq2Q+1}e^{-\frac{1}%
{4h}\frac{\left(  2n\rho+2\left\vert \xi_{3}\right\vert hL-\frac{\xi_{o3}%
}{\left\vert \xi_{o3}\right\vert }\left[  \left(  -1\right)  ^{n}x_{3}%
-x_{o3}\right]  \right)  ^{2}}{L^{2}+1}}\\
& \leq\sum\limits_{n\geq1}e^{-\frac{1}{4h}\frac{\left(  2n\rho+4Q\rho
+2\left\vert \xi_{3}\right\vert hL-\frac{\xi_{o3}}{\left\vert \xi
_{o3}\right\vert }\left[  \left(  -1\right)  ^{n}x_{3}-x_{o3}\right]  \right)
^{2}}{L^{2}+1}}\\
& \leq\sum\limits_{n\geq1}e^{-\frac{1}{4h}\frac{\left(  2n\rho\right)  ^{2}%
}{L^{2}+1}}e^{-\frac{1}{4h}\frac{\left(  4Q\rho+2\left\vert \xi_{3}\right\vert
hL-\frac{\xi_{o3}}{\left\vert \xi_{o3}\right\vert }\left[  \left(  -1\right)
^{n}x_{3}-x_{o3}\right]  \right)  ^{2}}{L^{2}+1}}\\
& \leq e^{-\frac{1}{4h}}\sum\limits_{n\geq1}e^{-\frac{1}{h}\frac{\left(
n\rho\right)  ^{2}}{L^{2}+1}}\leq e^{-\frac{1}{4h}}\left(  \frac{\sqrt{\pi}%
}{2}\frac{\sqrt{h}\sqrt{L^{2}+1}}{\rho}\right)  \text{ .}%
\end{array}
\tag{6.3.7}\label{6.3.7}%
\end{equation}
It implies from (\ref{6.3.3}), (\ref{6.3.5}) and (\ref{6.3.7}) that%
\begin{equation}%
\begin{array}
[c]{ll}
& \quad\sum\limits_{n\in\mathbb{Z}\left\backslash \left\{  -2Q,\cdot\cdot
\cdot,2P+1\right\}  \right.  }\left\vert \left(  \mathcal{A}\left(  x_{o}%
,\xi_{o3},n\right)  f\right)  \left(  x,t,L\right)  \right\vert \\
& \leq\frac{1}{\left(  2\pi\right)  ^{4}}\frac{1}{\sqrt{L^{2}+1}}\frac
{1}{\left(  \sqrt{L^{2}+1}\right)  ^{1/2}}\left(  \frac{\sqrt{\pi}\sqrt
{h}\sqrt{L^{2}+1}}{\rho}\right)  e^{-\frac{1}{4h}}\\
& \qquad%
{\displaystyle\int_{\mathbb{R}^{2}}}
{\displaystyle\int_{\xi_{o3}-1}^{\xi_{o3}+1}}
{\displaystyle\int_{\left\vert \tau\right\vert <\lambda}}
\left\vert \widehat{\varphi f}\left(  \xi,\tau\right)  \right\vert \left(
\frac{\sqrt{\theta}}{\left(  \sqrt{\left(  hL\right)  ^{2}+\theta^{2}}\right)
^{1/2}}~e^{-\frac{\theta}{4}\frac{\left(  t+2\tau hL\right)  ^{2}}{\left(
hL\right)  ^{2}+\theta^{2}}}\right)  d\xi d\tau\text{ .}%
\end{array}
\tag{6.3.8}\label{6.3.8}%
\end{equation}
Now,
\begin{equation}%
\begin{array}
[c]{ll}
& \quad\left\vert \sum\limits_{n\in\mathbb{Z}}\left(  -1\right)  ^{n}\left(
\mathcal{A}\left(  x_{o},\xi_{o3},n\right)  f\right)  \left(  x,t,L\right)
\right\vert \\
& \leq\frac{1}{\left(  2\pi\right)  ^{4}}%
{\displaystyle\int_{\mathbb{R}^{2}}}
{\displaystyle\int_{\xi_{o3}-1}^{\xi_{o3}+1}}
{\displaystyle\int_{\left\vert \tau\right\vert <\lambda}}
\left\vert \widehat{\varphi f}\left(  \xi,\tau\right)  \right\vert \left(
\frac{\sqrt{\theta}}{\left(  \sqrt{\left(  hL\right)  ^{2}+\theta^{2}}\right)
^{1/2}}~e^{-\frac{\theta}{4}\frac{\left(  t+2\tau hL\right)  ^{2}}{\left(
hL\right)  ^{2}+\theta^{2}}}\right) \\
& \qquad\left(  \frac{1}{\sqrt{L^{2}+1}}~e^{-\frac{1}{4h}\frac{\left\vert
\left(  x_{1}-x_{o1}-2\xi_{1}hL,x_{2}-x_{o2}-2\xi_{2}hL\right)  \right\vert
^{2}}{L^{2}+1}}\right) \\
& \qquad\left\vert \sum\limits_{n\in\mathbb{Z}}\left(  \frac{1}{\sqrt{iL+1}%
}~e^{i\left[  \left(  -1\right)  ^{n}x_{3}+2n\frac{\xi_{o3}}{\left\vert
\xi_{o3}\right\vert }\rho\right]  \xi_{3}}~e^{-\frac{1}{4h}\frac{\left(
\left(  -1\right)  ^{n}x_{3}+2n\frac{\xi_{o3}}{\left\vert \xi_{o3}\right\vert
}\rho-x_{o3}-2\xi_{3}hL\right)  ^{2}}{iL+1}}\right)  \right\vert d\xi d\tau\\
& \leq\frac{c}{\sqrt{L^{2}+1}}%
{\displaystyle\int_{\mathbb{R}^{2}}}
{\displaystyle\int_{\xi_{o3}-1}^{\xi_{o3}+1}}
{\displaystyle\int_{\left\vert \tau\right\vert <\lambda}}
\left\vert \widehat{\varphi f}\left(  \xi,\tau\right)  \right\vert \left(
\frac{\sqrt{\theta}}{\left(  \sqrt{\left(  hL\right)  ^{2}+\theta^{2}}\right)
^{1/2}}~e^{-\frac{\theta}{4}\frac{\left(  t+2\tau hL\right)  ^{2}}{\left(
hL\right)  ^{2}+\theta^{2}}}\right)  d\xi d\tau
\end{array}
\tag{6.3.9}\label{6.3.9}%
\end{equation}
because from Appendix B with $z=\frac{4h}{\rho^{2}}\left(  iL+1\right)  $, we
know that
\begin{equation}
\left\vert \frac{1}{\sqrt{iL+1}}\sum\limits_{n\in\mathbb{Z}}e^{i\left[
\left(  -1\right)  ^{n}x_{3}+2n\frac{\xi_{o3}}{\left\vert \xi_{o3}\right\vert
}\rho\right]  \xi_{3}}~e^{-\frac{1}{4h}\frac{\left(  \left(  -1\right)
^{n}x_{3}+2n\frac{\xi_{o3}}{\left\vert \xi_{o3}\right\vert }\rho-x_{o3}%
-2\xi_{3}hL\right)  ^{2}}{iL+1}}\right\vert \leq\frac{2\sqrt{h}}{\rho}\left(
\frac{\sqrt{\pi}}{2}+\frac{\rho}{\sqrt{h}}\right)  \text{ .} \tag{6.3.10}%
\label{6.3.10}%
\end{equation}
Finally, (\ref{6.3.2}), (\ref{6.3.8}) and (\ref{6.3.10}) imply that
\begin{equation}%
\begin{array}
[c]{ll}
& \quad\left\vert \left(  \mathbb{B}\left(  x_{o},\xi_{o3}\right)  f\right)
\left(  x,t,L\right)  \right\vert \\
& \leq c\left(  \frac{1}{\sqrt{L^{2}+1}}+\frac{1}{\left(  \sqrt{L^{2}%
+1}\right)  ^{1/2}}e^{-\frac{1}{4h}}\right)
{\displaystyle\int_{\mathbb{R}^{2}}}
{\displaystyle\int_{\xi_{o3}-1}^{\xi_{o3}+1}}
{\displaystyle\int_{\left\vert \tau\right\vert <\lambda}}
\left\vert \widehat{\varphi f}\left(  \xi,\tau\right)  \right\vert \left(
\frac{\sqrt{\theta}}{\left(  \sqrt{\left(  hL\right)  ^{2}+\theta^{2}}\right)
^{1/2}}~e^{-\frac{\theta}{4}\frac{\left(  t+2\tau hL\right)  ^{2}}{\left(
hL\right)  ^{2}+\theta^{2}}}\right)  d\xi d\tau
\end{array}
\tag{6.3.11}\label{6.3.11}%
\end{equation}
and we conclude that
\begin{equation}%
\begin{array}
[c]{ll}
& \quad\left\vert i%
{\displaystyle\int_{\Omega}}
{\displaystyle\int_{-T}^{T}}
\left(  \mathbb{B}\left(  x_{o},\xi_{o3}\right)  f\right)  \left(
x,t,L\right)  \ell\left(  x\right)  u_{3}\left(  x,t\right)  dxdt\right\vert
\\
& \leq c\left(  \frac{1}{\sqrt{L^{2}+1}}+\frac{1}{\left(  \sqrt{L^{2}%
+1}\right)  ^{1/2}}e^{-\frac{1}{4h}}\right) \\
& \quad%
{\displaystyle\int_{\mathbb{R}^{2}}}
{\displaystyle\int_{\xi_{o3}-1}^{\xi_{o3}+1}}
{\displaystyle\int_{\left\vert \tau\right\vert <\lambda}}
\left\vert \widehat{\varphi f}\left(  \xi,\tau\right)  \right\vert d\xi%
{\displaystyle\int_{-T}^{T}}
\left(  \frac{\sqrt{\theta}}{\left(  \sqrt{\left(  hL\right)  ^{2}+\theta^{2}%
}\right)  ^{1/2}}~e^{-\frac{\theta}{4}\frac{\left(  t+2\tau hL\right)  ^{2}%
}{\left(  hL\right)  ^{2}+\theta^{2}}}\right)
{\displaystyle\int_{\Omega}}
\left\vert \ell\left(  x\right)  u_{3}\left(  x,t\right)  \right\vert
dxdtd\tau\\
& \leq c\left(  \frac{1}{\sqrt{L}}+e^{-\frac{1}{4h}}\right)  \sqrt
{\mathcal{G}\left(  U,0\right)  }\left(
{\displaystyle\int_{\mathbb{R}^{2}}}
{\displaystyle\int_{\xi_{o3}-1}^{\xi_{o3}+1}}
{\displaystyle\int_{\left\vert \tau\right\vert <\lambda}}
\left\vert \widehat{\varphi f}\left(  \xi,\tau\right)  \right\vert d\xi
d\tau\right)  \text{ }%
\end{array}
\tag{6.3.12}\label{6.3.12}%
\end{equation}
where in the last line we have used the fact that the solution $U$ has the
following property, from Cauchy-Schwarz inequality and (\ref{6.1.12}),
\begin{equation}%
\begin{array}
[c]{ll}
& \quad%
{\displaystyle\int_{-T}^{T}}
e^{-\frac{\theta}{4}\frac{\left(  t+2\tau hL\right)  ^{2}}{\left(  hL\right)
^{2}+\theta^{2}}}\int_{\Omega}\left\vert \ell\left(  x\right)  u_{3}\left(
x,t\right)  \right\vert dxdt\\
& \leq c\sqrt{\left\vert \Omega\right\vert }\left(
{\displaystyle\int_{-\infty}^{\infty}}
e^{-\frac{\theta}{4}\frac{t^{2}}{\left(  hL\right)  ^{2}+\theta^{2}}%
}dt\right)  \sqrt{\mathcal{G}\left(  U,0\right)  }\\
& \leq c\sqrt{\left\vert \Omega\right\vert }\left(  \frac{2\sqrt{\pi}}%
{\sqrt{\theta}}\sqrt{\left(  hL\right)  ^{2}+\theta^{2}}\right)
\sqrt{\mathcal{G}\left(  U,0\right)  }\text{ .}%
\end{array}
\tag{6.3.13}\label{6.3.13}%
\end{equation}
Similarly,
\begin{equation}%
\begin{array}
[c]{ll}
& \quad\left\vert i%
{\displaystyle\int_{\Omega}}
{\displaystyle\int_{-T}^{T}}
\left(  \mathbb{A}\left(  x_{o},\xi_{o3}\right)  f\right)  \left(
x,t,L\right)  \ell\left(  x\right)  u_{j}\left(  x,t\right)  dxdt\right\vert
\\
& \leq c\left(  \frac{1}{\sqrt{L}}+e^{-\frac{1}{4h}}\right)  \sqrt
{\mathcal{G}\left(  U,0\right)  }\left(
{\displaystyle\int_{\mathbb{R}^{2}}}
{\displaystyle\int_{\xi_{o3}-1}^{\xi_{o3}+1}}
{\displaystyle\int_{\left\vert \tau\right\vert <\lambda}}
\left\vert \widehat{\varphi f}\left(  \xi,\tau\right)  \right\vert d\xi
d\tau\right)  \text{ ,}%
\end{array}
\tag{6.3.14}\label{6.3.14}%
\end{equation}
using the estimate%
\begin{equation}
\left\vert \left(  \mathbb{A}\left(  x_{o},\xi_{o3}\right)  f\right)  \left(
x,t,L\right)  \right\vert \leq\sum\limits_{n\in\mathbb{Z}\left\backslash
\left\{  -2Q,\cdot\cdot\cdot,2P+1\right\}  \right.  }\left\vert \left(
\mathcal{A}\left(  x_{o},\xi_{o3},n\right)  f\right)  \left(  x,t,L\right)
\right\vert +\left\vert \sum\limits_{n\in\mathbb{Z}}\left(  \mathcal{A}\left(
x_{o},\xi_{o3},n\right)  f\right)  \left(  x,t,L\right)  \right\vert
\tag{6.3.15}\label{6.3.15}%
\end{equation}
and
\begin{equation}
\left\vert \frac{1}{\sqrt{iL+1}}\sum\limits_{n\in\mathbb{Z}}\left(  -1\right)
^{n}~e^{i\left[  \left(  -1\right)  ^{n}x_{3}+2n\frac{\xi_{o3}}{\left\vert
\xi_{o3}\right\vert }\rho\right]  \xi_{3}}~e^{-\frac{1}{4h}\frac{\left(
\left(  -1\right)  ^{n}x_{3}+2n\frac{\xi_{o3}}{\left\vert \xi_{o3}\right\vert
}\rho-x_{o3}-2\xi_{3}hL\right)  ^{2}}{iL+1}}\right\vert \leq\frac{2\sqrt{h}%
}{\rho}\left(  \frac{\sqrt{\pi}}{2}+\frac{\rho}{\sqrt{h}}\right)
\tag{6.3.16}\label{6.3.16}%
\end{equation}
deduced from Appendix B with $z=\frac{4h}{\rho^{2}}\left(  iL+1\right)  $.
This completes the proof.

\bigskip

\bigskip

\subsection{Estimate for $\mathcal{I}_{3}$ (the boundary term with
$\mathbb{A}$)}

\bigskip

We estimate $\mathcal{I}_{3}=-h\int_{\Gamma_{1}\cup\Gamma_{2}}\int_{-T}%
^{T}\left\{  \left(  \int_{0}^{L}\mathbb{A}\left(  x_{o},\xi_{o3}\right)
f_{1}ds\right)  \ell\partial_{\nu}u_{1}+\left(  \int_{0}^{L}\mathbb{A}\left(
x_{o},\xi_{o3}\right)  f_{2}ds\right)  \ell\partial_{\nu}u_{2}\right\}
d\sigma dt$ as follows.

\bigskip

\textbf{Lemma 6.3 .-} \textit{There exists }$c>0$\textit{ such that for any
}$\left(  x_{o},\xi_{o3}\right)  \in\overline{\omega_{o}}\times\left(
2\mathbb{Z}+1\right)  $\textit{ and }$h\in\left(  0,1\right]  $\textit{,
}$L\geq1$\textit{, }$\lambda\geq1$\textit{, }$T>0$\textit{,} \textit{we have}%
\begin{equation}
\left\vert \mathcal{I}_{3}\right\vert \leq chL~e^{-\frac{1}{4h}}%
\sqrt{\mathcal{G}\left(  U,0\right)  }\left(  \int_{\mathbb{R}^{2}}\int
_{\xi_{o3}-1}^{\xi_{o3}+1}\int_{\left\vert \tau\right\vert <\lambda}\left\vert
\widehat{\varphi F}\left(  \xi,\tau\right)  \right\vert d\xi d\tau\right)
\text{ .} \tag{6.4.1}\label{6.4.1}%
\end{equation}

\bigskip

Proof .- First, by (\ref{6.1.7}), we deduce that%
\begin{equation}%
\begin{array}
[c]{ll}
& \quad%
{\displaystyle\int_{\Gamma_{1}\cup\Gamma_{2}}}
\left(  \mathbb{A}\left(  x_{o},\xi_{o3}\right)  f\right)  \left(
x,t,s\right)  \ell\left(  x\right)  \partial_{x_{3}}u_{j}\left(  x,t\right)
d\sigma\\
& \leq%
{\displaystyle\int_{\Gamma_{1}\cup\Gamma_{2}}}
\left\vert \left(  \mathcal{A}\left(  x_{o},\xi_{o3},-2Q\right)  f\right)
\left(  x_{1},x_{2},-\frac{\xi_{o3}}{\left\vert \xi_{o3}\right\vert }%
\rho,t,s\right)  \right\vert \left\vert \ell\partial_{x_{3}}u_{j}\left(
x_{1},x_{2},-\frac{\xi_{o3}}{\left\vert \xi_{o3}\right\vert }\rho,t\right)
\right\vert d\sigma\\
& \quad+%
{\displaystyle\int_{\Gamma_{1}\cup\Gamma_{2}}}
\left\vert \left(  \mathcal{A}\left(  x_{o},\xi_{o3},2P+1\right)  f\right)
\left(  x_{1},x_{2},-\frac{\xi_{o3}}{\left\vert \xi_{o3}\right\vert }%
\rho,t,s\right)  \right\vert \left\vert \ell\partial_{x_{3}}u_{j}\left(
x_{1},x_{2},-\frac{\xi_{o3}}{\left\vert \xi_{o3}\right\vert }\rho,t\right)
\right\vert d\sigma\text{ .}%
\end{array}
\tag{6.4.2}\label{6.4.2}%
\end{equation}
Next, recall that
\begin{equation}%
\begin{array}
[c]{ll}
& \quad\left\vert \left(  \mathcal{A}\left(  x_{o},\xi_{o3},n\right)
f\right)  \left(  x,t,s\right)  \right\vert \\
& \leq\frac{1}{\left(  2\pi\right)  ^{4}}%
{\displaystyle\int_{\mathbb{R}^{2}}}
{\displaystyle\int_{\xi_{o3}-1}^{\xi_{o3}+1}}
{\displaystyle\int_{\left\vert \tau\right\vert <\lambda}}
\left\vert \widehat{\varphi f}\left(  \xi,\tau\right)  \right\vert \left(
\frac{\sqrt{\theta}}{\left(  \sqrt{\left(  hs\right)  ^{2}+\theta^{2}}\right)
^{1/2}}~e^{-\frac{\theta}{4}\frac{\left(  t+2\tau hs\right)  ^{2}}{\left(
hs\right)  ^{2}+\theta^{2}}}\right) \\
& \qquad\left(  \frac{1}{\sqrt{s^{2}+1}}~e^{-\frac{1}{4h}\frac{\left\vert
\left(  x_{1}-x_{o1}-2\xi_{1}hs,x_{2}-x_{o2}-2\xi_{2}hs\right)  \right\vert
^{2}}{s^{2}+1}}\right) \\
& \qquad\frac{1}{\left(  \sqrt{s^{2}+1}\right)  ^{1/2}}\left(  ~e^{-\frac
{1}{4h}\frac{\left(  2n\rho-2\left\vert \xi_{3}\right\vert hs+\frac{\xi_{o3}%
}{\left\vert \xi_{o3}\right\vert }\left[  \left(  -1\right)  ^{n}x_{3}%
-x_{o3}\right]  \right)  ^{2}}{s^{2}+1}}\right)  d\xi d\tau\text{ .}%
\end{array}
\tag{6.4.3}\label{6.4.3}%
\end{equation}
Here, for any $s\in\left[  0,L\right]  $, $h\in\left(  0,1\right]  $,
$x_{3}\in\left[  -\rho,\rho\right]  $, $x_{o3}\in\left[  \rho-2r_{o}%
,\rho-r_{o}\right]  $, $\xi_{3}\in\left(  \xi_{o3}-1,\xi_{o3}+1\right)  $,
$\xi_{o3}\in\left(  2\mathbb{Z}+1\right)  $, we have chosen $\left(
P,Q\right)  \in\mathbb{N}^{2}$ (only depending on $\left(  \xi_{o3},L\right)
$) such that%
\begin{equation}
s^{2}+1\leq\left(  2n\rho-2\left\vert \xi_{3}\right\vert hs+\frac{\xi_{o3}%
}{\left\vert \xi_{o3}\right\vert }\left[  \left(  -1\right)  ^{n}x_{3}%
-x_{o3}\right]  \right)  ^{2}\quad\text{when }n\in\left\{  -2Q,2P+1\right\}
\text{ .} \tag{6.4.4}\label{6.4.4}%
\end{equation}
Indeed, for any $x_{3}\in\left[  -\rho,\rho\right]  $ and $x_{o3}\in\left[
\rho-2r_{o},\rho-r_{o}\right]  $
\begin{equation}%
\begin{array}
[c]{ll}%
\sqrt{s^{2}+1} & \leq\sqrt{L^{2}+1}\leq4P\rho-2\left(  \left\vert \xi
_{o3}\right\vert +1\right)  L\quad\text{from our choice of }P\\
& \leq4P\rho-2\left\vert \xi_{3}\right\vert hs+2\rho-\left\vert x_{3}%
+x_{o3}\right\vert \quad\text{because }2r_{o}\leq2\rho-\left\vert x_{3}%
+x_{o3}\right\vert \\
& \leq\left\vert 4P\rho-2\left\vert \xi_{3}\right\vert hs+2\rho+\frac{\xi
_{o3}}{\left\vert \xi_{o3}\right\vert }\left[  -x_{3}-x_{o3}\right]
\right\vert
\end{array}
\tag{6.4.5}\label{6.4.5}%
\end{equation}
and%
\begin{equation}%
\begin{array}
[c]{ll}%
\sqrt{s^{2}+1} & \leq\sqrt{L^{2}+1}\leq4Q\rho-2\left(  \rho-r_{o}\right)
\quad\text{from our choice of }Q\\
& \leq4Q\rho+2\left\vert \xi_{3}\right\vert hs-\left\vert x_{3}-x_{o3}%
\right\vert \quad\text{because }\left\vert x_{3}-x_{o3}\right\vert
\leq2\left(  \rho-r_{o}\right) \\
& \leq\left\vert -4Q\rho-2\left\vert \xi_{3}\right\vert hs+\frac{\xi_{o3}%
}{\left\vert \xi_{o3}\right\vert }\left[  x_{3}-x_{o3}\right]  \right\vert
\text{ .}%
\end{array}
\tag{6.4.6}\label{6.4.6}%
\end{equation}
So (\ref{6.4.4}) implies that
\begin{equation}
e^{-\frac{1}{4h}\frac{\left(  2n\rho-2\left\vert \xi_{3}\right\vert
hs+\frac{\xi_{o3}}{\left\vert \xi_{o3}\right\vert }\left[  \left(  -1\right)
^{n}x_{3}-x_{o3}\right]  \right)  ^{2}}{s^{2}+1}}\leq e^{-\frac{1}{4h}}%
\quad\text{when }n\in\left\{  -2Q,2P+1\right\}  \text{ .} \tag{6.4.7}%
\label{6.4.7}%
\end{equation}
Therefore, from (\ref{6.4.3}) and (\ref{6.4.7}), for any $s\in\left[
0,L\right]  $, $h\in\left(  0,1\right]  $, $x_{3}\in\left[  -\rho,\rho\right]
$, $x_{o3}\in\left[  \rho-2r_{o},\rho-r_{o}\right]  $, $\xi_{3}\in\left(
\xi_{o3}-1,\xi_{o3}+1\right)  $,
\begin{equation}%
\begin{array}
[c]{ll}
& \quad\left\vert \left(  \mathcal{A}\left(  x_{o},\xi_{o3},n\right)
f\right)  \left(  x,t,s\right)  \right\vert \\
& \leq\frac{1}{\left(  2\pi\right)  ^{4}}\frac{1}{\sqrt{s^{2}+1}}\frac
{1}{\left(  \sqrt{s^{2}+1}\right)  ^{1/2}}e^{-\frac{1}{4h}}\\
& \qquad%
{\displaystyle\int_{\mathbb{R}^{2}}}
{\displaystyle\int_{\xi_{o3}-1}^{\xi_{o3}+1}}
{\displaystyle\int_{\left\vert \tau\right\vert <\lambda}}
\left\vert \widehat{\varphi f}\left(  \xi,\tau\right)  \right\vert \left(
\frac{\sqrt{\theta}}{\left(  \sqrt{\left(  hs\right)  ^{2}+\theta^{2}}\right)
^{1/2}}~e^{-\frac{\theta}{4}\frac{\left(  t+2\tau hs\right)  ^{2}}{\left(
hs\right)  ^{2}+\theta^{2}}}\right)  ~d\xi d\tau\quad\text{when }n\in\left\{
-2Q,2P+1\right\}  \text{ .}%
\end{array}
\tag{6.4.8}\label{6.4.8}%
\end{equation}
On the other hand, by Cauchy-Schwarz inequality
\begin{equation}%
\begin{array}
[c]{ll}
& \quad%
{\displaystyle\int_{-T}^{T}}
e^{-\frac{\theta}{4}\frac{\left(  t+2\tau hs\right)  ^{2}}{\left(  hs\right)
^{2}+\theta^{2}}}\left(
{\displaystyle\int_{\Gamma_{1}\cup\Gamma_{2}}}
\left\vert \ell\left(  x\right)  \partial_{x_{3}}u_{j}\left(  x,t\right)
\right\vert d\sigma\right)  dt\\
& \leq c\left(
{\displaystyle\int_{-T}^{T}}
e^{-\frac{\theta}{4}\frac{\left(  t+2\tau hs\right)  ^{2}}{\left(  hs\right)
^{2}+\theta^{2}}}dt\right)  ^{1/2}\left(
{\displaystyle\int_{-T}^{T}}
{\displaystyle\int_{\Gamma_{1}\cup\Gamma_{2}}}
e^{-\frac{\theta}{4}\frac{\left(  t+2\tau hs\right)  ^{2}}{\left(  hs\right)
^{2}+\theta^{2}}}\left\vert \ell\left(  x\right)  \partial_{x_{3}}u_{j}\left(
x,t\right)  \right\vert ^{2}d\sigma dt\right)  ^{1/2}\\
& \leq c\left(  \frac{2\sqrt{\pi}}{\sqrt{\theta}}\sqrt{\left(  hs\right)
^{2}+\theta^{2}}\right)  ^{1/2}\left(
{\displaystyle\int_{-T}^{T}}
{\displaystyle\int_{\Gamma_{1}\cup\Gamma_{2}}}
e^{-\frac{\theta}{4}\frac{\left(  t+2\tau hs\right)  ^{2}}{\left(  hs\right)
^{2}+\theta^{2}}}\left\vert \ell\left(  x\right)  \partial_{x_{3}}u_{j}\left(
x,t\right)  \right\vert ^{2}d\sigma dt\right)  ^{1/2}\text{ .}%
\end{array}
\tag{6.4.9}\label{6.4.9}%
\end{equation}
Next, by multiplying the equation $\partial_{t}^{2}u_{j}-\Delta u_{j}=0$ by
$g\ell^{2}\nabla u_{j}\cdot W$ where $g\left(  t\right)  =e^{-\frac{\theta}%
{4}\frac{\left(  t+2\tau hs\right)  ^{2}}{\left(  hs\right)  ^{2}+\theta^{2}}%
}$ and $W=W\left(  x\right)  $ is a smooth vector field such that $W=\nu$ on
$\partial\Omega$ (see \cite[page 29]{Li}), we get, after integrations by parts
and by Cauchy-Schwarz inequality, observing that $\ell u_{j}=0$ on
$\partial\Omega$,%
\begin{equation}%
\begin{array}
[c]{ll}%
{\displaystyle\int_{\mathbb{R}}}
{\displaystyle\int_{\Gamma_{1}\cup\Gamma_{2}}}
g\left(  t\right)  \left\vert \ell\left(  x\right)  \partial_{x_{3}}%
u_{j}\left(  x,t\right)  \right\vert ^{2}d\sigma dt & \leq c%
{\displaystyle\int_{\mathbb{R}}}
\left(  g+\frac{d}{dt}g\right)
{\displaystyle\int_{\Omega}}
\left(  \left\vert u_{j}\right\vert ^{2}+\left\vert \nabla u_{j}\right\vert
^{2}+\left\vert \partial_{t}u_{j}\right\vert ^{2}\right)  dxdt\\
& \leq c\left(  \frac{2\sqrt{\pi}}{\sqrt{\theta}}\sqrt{\left(  hs\right)
^{2}+\theta^{2}}\right)  \mathcal{G}\left(  U,0\right)  \text{ .}%
\end{array}
\tag{6.4.10}\label{6.4.10}%
\end{equation}
Therefore, (\ref{6.4.9}) and (\ref{6.4.10}) imply that%
\begin{equation}
\int_{-T}^{T}e^{-\frac{\theta}{4}\frac{\left(  t+2\tau hs\right)  ^{2}%
}{\left(  hs\right)  ^{2}+\theta^{2}}}\left(  \int_{\Gamma_{1}\cup\Gamma_{2}%
}\left\vert \ell\left(  x\right)  \partial_{x_{3}}u_{j}\left(  x,t\right)
\right\vert d\sigma\right)  dt\leq c\sqrt{\left(  hs\right)  ^{2}+\theta^{2}%
}\sqrt{\mathcal{G}\left(  U,0\right)  }\text{ .} \tag{6.4.11}\label{6.4.11}%
\end{equation}
We conclude from (\ref{6.4.2}), (\ref{6.4.8}) and (\ref{6.4.11}) that
\begin{equation}%
\begin{array}
[c]{ll}
& \quad\left\vert h%
{\displaystyle\int_{\Gamma_{1}\cup\Gamma_{2}}}
{\displaystyle\int_{-T}^{T}}
\left(
{\displaystyle\int_{0}^{L}}
\mathbb{A}\left(  x_{o},\xi_{o3}\right)  f\left(  x,t,s\right)  ds\right)
\ell\left(  x\right)  \partial_{\nu}u_{j}\left(  x,t\right)  d\sigma
dt\right\vert \\
& \leq\frac{h}{\left(  2\pi\right)  ^{4}}e^{-\frac{1}{4h}}%
{\displaystyle\int_{\mathbb{R}^{2}}}
{\displaystyle\int_{\xi_{o3}-1}^{\xi_{o3}+1}}
{\displaystyle\int_{\left\vert \tau\right\vert <\lambda}}
\left\vert \widehat{\varphi f}\left(  \xi,\tau\right)  \right\vert d\xi
d\tau\\
& \quad%
{\displaystyle\int_{0}^{L}}
\frac{1}{\sqrt{s^{2}+1}}\frac{1}{\left(  \sqrt{s^{2}+1}\right)  ^{1/2}}%
\frac{\sqrt{\theta}}{\left(  \sqrt{\left(  hs\right)  ^{2}+\theta^{2}}\right)
^{1/2}}%
{\displaystyle\int_{-T}^{T}}
e^{-\frac{\theta}{4}\frac{\left(  t+2\tau hs\right)  ^{2}}{\left(  hs\right)
^{2}+\theta^{2}}}%
{\displaystyle\int_{\Gamma_{1}\cup\Gamma_{2}}}
\left\vert \ell\left(  x\right)  \partial_{x_{3}}u_{j}\left(  x,t\right)
\right\vert d\sigma dtds\\
& \leq chL~e^{-\frac{1}{4h}}\sqrt{\mathcal{G}\left(  U,0\right)  }\left(
{\displaystyle\int_{\mathbb{R}^{2}}}
{\displaystyle\int_{\xi_{o3}-1}^{\xi_{o3}+1}}
{\displaystyle\int_{\left\vert \tau\right\vert <\lambda}}
\left\vert \widehat{\varphi f}\left(  \xi,\tau\right)  \right\vert d\xi
d\tau\right)  \text{ .}%
\end{array}
\tag{6.4.12}\label{6.4.12}%
\end{equation}

\bigskip

\bigskip

\subsection{Estimate for $\mathcal{I}_{4}$ (the boundary term with
$\partial_{x_{3}}\mathbb{B}$)}

\bigskip

We estimate $\mathcal{I}_{4}=h\int_{\Gamma_{1}\cup\Gamma_{2}}\int_{-T}%
^{T}\left(  \int_{0}^{L}\partial_{\nu}\left(  \mathbb{B}\left(  x_{o},\xi
_{o3}\right)  f_{3}\right)  ds\right)  \ell u_{3}d\sigma dt$ as follows.

\bigskip

\textbf{Lemma 6.4 .-} \textit{There exists }$c>0$\textit{ such that for any
}$\left(  x_{o},\xi_{o3}\right)  \in\overline{\omega_{o}}\times\left(
2\mathbb{Z}+1\right)  $\textit{ and }$h\in\left(  0,1\right]  $\textit{,
}$L\geq1$\textit{, }$\lambda\geq1$\textit{, }$T>0$\textit{,} \textit{we have}%
\begin{equation}
\left\vert \mathcal{I}_{4}\right\vert \leq chL~e^{-\frac{1}{ch}}%
\sqrt{\mathcal{G}\left(  U,0\right)  }\left(  \int_{\mathbb{R}^{2}}\int
_{\xi_{o3}-1}^{\xi_{o3}+1}\int_{\left\vert \tau\right\vert <\lambda}\left\vert
\widehat{\varphi F}\left(  \xi,\tau\right)  \right\vert d\xi d\tau\right)
\text{ .} \tag{6.5.1}\label{6.5.1}%
\end{equation}

\bigskip

Proof .- First, by (\ref{6.1.8}), we deduce that
\begin{equation}%
\begin{array}
[c]{ll}
& \quad%
{\displaystyle\int_{\Gamma_{1}\cup\Gamma_{2}}}
\partial_{\nu}\left(  \mathbb{B}\left(  x_{o},\xi_{o3}\right)  f\right)
\left(  x,t,s\right)  \ell\left(  x\right)  u_{3}\left(  x,t\right)  d\sigma\\
& \leq%
{\displaystyle\int_{\Gamma_{1}\cup\Gamma_{2}}}
\left\vert \partial_{x_{3}}\left(  \mathcal{A}\left(  x_{o},\xi_{o3}%
,-2Q\right)  f\right)  \left(  x_{1},x_{2},-\frac{\xi_{o3}}{\left\vert
\xi_{o3}\right\vert }\rho,t,s\right)  \right\vert \left\vert \ell u_{3}\left(
x_{1},x_{2},-\frac{\xi_{o3}}{\left\vert \xi_{o3}\right\vert }\rho,t\right)
\right\vert d\sigma\\
& \quad+%
{\displaystyle\int_{\Gamma_{1}\cup\Gamma_{2}}}
\left\vert \partial_{x_{3}}\left(  \mathcal{A}\left(  x_{o},\xi_{o3}%
,2P+1\right)  f\right)  \left(  x_{1},x_{2},-\frac{\xi_{o3}}{\left\vert
\xi_{o3}\right\vert }\rho,t,s\right)  \right\vert \left\vert \ell u_{3}\left(
x_{1},x_{2},-\frac{\xi_{o3}}{\left\vert \xi_{o3}\right\vert }\rho,t\right)
\right\vert d\sigma\text{ .}%
\end{array}
\tag{6.5.2}\label{6.5.2}%
\end{equation}
Next, recall that
\begin{equation}%
\begin{array}
[c]{ll}
& \quad\left\vert \partial_{x_{3}}\left(  \mathcal{A}\left(  x_{o},\xi
_{o3},n\right)  f\right)  \left(  x,t,s\right)  \right\vert \\
& \leq\frac{1}{\left(  2\pi\right)  ^{4}}%
{\displaystyle\int_{\mathbb{R}^{2}}}
{\displaystyle\int_{\xi_{o3}-1}^{\xi_{o3}+1}}
{\displaystyle\int_{\left\vert \tau\right\vert <\lambda}}
\left\vert \widehat{\varphi f}\left(  \xi,\tau\right)  \right\vert \left(
\frac{\sqrt{\theta}}{\left(  \sqrt{\left(  hs\right)  ^{2}+\theta^{2}}\right)
^{1/2}}~e^{-\frac{\theta}{4}\frac{\left(  t+2\tau hs\right)  ^{2}}{\left(
hs\right)  ^{2}+\theta^{2}}}\right) \\
& \qquad\left(  \frac{1}{\sqrt{s^{2}+1}}~e^{-\frac{1}{4h}\frac{\left\vert
\left(  x_{1}-x_{o1}-2\xi_{1}hs,x_{2}-x_{o2}-2\xi_{2}hs\right)  \right\vert
^{2}}{s^{2}+1}}\right) \\
& \qquad\frac{1}{\left(  \sqrt{s^{2}+1}\right)  ^{1/2}}\left[  \left\vert
\xi_{3}\right\vert +\frac{1}{\sqrt{h}}\right]  \left(  e^{-\frac{1}{8h}%
\frac{\left(  2n\rho-2\left\vert \xi_{3}\right\vert hs+\frac{\xi_{o3}%
}{\left\vert \xi_{o3}\right\vert }\left[  \left(  -1\right)  ^{n}x_{3}%
-x_{o3}\right]  \right)  ^{2}}{s^{2}+1}}\right)  d\xi d\tau\text{ .}%
\end{array}
\tag{6.5.3}\label{6.5.3}%
\end{equation}
Here, for any $s\in\left[  0,L\right]  $, $h\in\left(  0,1\right]  $,
$x_{3}\in\left[  -\rho,\rho\right]  $, $x_{o3}\in\left[  \rho-2r_{o}%
,\rho-r_{o}\right]  $, $\xi_{3}\in\left(  \xi_{o3}-1,\xi_{o3}+1\right)  $,
$\xi_{o3}\in\left(  2\mathbb{Z}+1\right)  $, we have chosen $\left(
P,Q\right)  \in\mathbb{N}^{2}$ (only depending on $\left(  \xi_{o3},L\right)
$) such that%
\begin{equation}
\left(  h\left\vert \xi_{3}\right\vert +1\right)  \left(  s^{2}+1\right)
\leq\left(  2n\rho-2\left\vert \xi_{3}\right\vert hs+\frac{\xi_{o3}%
}{\left\vert \xi_{o3}\right\vert }\left[  \left(  -1\right)  ^{n}x_{3}%
-x_{o3}\right]  \right)  ^{2}\quad\text{when }n\in\left\{  -2Q,2P+1\right\}
\text{ .} \tag{6.5.4}\label{6.5.4}%
\end{equation}
Indeed, for any $x_{3}\in\left[  -\rho,\rho\right]  $ and $x_{o3}\in\left[
\rho-2r_{o},\rho-r_{o}\right]  $
\begin{equation}%
\begin{array}
[c]{ll}%
\sqrt{\left(  h\left\vert \xi_{3}\right\vert +1\right)  \left(  s^{2}%
+1\right)  } & \leq\sqrt{\left(  \left\vert \xi_{o3}\right\vert +2\right)
\left(  L^{2}+1\right)  }\leq4P\rho-2\left(  \left\vert \xi_{o3}\right\vert
+1\right)  L\\
& \leq4P\rho-2\left\vert \xi_{3}\right\vert hs+2\rho-\left\vert x_{3}%
+x_{o3}\right\vert \quad\text{because }0\leq2\rho-\left\vert x_{3}%
+x_{o3}\right\vert \\
& \leq\left\vert 4P\rho-2\left\vert \xi_{3}\right\vert hs+2\rho+\frac{\xi
_{o3}}{\left\vert \xi_{o3}\right\vert }\left[  -x_{3}-x_{o3}\right]
\right\vert
\end{array}
\tag{6.5.5}\label{6.5.5}%
\end{equation}
and%
\begin{equation}%
\begin{array}
[c]{ll}%
\sqrt{\left(  h\left\vert \xi_{3}\right\vert +1\right)  \left(  s^{2}%
+1\right)  } & \leq\sqrt{\left(  \left\vert \xi_{o3}\right\vert +2\right)
\left(  L^{2}+1\right)  }\leq4Q\rho-2\left(  \rho-r_{o}\right) \\
& \leq4Q\rho+2\left\vert \xi_{3}\right\vert hs-\left\vert x_{3}-x_{o3}%
\right\vert \quad\text{because }\left\vert x_{3}-x_{o3}\right\vert
\leq2\left(  \rho-r_{o}\right) \\
& \leq\left\vert -4Q\rho-2\left\vert \xi_{3}\right\vert hs+\frac{\xi_{o3}%
}{\left\vert \xi_{o3}\right\vert }\left[  x_{3}-x_{o3}\right]  \right\vert
\text{ .}%
\end{array}
\tag{6.5.6}\label{6.5.6}%
\end{equation}
So (\ref{6.5.4}) implies that when $n\in\left\{  -2Q,2P+1\right\}  $
\begin{equation}
\left[  \left\vert \xi_{3}\right\vert +\frac{1}{\sqrt{h}}\right]  e^{-\frac
{1}{8h}\frac{\left(  2n\rho-2\left\vert \xi_{3}\right\vert hs+\frac{\xi_{o3}%
}{\left\vert \xi_{o3}\right\vert }\left[  \left(  -1\right)  ^{n}x_{3}%
-x_{o3}\right]  \right)  ^{2}}{s^{2}+1}}\leq\left[  \left\vert \xi
_{3}\right\vert +\frac{1}{h}\right]  e^{-\frac{1}{8}\left(  \left\vert \xi
_{3}\right\vert +\frac{1}{h}\right)  }\leq16e^{-\frac{1}{16h}}\text{ .}
\tag{6.5.7}\label{6.5.7}%
\end{equation}
Therefore, from (\ref{6.5.3}) and (\ref{6.5.7}), for any $s\in\left[
0,L\right]  $, $h\in\left(  0,1\right]  $, $x_{3}\in\left[  -\rho,\rho\right]
$, $x_{o3}\in\left[  \rho-2r_{o},\rho-r_{o}\right]  $, $\xi_{3}\in\left(
\xi_{o3}-1,\xi_{o3}+1\right)  $,
\begin{equation}%
\begin{array}
[c]{ll}
& \quad\left\vert \partial_{x_{3}}\left(  \mathcal{A}\left(  x_{o},\xi
_{o3},n\right)  f\right)  \left(  x,t,s\right)  \right\vert \\
& \leq\frac{1}{\left(  2\pi\right)  ^{4}}\frac{1}{\sqrt{s^{2}+1}}\frac
{1}{\left(  \sqrt{s^{2}+1}\right)  ^{1/2}}ce^{-\frac{1}{ch}}\\
& \qquad%
{\displaystyle\int_{\mathbb{R}^{2}}}
{\displaystyle\int_{\xi_{o3}-1}^{\xi_{o3}+1}}
{\displaystyle\int_{\left\vert \tau\right\vert <\lambda}}
\left\vert \widehat{\varphi f}\left(  \xi,\tau\right)  \right\vert \left(
\frac{\sqrt{\theta}}{\left(  \sqrt{\left(  hs\right)  ^{2}+\theta^{2}}\right)
^{1/2}}~e^{-\frac{\theta}{4}\frac{\left(  t+2\tau hs\right)  ^{2}}{\left(
hs\right)  ^{2}+\theta^{2}}}\right)  ~d\xi d\tau\quad\text{when }n\in\left\{
-2Q,2P+1\right\}  \text{ .}%
\end{array}
\tag{6.5.8}\label{6.5.8}%
\end{equation}
On the other hand, by Cauchy-Schwarz inequality, a trace theorem and
(\ref{6.1.12}), we have
\begin{equation}%
\begin{array}
[c]{ll}
& \quad%
{\displaystyle\int_{-T}^{T}}
e^{-\frac{\theta}{4}\frac{\left(  t+2\tau hs\right)  ^{2}}{\left(  hs\right)
^{2}+\theta^{2}}}\left(
{\displaystyle\int_{\Gamma_{1}\cup\Gamma_{2}}}
\left\vert \ell\left(  x\right)  u_{3}\left(  x,t\right)  \right\vert
d\sigma\right)  dt\\
& \leq c\left(
{\displaystyle\int_{-T}^{T}}
e^{-\frac{\theta}{4}\frac{\left(  t+2\tau hs\right)  ^{2}}{\left(  hs\right)
^{2}+\theta^{2}}}dt\right)  ^{1/2}\left(
{\displaystyle\int_{-T}^{T}}
e^{-\frac{\theta}{4}\frac{\left(  t+2\tau hs\right)  ^{2}}{\left(  hs\right)
^{2}+\theta^{2}}}%
{\displaystyle\int_{\Gamma_{1}\cup\Gamma_{2}}}
\left\vert u_{3}\left(  x,t\right)  \right\vert ^{2}d\sigma dt\right)
^{1/2}\\
& \leq c\left(  \frac{2\sqrt{\pi}}{\sqrt{\theta}}\sqrt{\left(  hs\right)
^{2}+\theta^{2}}\right)  ^{1/2}\left(
{\displaystyle\int_{-T}^{T}}
e^{-\frac{\theta}{4}\frac{\left(  t+2\tau hs\right)  ^{2}}{\left(  hs\right)
^{2}+\theta^{2}}}\left\Vert u_{3}\left(  \cdot,t\right)  \right\Vert
_{H^{1}\left(  \Omega\right)  }^{2}dt\right)  ^{1/2}\\
& \leq c\left(  \sqrt{\left(  hs\right)  ^{2}+\theta^{2}}\right)
\sqrt{\mathcal{G}\left(  U,0\right)  }\text{ .}%
\end{array}
\tag{6.5.9}\label{6.5.9}%
\end{equation}
We conclude from (\ref{6.5.2}), (\ref{6.5.8}) and (\ref{6.5.9}) that
\begin{equation}%
\begin{array}
[c]{ll}
& \quad\left\vert h%
{\displaystyle\int_{\Gamma_{1}\cup\Gamma_{2}}}
{\displaystyle\int_{-T}^{T}}
\left(
{\displaystyle\int_{0}^{L}}
\partial_{\nu}\left(  \mathbb{B}\left(  x_{o},\xi_{o3}\right)  f\right)
\left(  x,t,s\right)  ds\right)  \ell\left(  x\right)  u_{3}\left(
x,t\right)  d\sigma dt\right\vert \\
& \leq\frac{h}{\left(  2\pi\right)  ^{4}}ce^{-\frac{1}{ch}}%
{\displaystyle\int_{\mathbb{R}^{2}}}
{\displaystyle\int_{\xi_{o3}-1}^{\xi_{o3}+1}}
{\displaystyle\int_{\left\vert \tau\right\vert <\lambda}}
\left\vert \widehat{\varphi f}\left(  \xi,\tau\right)  \right\vert d\xi
d\tau\\
& \quad%
{\displaystyle\int_{0}^{L}}
\frac{1}{\sqrt{s^{2}+1}}\frac{1}{\left(  \sqrt{s^{2}+1}\right)  ^{1/2}}%
\frac{\sqrt{\theta}}{\left(  \sqrt{\left(  hs\right)  ^{2}+\theta^{2}}\right)
^{1/2}}%
{\displaystyle\int_{-T}^{T}}
e^{-\frac{\theta}{4}\frac{\left(  t+2\tau hs\right)  ^{2}}{\left(  hs\right)
^{2}+\theta^{2}}}%
{\displaystyle\int_{\Gamma_{1}\cup\Gamma_{2}}}
\left\vert \ell\left(  x\right)  u_{3}\left(  x,t\right)  \right\vert d\sigma
dtds\\
& \leq chL~e^{-\frac{1}{ch}}\sqrt{\mathcal{G}\left(  U,0\right)  }\left(
{\displaystyle\int_{\mathbb{R}^{2}}}
{\displaystyle\int_{\xi_{o3}-1}^{\xi_{o3}+1}}
{\displaystyle\int_{\left\vert \tau\right\vert <\lambda}}
\left\vert \widehat{\varphi f}\left(  \xi,\tau\right)  \right\vert d\xi
d\tau\right)  \text{ .}%
\end{array}
\tag{6.5.10}\label{6.5.10}%
\end{equation}

\bigskip

\bigskip

\subsection{Estimate for $\mathcal{I}_{5}$ (the boundary term on $\left(
\Gamma_{1}\cup\Gamma_{2}\right)  \cap\Theta$)}

\bigskip

We estimate $\mathcal{I}_{5}=-h\int_{\left(  \Gamma_{1}\cup\Gamma_{2}\right)
\cap\Theta}\int_{-T}^{T}\left(  \int_{0}^{L}\mathbb{B}\left(  x_{o},\xi
_{o3}\right)  f_{3}ds\right)  \partial_{\nu}\ell u_{3}d\sigma dt$ as follows.

\bigskip

\textbf{Lemma 6.5 .-} \textit{There exists }$c>0$\textit{ such that for any
}$\left(  x_{o},\xi_{o3}\right)  \in\overline{\omega_{o}}\times\left(
2\mathbb{Z}+1\right)  $\textit{ and }$h\in\left(  0,1\right]  $\textit{,
}$L\geq1$\textit{, }$\lambda\geq1$\textit{, }$T>0$\textit{,} \textit{we have}%
\begin{equation}
\left\vert \mathcal{I}_{5}\right\vert \leq ch\left(  1+\sqrt{hL}\right)
\left\Vert \left(  u_{3},\partial_{t}u_{3}\right)  \right\Vert _{L^{2}\left(
\omega\times\left(  -1-T,T+1\right)  \right)  ^{2}}\left(  \int_{\mathbb{R}%
^{2}}\int_{\xi_{o3}-1}^{\xi_{o3}+1}\int_{\left\vert \tau\right\vert <\lambda
}\left\vert \widehat{\varphi F}\left(  \xi,\tau\right)  \right\vert d\xi
d\tau\right)  \text{ .} \tag{6.6.1}\label{6.6.1}%
\end{equation}

\bigskip

Proof .- Since%
\begin{equation}%
\begin{array}
[c]{ll}
& \quad\left\vert \left(  \mathbb{B}\left(  x_{o},\xi_{o3}\right)  f\right)
\left(  x,t,s\right)  \right\vert \leq\sum\limits_{n\in\mathbb{Z}}\left\vert
\left(  \mathcal{A}\left(  x_{o},\xi_{o3},n\right)  f\right)  \left(
x,t,s\right)  \right\vert \\
& \leq\frac{1}{\left(  2\pi\right)  ^{4}}%
{\displaystyle\int_{\mathbb{R}^{2}}}
{\displaystyle\int_{\xi_{o3}-1}^{\xi_{o3}+1}}
{\displaystyle\int_{\left\vert \tau\right\vert <\lambda}}
\left\vert \widehat{\varphi f}\left(  \xi,\tau\right)  \right\vert \\
& \qquad\left(  \frac{1}{\sqrt{s^{2}+1}}~e^{-\frac{1}{4h}\frac{\left\vert
\left(  x_{1}-x_{o1}-2\xi_{1}hs,x_{2}-x_{o2}-2\xi_{2}hs\right)  \right\vert
^{2}}{s^{2}+1}}\right)  \left(  \frac{\sqrt{\theta}}{\left(  \sqrt{\left(
hs\right)  ^{2}+\theta^{2}}\right)  ^{1/2}}~e^{-\frac{\theta}{4}\frac{\left(
t+2\tau hs\right)  ^{2}}{\left(  hs\right)  ^{2}+\theta^{2}}}\right) \\
& \qquad\left(  \frac{1}{\left(  \sqrt{s^{2}+1}\right)  ^{1/2}}~\sum
\limits_{n\in\mathbb{Z}}e^{-\frac{1}{4h}\frac{\left(  \left(  -1\right)
^{n}x_{3}+2n\frac{\xi_{o3}}{\left\vert \xi_{o3}\right\vert }\rho-x_{o3}%
-2\xi_{3}hs\right)  ^{2}}{s^{2}+1}}\right)  d\xi d\tau\text{ }%
\end{array}
\tag{6.6.2}\label{6.6.2}%
\end{equation}
and
\begin{equation}
\sum\limits_{n\in\mathbb{Z}}e^{-\frac{1}{4h}\frac{\left(  \left(  -1\right)
^{n}x_{3}+2n\frac{\xi_{o3}}{\left\vert \xi_{o3}\right\vert }\rho-x_{o3}%
-2\xi_{3}hs\right)  ^{2}}{s^{2}+1}}\leq2+\frac{\sqrt{\pi}}{\rho}\sqrt{h}%
\sqrt{s^{2}+1} \tag{6.6.3}\label{6.6.3}%
\end{equation}
(see Appendix B with $z=\frac{4h}{\rho^{2}}\left(  s^{2}+1\right)  $), we have%
\begin{equation}%
\begin{array}
[c]{ll}
& \quad\left\vert h%
{\displaystyle\int_{\left(  \Gamma_{1}\cup\Gamma_{2}\right)  \cap\Theta}}
{\displaystyle\int_{-T}^{T}}
\left(
{\displaystyle\int_{0}^{L}}
\mathbb{B}\left(  x_{o},\xi_{o3}\right)  fds\right)  \partial_{\nu}\ell
u_{3}d\sigma dt\right\vert \\
& \leq h\frac{1}{\left(  2\pi\right)  ^{4}}%
{\displaystyle\int_{\mathbb{R}^{2}}}
{\displaystyle\int_{\xi_{o3}-1}^{\xi_{o3}+1}}
{\displaystyle\int_{\left\vert \tau\right\vert <\lambda}}
\left\vert \widehat{\varphi f}\left(  \xi,\tau\right)  \right\vert
{\displaystyle\int_{0}^{L}}
\frac{1}{\sqrt{s^{2}+1}}\frac{\left(  2+\frac{\sqrt{\pi}}{\rho}\sqrt{h}%
\sqrt{s^{2}+1}\right)  }{\left(  \sqrt{s^{2}+1}\right)  ^{1/2}}\\
& \qquad\left(  \frac{\sqrt{\theta}}{\left(  \sqrt{\left(  hs\right)
^{2}+\theta^{2}}\right)  ^{1/2}}%
{\displaystyle\int_{-T}^{T}}
e^{-\frac{\theta}{4}\frac{\left(  t+2\tau hs\right)  ^{2}}{\left(  hs\right)
^{2}+\theta^{2}}}%
{\displaystyle\int_{\left(  \Gamma_{1}\cup\Gamma_{2}\right)  \cap\Theta}}
\left\vert \partial_{\nu}\ell u_{3}\left(  \cdot,t\right)  \right\vert d\sigma
dt\right)  dsd\xi d\tau\\
& \leq ch%
{\displaystyle\int_{\mathbb{R}^{2}}}
{\displaystyle\int_{\xi_{o3}-1}^{\xi_{o3}+1}}
{\displaystyle\int_{\left\vert \tau\right\vert <\lambda}}
\left\vert \widehat{\varphi f}\left(  \xi,\tau\right)  \right\vert d\xi
d\tau\\
& \qquad%
{\displaystyle\int_{0}^{L}}
\frac{1}{\sqrt{s^{2}+1}}\frac{\left(  1+\sqrt{h}\sqrt{s^{2}+1}\right)
}{\left(  \sqrt{s^{2}+1}\right)  ^{1/2}}\frac{\sqrt{\theta}}{\left(
\sqrt{\left(  hs\right)  ^{2}+\theta^{2}}\right)  ^{1/2}}\left(
{\displaystyle\int_{-\infty}^{\infty}}
e^{-\frac{\theta}{2}\frac{t^{2}}{\left(  hs\right)  ^{2}+\theta^{2}}%
}dt\right)  ^{1/2}ds\\
& \qquad\left(
{\displaystyle\int_{-T}^{T}}
{\displaystyle\int_{\left(  \Gamma_{1}\cup\Gamma_{2}\right)  \cap\Theta}}
\left\vert \partial_{\nu}\ell u_{3}\right\vert ^{2}d\sigma dt\right)
^{1/2}\text{ .}%
\end{array}
\tag{6.6.4}\label{6.6.4}%
\end{equation}
Now, we shall treat the term $\left(  \int_{-T}^{T}\int_{\left(  \Gamma
_{1}\cup\Gamma_{2}\right)  \cap\Theta}\left\vert \partial_{\nu}\ell
u_{3}\right\vert ^{2}d\sigma dt\right)  ^{1/2}$ as follows. Let $W=W\left(
x\right)  $ be a smooth vector field such that $W=\nu$ on $\partial\Omega$
(see \cite[page 29]{Li}). Since
\begin{equation}
\operatorname{div}\left(  Wu^{2}\left(  \nabla\ell\cdot W\right)  ^{2}\right)
=2u\left(  \nabla u\cdot W\right)  \left(  \nabla\ell\cdot W\right)
^{2}+u^{2}\nabla\left[  \left(  \nabla\ell\cdot W\right)  ^{2}\right]  \cdot
W+u^{2}\left(  \nabla\ell\cdot W\right)  ^{2}\operatorname{div}W\text{ ,}
\tag{6.6.5}\label{6.6.5}%
\end{equation}
we have the following trace theorem%
\begin{equation}
\int_{\partial\Omega}\left\vert \partial_{\nu}\ell u\right\vert ^{2}%
d\sigma\leq c\int_{\omega}\left\vert u\right\vert ^{2}dx+c\int_{\omega
}\left\vert \nabla u\right\vert ^{2}\left(  \nabla\ell\cdot W\right)
^{2}dx\text{ .} \tag{6.6.6}\label{6.6.6}%
\end{equation}
Next, by multiplying the equation $\partial_{t}^{2}u_{3}-\Delta u_{3}=0$ by
$u_{3}\left(  \nabla\ell\cdot W\right)  ^{2}g$ where $g\in C_{0}^{\infty
}\left(  -1-T,T+1\right)  $ and $g=1$ in $\left(  -T,T\right)  $, we get,
after integrations by parts and by Cauchy-Schwarz inequality, observing that
$\partial_{\nu}u_{3}\partial_{\nu}\ell=0$ on $\partial\Omega$,%
\begin{equation}
\int_{\Omega}\int_{-T}^{T}\left\vert \nabla u_{3}\right\vert ^{2}\left(
\nabla\ell\cdot W\right)  ^{2}dxdt\leq c\int_{\omega}\int_{-1-T}^{T+1}\left(
\left\vert u_{3}\right\vert ^{2}+\left\vert \partial_{t}u_{3}\right\vert
^{2}\right)  dxdt\text{ .} \tag{6.6.7}\label{6.6.7}%
\end{equation}
Therefore, combining (\ref{6.6.6}), (\ref{6.6.7}) and (\ref{6.6.4}), we
conclude that%
\begin{equation}%
\begin{array}
[c]{ll}
& \quad\left\vert h%
{\displaystyle\int_{\left(  \Gamma_{1}\cup\Gamma_{2}\right)  \cap\Theta}}
{\displaystyle\int_{-T}^{T}}
\left(
{\displaystyle\int_{0}^{L}}
\mathbb{B}\left(  x_{o},\xi_{o3}\right)  fds\right)  \partial_{\nu}\ell
u_{3}d\sigma dt\right\vert \\
& \leq ch\left(  1+\sqrt{hL}\right)  \left\Vert \left(  u_{3},\partial
_{t}u_{3}\right)  \right\Vert _{L^{2}\left(  \omega\times\left(
-1-T,T+1\right)  \right)  ^{2}}\left(
{\displaystyle\int_{\mathbb{R}^{2}}}
{\displaystyle\int_{\xi_{o3}-1}^{\xi_{o3}+1}}
{\displaystyle\int_{\left\vert \tau\right\vert <\lambda}}
\left\vert \widehat{\varphi f}\left(  \xi,\tau\right)  \right\vert d\xi
d\tau\right)  \text{ .}%
\end{array}
\tag{6.6.8}\label{6.6.8}%
\end{equation}

\bigskip

\bigskip

\subsection{Estimate for $\mathcal{I}_{6}$ (the term at $t=\pm T$)}

\bigskip

We estimate $\mathcal{I}_{6}=h\int_{\Omega}\left[  \left(  \int_{0}%
^{L}\partial_{t}\left(  V\left(  x_{o},\xi_{o3}\right)  F\right)  \left(
\cdot,t,\cdot\right)  ds\right)  \cdot\ell U\left(  \cdot,t\right)  -\left(
\int_{0}^{L}\left(  V\left(  x_{o},\xi_{o3}\right)  F\right)  \left(
\cdot,t,\cdot\right)  ds\right)  \cdot\ell\partial_{t}U\left(  \cdot,t\right)
\right]  _{-T}^{T}dx$ as follows.

\bigskip

\textbf{Lemma 6.6 .-} \textit{There exists }$c>0$\textit{ such that for any
}$\left(  x_{o},\xi_{o3}\right)  \in\overline{\omega_{o}}\times\left(
2\mathbb{Z}+1\right)  $\textit{ and }$h\in\left(  0,1\right]  $\textit{,
}$L\geq1$\textit{, }$\lambda\geq1$\textit{,} \textit{we have}%
\begin{equation}
\left\vert \mathcal{I}_{6}\right\vert \leq chL\lambda~e^{-\frac{1}{h}}%
\sqrt{\mathcal{G}\left(  U,0\right)  }\left(  \int_{\mathbb{R}^{2}}\int
_{\xi_{o3}-1}^{\xi_{o3}+1}\int_{\left\vert \tau\right\vert <\lambda}\left\vert
\widehat{\varphi F}\left(  \xi,\tau\right)  \right\vert d\xi d\tau\right)
\tag{6.7.1}\label{6.7.1}%
\end{equation}
\textit{when}
\begin{equation}
T=4\left[  \frac{\lambda hL}{\sqrt{2}}+\sqrt{h}L+\frac{\sqrt{2}}{\sqrt{h}%
}\right]  \text{ .} \tag{6.7.2}\label{6.7.2}%
\end{equation}

\bigskip

Proof .- Since%
\begin{equation}%
\begin{array}
[c]{ll}
& \quad\left\vert \partial_{t}\left(  \mathcal{A}\left(  x_{o},\xi
_{o3},n\right)  f\right)  \left(  x,t,s\right)  \right\vert \\
& \leq\frac{1}{\left(  2\pi\right)  ^{4}}%
{\displaystyle\int_{\mathbb{R}^{2}}}
{\displaystyle\int_{\xi_{o3}-1}^{\xi_{o3}+1}}
{\displaystyle\int_{\left\vert \tau\right\vert <\lambda}}
\left\vert \widehat{\varphi f}\left(  \xi,\tau\right)  \right\vert \\
& \qquad\left(  \frac{1}{\sqrt{s^{2}+1}}~e^{-\frac{1}{4h}\frac{\left\vert
\left(  x_{1}-x_{o1}-2\xi_{1}hs,x_{2}-x_{o2}-2\xi_{2}hs\right)  \right\vert
^{2}}{s^{2}+1}}\right)  \left(  \frac{1}{\left(  \sqrt{s^{2}+1}\right)
^{1/2}}~e^{-\frac{1}{4h}\frac{\left(  \left(  -1\right)  ^{n}x_{3}+2n\frac
{\xi_{o3}}{\left\vert \xi_{o3}\right\vert }\rho-x_{o3}-2\xi_{3}hs\right)
^{2}}{s^{2}+1}}\right) \\
& \qquad\left[  \left\vert \tau\right\vert +\frac{1}{2}\frac{\left\vert
t+2\tau hs\right\vert }{\sqrt{\left(  hs\right)  ^{2}+\theta^{2}}}\right]
\left(  \frac{\sqrt{\theta}}{\left(  \sqrt{\left(  hs\right)  ^{2}+\theta^{2}%
}\right)  ^{1/2}}~e^{-\frac{\theta}{4}\frac{\left(  t+2\tau hs\right)  ^{2}%
}{\left(  hs\right)  ^{2}+\theta^{2}}}\right)  d\xi d\tau\text{ ,}%
\end{array}
\tag{6.7.3}\label{6.7.3}%
\end{equation}
we have%
\begin{equation}%
\begin{array}
[c]{ll}
& \quad\left\vert \left(  \mathbb{A}\left(  x_{o},\xi_{o3}\right)  f\right)
\left(  x,\pm T,s\right)  \right\vert +\left\vert \partial_{t}\left(
\mathbb{A}\left(  x_{o},\xi_{o3}\right)  f\right)  \left(  x,\pm T,s\right)
\right\vert \\
& \quad+\left\vert \left(  \mathbb{B}\left(  x_{o},\xi_{o3}\right)  f\right)
\left(  x,\pm T,s\right)  \right\vert +\left\vert \partial_{t}\left(
\mathbb{B}\left(  x_{o},\xi_{o3}\right)  f\right)  \left(  x,\pm T,s\right)
\right\vert \\
& \leq c%
{\displaystyle\int_{\mathbb{R}^{2}}}
{\displaystyle\int_{\xi_{o3}-1}^{\xi_{o3}+1}}
{\displaystyle\int_{\left\vert \tau\right\vert <\lambda}}
\left\vert \widehat{\varphi f}\left(  \xi,\tau\right)  \right\vert \frac
{1}{\sqrt{s^{2}+1}}\left(  \frac{1}{\left(  \sqrt{s^{2}+1}\right)  ^{1/2}}%
\sum\limits_{n\in\mathbb{Z}}e^{-\frac{1}{4h}\frac{\left(  \left(  -1\right)
^{n}x_{3}+2n\frac{\xi_{o3}}{\left\vert \xi_{o3}\right\vert }\rho-x_{o3}%
-2\xi_{3}hs\right)  ^{2}}{s^{2}+1}}\right) \\
& \qquad\left(  1+\lambda+\frac{\sqrt{\theta}}{2}\frac{\left\vert \pm T+2\tau
hs\right\vert }{\sqrt{\left(  hs\right)  ^{2}+\theta^{2}}}\right)  \left(
\frac{\sqrt{\theta}}{\left(  \sqrt{\left(  hs\right)  ^{2}+\theta^{2}}\right)
^{1/2}}~e^{-\frac{\theta}{4}\frac{\left(  \pm T+2\tau hs\right)  ^{2}}{\left(
hs\right)  ^{2}+\theta^{2}}}\right)  d\xi d\tau\\
& \leq c\left(  \frac{1}{\sqrt{s^{2}+1}}\frac{1+\sqrt{h}\sqrt{s^{2}+1}%
}{\left(  \sqrt{s^{2}+1}\right)  ^{1/2}}\frac{1}{\left(  \sqrt{\left(
hs\right)  ^{2}+\theta^{2}}\right)  ^{1/2}}\right)  \lambda\left(
{\displaystyle\int_{\mathbb{R}^{2}}}
{\displaystyle\int_{\xi_{o3}-1}^{\xi_{o3}+1}}
{\displaystyle\int_{\left\vert \tau\right\vert <\lambda}}
\left\vert \widehat{\varphi f}\left(  \xi,\tau\right)  \right\vert
~e^{-\frac{\theta}{8}\frac{\left(  \pm T+2\tau hs\right)  ^{2}}{\left(
hs\right)  ^{2}+\theta^{2}}}d\xi d\tau\right)  \text{ .}%
\end{array}
\tag{6.7.4}\label{6.7.4}%
\end{equation}
On the other hand,
\begin{equation}
e^{-\frac{\theta}{8}\frac{\left(  \pm T+2\tau hs\right)  ^{2}}{\left(
hs\right)  ^{2}+\theta^{2}}}\leq e^{-\frac{\theta T^{2}}{16}\frac{1}{\left(
hL\right)  ^{2}+\theta^{2}}}e^{\frac{\theta}{8}\frac{\left(  2\lambda
hL\right)  ^{2}}{\left(  hL\right)  ^{2}+\theta^{2}}}\quad\forall s\in\left[
0,L\right]  ,\left\vert \tau\right\vert <\lambda\text{ .} \tag{6.7.5}%
\label{6.7.5}%
\end{equation}
Now, when $T=4\left[  \frac{\lambda hL}{\sqrt{2}}+\sqrt{h}L+\frac{\sqrt{2}%
}{\sqrt{h}}\right]  $, then $\frac{1}{2}\left(  \lambda hL\right)  ^{2}%
+\frac{1}{\theta h}\left(  \left(  hL\right)  ^{2}+\theta^{2}\right)
\leq\frac{T^{2}}{16}$ which implies that
\begin{equation}
e^{-\frac{\theta T^{2}}{16}\frac{1}{\left(  hL\right)  ^{2}+\theta^{2}}%
}e^{\frac{\theta}{8}\frac{\left(  2\lambda hL\right)  ^{2}}{\left(  hL\right)
^{2}+\theta^{2}}}\leq e^{-\frac{1}{h}}\text{ .} \tag{6.7.6}\label{6.7.6}%
\end{equation}
In conclusion, combining (\ref{6.7.4}), (\ref{6.7.5}), (\ref{6.7.6}) and
(\ref{6.1.12}), we get
\begin{equation}%
\begin{array}
[c]{ll}
& \left\vert h%
{\displaystyle\int_{\Omega}}
\left[  \left(
{\displaystyle\int_{0}^{L}}
\partial_{t}\left(  V\left(  x_{o},\xi_{o3}\right)  F\right)  \left(
\cdot,t,\cdot\right)  ds\right)  \cdot\ell U\left(  \cdot,t\right)  -\left(
{\displaystyle\int_{0}^{L}}
\left(  V\left(  x_{o},\xi_{o3}\right)  F\right)  \left(  \cdot,t,\cdot
\right)  ds\right)  \cdot\ell\partial_{t}U\left(  \cdot,t\right)  \right]
_{-T}^{T}dx\right\vert \\
& \leq chL\lambda~e^{-\frac{1}{h}}\sqrt{\mathcal{G}\left(  U,0\right)
}\left(
{\displaystyle\int_{\mathbb{R}^{2}}}
{\displaystyle\int_{\xi_{o3}-1}^{\xi_{o3}+1}}
{\displaystyle\int_{\left\vert \tau\right\vert <\lambda}}
\left\vert \widehat{\varphi F}\left(  \xi,\tau\right)  \right\vert d\xi
d\tau\right)  \text{ .}%
\end{array}
\tag{6.7.7}\label{6.7.7}%
\end{equation}

\bigskip

\bigskip

\subsection{Estimate for $\mathcal{I}_{7}$ (the internal term in $\omega$)}

\bigskip

We estimate $\mathcal{I}_{7}=h\int_{\omega}\int_{-T}^{T}\left(  \int_{0}%
^{L}V\left(  x_{o},\xi_{o3}\right)  Fds\right)  \cdot\left[  2\left(
\nabla\ell\cdot\nabla\right)  U+\Delta\ell U\right]  dxdt$ as follows.

\bigskip

\textbf{Lemma 6.7 .-} \textit{There exists }$c>0$\textit{ such that for any
}$\left(  x_{o},\xi_{o3}\right)  \in\overline{\omega_{o}}\times\left(
2\mathbb{Z}+1\right)  $\textit{ and }$h\in\left(  0,1\right]  $\textit{,
}$L\geq1$\textit{, }$\lambda\geq1$\textit{, }$T>0$\textit{,} \textit{we have}%
\begin{equation}
\left\vert \mathcal{I}_{7}\right\vert \leq ch\left(  1+\sqrt{hL}\right)
\left(  \left\Vert \left(  U,\partial_{t}U\right)  \right\Vert _{L^{2}\left(
\omega\times\left(  -1-T,T+1\right)  \right)  ^{6}}\right)  \left(
\int_{\mathbb{R}^{2}}\int_{\xi_{o3}-1}^{\xi_{o3}+1}\int_{\left\vert
\tau\right\vert <\lambda}\left\vert \widehat{\varphi F}\left(  \xi
,\tau\right)  \right\vert d\xi d\tau\right)  \text{ .} \tag{6.8.1}%
\label{6.8.1}%
\end{equation}

\bigskip

Proof .- We start with the third component of $V\left(  x_{o},\xi_{o3}\right)
F$. Since
\begin{equation}%
\begin{array}
[c]{ll}
& \quad\left\vert \left(  \mathbb{B}\left(  x_{o},\xi_{o3}\right)  f\right)
\left(  x,t,s\right)  \right\vert \leq\sum\limits_{n\in\mathbb{Z}}\left\vert
\left(  \mathcal{A}\left(  x_{o},\xi_{o3},n\right)  f\right)  \left(
x,t,s\right)  \right\vert \\
& \leq\frac{1}{\left(  2\pi\right)  ^{4}}%
{\displaystyle\int_{\mathbb{R}^{2}}}
{\displaystyle\int_{\xi_{o3}-1}^{\xi_{o3}+1}}
{\displaystyle\int_{\left\vert \tau\right\vert <\lambda}}
\left\vert \widehat{\varphi f}\left(  \xi,\tau\right)  \right\vert \left(
\frac{\sqrt{\theta}}{\left(  \sqrt{\left(  hs\right)  ^{2}+\theta^{2}}\right)
^{1/2}}~e^{-\frac{\theta}{4}\frac{\left(  t+2\tau hs\right)  ^{2}}{\left(
hs\right)  ^{2}+\theta^{2}}}\right)  d\xi d\tau\text{ }\\
& \qquad\left(  \frac{1}{\sqrt{s^{2}+1}}\right)  \left(  \frac{2+\frac
{\sqrt{\pi}}{\rho}\sqrt{h}\sqrt{s^{2}+1}}{\left(  \sqrt{s^{2}+1}\right)
^{1/2}}\right)  \text{ }%
\end{array}
\tag{6.8.2}\label{6.8.2}%
\end{equation}
(see (\ref{6.6.2})-(\ref{6.6.3})), we have%
\begin{equation}%
\begin{array}
[c]{ll}
& \quad\left\vert h%
{\displaystyle\int_{\omega}}
{\displaystyle\int_{-T}^{T}}
\left(
{\displaystyle\int_{0}^{L}}
\left(  \mathbb{B}\left(  x_{o},\xi_{o3}\right)  f\right)  \left(
x,t,s\right)  ds\right)  \left[  2\nabla\ell\nabla u_{3}+\Delta\ell
u_{3}\right]  \left(  x,t\right)  dxdt\right\vert \\
& \leq ch%
{\displaystyle\int_{\omega}}
{\displaystyle\int_{-T}^{T}}
\left(
{\displaystyle\int_{0}^{L}}
\left\vert \left(  \mathbb{B}\left(  x_{o},\xi_{o3}\right)  f\right)  \left(
x,t,s\right)  \right\vert ds\right)  \left(  \left\vert \nabla\ell\right\vert
\left\vert \nabla u_{3}\right\vert +\left\vert u_{3}\right\vert \right)
\left(  x,t\right)  dxdt\\
& \leq ch\frac{1}{\left(  2\pi\right)  ^{4}}%
{\displaystyle\int_{\mathbb{R}^{2}}}
{\displaystyle\int_{\xi_{o3}-1}^{\xi_{o3}+1}}
{\displaystyle\int_{\left\vert \tau\right\vert <\lambda}}
\left\vert \widehat{\varphi f}\left(  \xi,\tau\right)  \right\vert \int
_{0}^{L}\frac{1}{\sqrt{s^{2}+1}}\frac{\left(  2+\frac{\sqrt{\pi}}{\rho}%
\sqrt{h}\sqrt{s^{2}+1}\right)  }{\left(  \sqrt{s^{2}+1}\right)  ^{1/2}}\\
& \qquad\qquad\left(
{\displaystyle\int_{-T}^{T}}
\left(  \frac{\sqrt{\theta}}{\left(  \sqrt{\left(  hs\right)  ^{2}+\theta^{2}%
}\right)  ^{1/2}}~e^{-\frac{\theta}{4}\frac{\left(  t+2\tau hs\right)  ^{2}%
}{\left(  hs\right)  ^{2}+\theta^{2}}}\right)  \int_{\omega}\left(  \left\vert
\nabla\ell\right\vert \left\vert \nabla u_{3}\right\vert +\left\vert
u_{3}\right\vert \right)  \left(  x,t\right)  dxdt\right)  dsd\xi d\tau
\end{array}
\tag{6.8.3}\label{6.8.3}%
\end{equation}
which implies using Cauchy-Schwarz inequality%
\begin{equation}%
\begin{array}
[c]{ll}
& \quad\left\vert h%
{\displaystyle\int_{\omega}}
{\displaystyle\int_{-T}^{T}}
\left(
{\displaystyle\int_{0}^{L}}
\left(  \mathbb{B}\left(  x_{o},\xi_{o3}\right)  f\right)  \left(
x,t,s\right)  ds\right)  \left[  2\nabla\ell\nabla u_{3}+\Delta\ell
u_{3}\right]  \left(  x,t\right)  dxdt\right\vert \\
& \leq ch\frac{1}{\left(  2\pi\right)  ^{4}}%
{\displaystyle\int_{\mathbb{R}^{2}}}
{\displaystyle\int_{\xi_{o3}-1}^{\xi_{o3}+1}}
{\displaystyle\int_{\left\vert \tau\right\vert <\lambda}}
\left\vert \widehat{\varphi f}\left(  \xi,\tau\right)  \right\vert d\xi
d\tau\\
& \qquad%
{\displaystyle\int_{0}^{L}}
\frac{1}{\sqrt{s^{2}+1}}\frac{\left(  1+\sqrt{h}\sqrt{s^{2}+1}\right)
}{\left(  \sqrt{s^{2}+1}\right)  ^{1/2}}\left(  \frac{\sqrt{\theta}}{\left(
\sqrt{\left(  hs\right)  ^{2}+\theta^{2}}\right)  ^{1/2}}\right)  \left(
{\displaystyle\int_{-\infty}^{\infty}}
e^{-\frac{\theta}{2}\frac{t^{2}}{\left(  hs\right)  ^{2}+\theta^{2}}%
}dt\right)  ^{1/2}ds\\
& \qquad\left(
{\displaystyle\int_{\omega}}
{\displaystyle\int_{-T}^{T}}
\left(  \left\vert \nabla\ell\right\vert ^{2}\left\vert \nabla u_{3}%
\right\vert ^{2}+\left\vert u_{3}\right\vert ^{2}\right)  dxdt\right)
^{1/2}\\
& \leq ch\left(  1+\sqrt{hL}\right)  \left(
{\displaystyle\int_{\omega}}
{\displaystyle\int_{-T}^{T}}
\left(  \left\vert \nabla\ell\right\vert ^{2}\left\vert \nabla u_{3}%
\right\vert ^{2}+\left\vert u_{3}\right\vert ^{2}\right)  dxdt\right)
^{1/2}\left(
{\displaystyle\int_{\mathbb{R}^{2}}}
{\displaystyle\int_{\xi_{o3}-1}^{\xi_{o3}+1}}
{\displaystyle\int_{\left\vert \tau\right\vert <\lambda}}
\left\vert \widehat{\varphi f}\left(  \xi,\tau\right)  \right\vert d\xi
d\tau\right)  \text{ .}%
\end{array}
\tag{6.8.4}\label{6.8.4}%
\end{equation}
Similarly, for any $j\in\left\{  1,2\right\}  $,%
\begin{equation}%
\begin{array}
[c]{ll}
& \quad\left\vert h%
{\displaystyle\int_{\omega}}
{\displaystyle\int_{-T}^{T}}
\left(
{\displaystyle\int_{0}^{L}}
\left(  \mathbb{A}\left(  x_{o},\xi_{o3}\right)  f\right)  \left(
x,t,s\right)  ds\right)  \left[  2\nabla\ell\nabla u_{j}+\Delta\ell
u_{j}\right]  \left(  x,t\right)  dxdt\right\vert \\
& \leq ch\left(  1+\sqrt{hL}\right)  \left(
{\displaystyle\int_{\omega}}
{\displaystyle\int_{-T}^{T}}
\left(  \left\vert \nabla\ell\right\vert ^{2}\left\vert \nabla u_{j}%
\right\vert ^{2}+\left\vert u_{j}\right\vert ^{2}\right)  dxdt\right)
^{1/2}\left(
{\displaystyle\int_{\mathbb{R}^{2}}}
{\displaystyle\int_{\xi_{o3}-1}^{\xi_{o3}+1}}
{\displaystyle\int_{\left\vert \tau\right\vert <\lambda}}
\left\vert \widehat{\varphi f}\left(  \xi,\tau\right)  \right\vert d\xi
d\tau\right)  \text{ .}%
\end{array}
\tag{6.8.5}\label{6.8.5}%
\end{equation}
Now, we shall bound the term $\left(  \int_{\omega}\int_{-T}^{T}\left(
\left\vert \nabla\ell\right\vert ^{2}\left\vert \nabla u_{j}\right\vert
^{2}+\left\vert u_{j}\right\vert ^{2}\right)  dxdt\right)  ^{1/2}$\ for any
$j\in\left\{  1,2,3\right\}  $ by the quantity $\left\Vert \left(
U,\partial_{t}U\right)  \right\Vert _{L^{2}\left(  \omega\times\left(
-1-T,T+1\right)  \right)  ^{6}}$. By multiplying the equation $\partial
_{t}^{2}u_{j}-\Delta u_{j}=0$ by $u_{j}\left\vert \nabla\ell\right\vert ^{2}g$
where $g\in C_{0}^{\infty}\left(  -1-T,T+1\right)  $ and $g=1$ in $\left(
-T,T\right)  $, we get, after integrations by parts and by Cauchy-Schwarz
inequality, observing that $u_{j}\partial_{\nu}u_{j}\left\vert \nabla
\ell\right\vert =0$ on $\partial\Omega$,%
\begin{equation}
\int_{\Omega}\int_{-T}^{T}\left\vert \nabla u_{j}\right\vert ^{2}\left\vert
\nabla\ell\right\vert ^{2}dxdt\leq c\int_{\omega}\int_{-1-T}^{T+1}\left(
\left\vert u_{j}\right\vert ^{2}+\left\vert \partial_{t}u_{j}\right\vert
^{2}\right)  dxdt\text{ ,} \tag{6.8.6}\label{6.8.6}%
\end{equation}
for any $j\in\left\{  1,2,3\right\}  $. This completes the proof.

\bigskip

\bigskip

\subsection{Key inequality}

\bigskip

From now,
\begin{equation}
T=4\left[  \frac{\lambda hL}{\sqrt{2}}+\sqrt{h}L+\frac{\sqrt{2}}{\sqrt{h}%
}\right]  \text{ .} \tag{6.9.1}\label{6.9.1}%
\end{equation}
By (\ref{6.1.13}), (\ref{6.2.1}), (\ref{6.3.1}), (\ref{6.4.1}), (\ref{6.5.1}),
(\ref{6.6.1}), (\ref{6.7.1}) and (\ref{6.8.1}), there exists $c>0$ such that
for any $\left(  x_{o},\xi_{o3}\right)  \in\overline{\omega_{o}}\times\left(
2\mathbb{Z}+1\right)  $ and $h\in\left(  0,1\right]  $, $L\geq1$, $\lambda
\geq1$, we have
\begin{equation}%
\begin{array}
[c]{ll}
& \quad\left\vert
{\displaystyle\int_{\Omega\times\mathbb{R}}}
\left(  \frac{1}{\left(  2\pi\right)  ^{4}}%
{\displaystyle\int_{\mathbb{R}^{2}}}
{\displaystyle\int_{\xi_{o3}-1}^{\xi_{o3}+1}}
{\displaystyle\int_{\left\vert \tau\right\vert <\lambda}}
e^{i\left(  x\xi+t\tau\right)  }~\widehat{\varphi F}\left(  \xi,\tau\right)
d\xi d\tau\right)  \cdot a_{o,\theta}\left(  x,t\right)  \ell\left(  x\right)
U\left(  x,t\right)  dxdt\right\vert \\
& \leq c\left[  \left(  1+hL\lambda\right)  e^{-\frac{1}{ch}}+\frac{1}%
{\sqrt{L}}\right]  \sqrt{\mathcal{G}\left(  U,0\right)  }\left(
{\displaystyle\int_{\mathbb{R}^{2}}}
{\displaystyle\int_{\xi_{o3}-1}^{\xi_{o3}+1}}
{\displaystyle\int_{\left\vert \tau\right\vert <\lambda}}
\left\vert \widehat{\varphi F}\left(  \xi,\tau\right)  \right\vert d\xi
d\tau\right) \\
& \quad+ch\left(  1+\sqrt{hL}\right)  \left(  \left\Vert \left(
U,\partial_{t}U\right)  \right\Vert _{L^{2}\left(  \omega\times\left(
-1-T,T+1\right)  \right)  ^{6}}\right)  \left(
{\displaystyle\int_{\mathbb{R}^{2}}}
{\displaystyle\int_{\xi_{o3}-1}^{\xi_{o3}+1}}
{\displaystyle\int_{\left\vert \tau\right\vert <\lambda}}
\left\vert \widehat{\varphi F}\left(  \xi,\tau\right)  \right\vert d\xi
d\tau\right)  \text{ .}%
\end{array}
\tag{6.9.2}\label{6.9.2}%
\end{equation}
By summing over $\xi_{o3}\in\left(  2\mathbb{Z}+1\right)  $, it implies that
\begin{equation}%
\begin{array}
[c]{ll}
& \quad\left\vert
{\displaystyle\int_{\Omega\times\mathbb{R}}}
\left(  \frac{1}{\left(  2\pi\right)  ^{4}}%
{\displaystyle\int_{\mathbb{R}^{3}}}
{\displaystyle\int_{\left\vert \tau\right\vert <\lambda}}
e^{i\left(  x\xi+t\tau\right)  }~\widehat{\varphi F}\left(  \xi,\tau\right)
d\xi d\tau\right)  \cdot a_{o,\theta}\left(  x,t\right)  \ell\left(  x\right)
U\left(  x,t\right)  dxdt\right\vert \\
& \leq c\left[  \left(  1+hL\lambda\right)  e^{-\frac{1}{ch}}+\frac{1}%
{\sqrt{L}}\right]  \sqrt{\mathcal{G}\left(  U,0\right)  }\left(
{\displaystyle\int_{\mathbb{R}^{3}}}
{\displaystyle\int_{\left\vert \tau\right\vert <\lambda}}
\left\vert \widehat{\varphi F}\left(  \xi,\tau\right)  \right\vert d\xi
d\tau\right) \\
& \quad+ch\left(  1+\sqrt{hL}\right)  \left(  \left\Vert \left(
U,\partial_{t}U\right)  \right\Vert _{L^{2}\left(  \omega\times\left(
-1-T,T+1\right)  \right)  ^{6}}\right)  \left(
{\displaystyle\int_{\mathbb{R}^{3}}}
{\displaystyle\int_{\left\vert \tau\right\vert <\lambda}}
\left\vert \widehat{\varphi F}\left(  \xi,\tau\right)  \right\vert d\xi
d\tau\right)  \text{ .}%
\end{array}
\tag{6.9.3}\label{6.9.3}%
\end{equation}
On the other hand, from (A2) of Appendix A,
\begin{equation}
\int_{\mathbb{R}^{3}}\int_{\left\vert \tau\right\vert <\lambda}\left\vert
\widehat{\varphi F}\left(  \xi,\tau\right)  \right\vert d\xi d\tau\leq
c\sqrt{\lambda}\left(  \frac{\lambda^{2}}{\sqrt{h}}+\frac{1}{h}\right)
\sqrt{\mathcal{G}\left(  \partial_{t}U,0\right)  } \tag{6.9.4}\label{6.9.4}%
\end{equation}
whenever $\varphi F=\varphi_{2}U$ or $\varphi F=\varphi_{1}\partial_{t}^{2}U$.
Therefore, by (\ref{6.9.3}) with (\ref{6.9.4}), we obtain that
\begin{equation}%
\begin{array}
[c]{ll}
& \quad2\left\vert
{\displaystyle\int_{\Omega\times\mathbb{R}}}
\left(  \frac{1}{\left(  2\pi\right)  ^{4}}%
{\displaystyle\int_{\mathbb{R}^{3}}}
{\displaystyle\int_{\left\vert \tau\right\vert <\lambda}}
e^{i\left(  x\xi+t\tau\right)  }~\widehat{\varphi_{2}U}\left(  \xi
,\tau\right)  d\xi d\tau\right)  \cdot a_{o,2}\left(  x,t\right)  \ell\left(
x\right)  U\left(  x,t\right)  dxdt\right\vert \\
& \quad+\left\vert
{\displaystyle\int_{\Omega\times\mathbb{R}}}
\left(  \frac{1}{\left(  2\pi\right)  ^{4}}%
{\displaystyle\int_{\mathbb{R}^{3}}}
{\displaystyle\int_{\left\vert \tau\right\vert <\lambda}}
e^{i\left(  x\xi+t\tau\right)  }~\widehat{\varphi_{1}\partial_{t}^{2}U}\left(
\xi,\tau\right)  d\xi d\tau\right)  \cdot a_{o,1}\left(  x,t\right)
\ell\left(  x\right)  U\left(  x,t\right)  dxdt\right\vert \\
& \leq c\left[  \left(  1+hL\lambda\right)  e^{-\frac{1}{ch}}+\frac{1}%
{\sqrt{L}}\right]  \sqrt{\lambda}\left(  \frac{\lambda^{2}}{\sqrt{h}}+\frac
{1}{h}\right)  \mathcal{G}\left(  \partial_{t}U,0\right) \\
& \quad+ch\left(  1+\sqrt{hL}\right)  \left(  \left\Vert \left(
U,\partial_{t}U\right)  \right\Vert _{L^{2}\left(  \omega\times\left(
-1-T,T+1\right)  \right)  ^{6}}\right)  \sqrt{\lambda}\left(  \frac
{\lambda^{2}}{\sqrt{h}}+\frac{1}{h}\right)  \sqrt{\mathcal{G}\left(
\partial_{t}U,0\right)  }\text{ ,}%
\end{array}
\tag{6.9.5}\label{6.9.5}%
\end{equation}
which is our claim (\ref{6.8}). This completes the proof.

\bigskip

\bigskip

\section*{Appendix A}

\bigskip

The goal of this Appendix A is to prove the two following inequalities (A1)
and (A2) below.

\bigskip

\textbf{Lemma A .-} \textit{Let}
\[
a_{o}\left(  x,t\right)  =e^{-\frac{1}{c_{1}h}\left\vert x-x_{o}\right\vert
^{2}}e^{-\frac{1}{c_{2}}t^{2}}\quad\text{\textit{and}\quad}\varphi\left(
x,t\right)  =\phi\left(  x\right)  e^{-\frac{1}{c_{3}h}\left\vert
x-x_{o}\right\vert ^{2}}e^{-\frac{1}{c_{4}}t^{2}}%
\]
\textit{for some }$c_{1},c_{2},c_{3},c_{4}>0$\textit{ and }$\phi\in
C_{0}^{\infty}\left(  \Omega\right)  $\textit{. Let }$\ell\in C^{\infty
}\left(  \mathbb{R}^{3}\right)  $\textit{ be such that }$0\leq\ell\left(
x\right)  \leq1$\textit{. There exists }$c>0$\textit{ such that for any }%
$h\in\left(  0,1\right]  $\textit{, }$\lambda\geq1$\textit{ and any }$u\in
C^{1}\left(  \mathbb{R},H^{1}\left(  \Omega\right)  \right)  \cap C^{2}\left(
\mathbb{R},L^{2}\left(  \Omega\right)  \right)  $\textit{ satisfying }%
\[
\partial_{t}^{2}u-\Delta u=0\quad\text{\textit{in}}~\Omega\times
\mathbb{R}\text{ ,}%
\]
\textit{and}%
\[
\exists R_{j}>0\qquad\left\Vert \partial_{t}^{j}u\left(  \cdot,t\right)
\right\Vert _{L^{2}\left(  \Omega\right)  }^{2}\leq R_{j}\quad
\text{\textit{for} }j\in\left\{  0,1,2\right\}  \text{,\quad}\left\Vert \nabla
u\left(  \cdot,t\right)  \right\Vert _{L^{2}\left(  \Omega\right)  }^{2}\leq
R_{1}\text{ ,}%
\]
\textit{we have}%
\begin{equation}%
\begin{array}
[c]{ll}
& \quad\left\vert
{\displaystyle\int_{\Omega\times\mathbb{R}}}
\left(  \frac{1}{\left(  2\pi\right)  ^{4}}%
{\displaystyle\int_{\mathbb{R}^{3}}}
{\displaystyle\int_{\left\vert \tau\right\vert \geq\lambda}}
e^{i\left(  x\xi+t\tau\right)  }~\widehat{\varphi f}\left(  \xi,\tau\right)
d\xi d\tau\right)  a_{o}\left(  x,t\right)  \ell\left(  x\right)  u\left(
x,t\right)  dxdt\right\vert \\
& \leq c\sqrt{\frac{1}{\lambda}}\left(  R_{0}+R_{1}\right)  \left(
R_{0}+R_{2}\right)
\end{array}
\tag{A1}\label{A1}%
\end{equation}
\textit{and}
\begin{equation}
\left(  \int_{\mathbb{R}^{3}}\int_{\left\vert \tau\right\vert <\lambda
}\left\vert \widehat{\varphi f}\left(  \xi,\tau\right)  \right\vert d\xi
d\tau\right)  \leq c\sqrt{\lambda}\left(  \left(  \lambda^{2}+\frac{1}%
{h}\right)  R_{2}+\frac{\lambda^{2}}{\sqrt{h}}R_{1}+\frac{1}{h}R_{0}\right)
\tag{A2}\label{A2}%
\end{equation}
\textit{whenever }$f=u$\textit{ or }$f=\partial_{t}^{2}u$\textit{.}

\bigskip

\bigskip

\textbf{Proof of (\ref{A1}).} Introduce
\[
\mathcal{R}\left(  f\right)  =\int_{\Omega\times\mathbb{R}}\left(  \frac
{1}{\left(  2\pi\right)  ^{4}}\int_{\mathbb{R}^{3}}\int_{\left\vert
\tau\right\vert \geq\lambda}e^{i\left(  x\xi+t\tau\right)  }~\widehat{\varphi
f}\left(  \xi,\tau\right)  d\xi d\tau\right)  a_{o}\left(  x,t\right)
\ell\left(  x\right)  u\left(  x,t\right)  dxdt\text{ .}%
\]
Thus,
\[%
\begin{array}
[c]{ll}%
\left\vert \mathcal{R}\left(  f\right)  \right\vert  & =\left\vert
{\displaystyle\int_{\Omega\times\mathbb{R}}}
a_{o}\left(  x,t\right)  \ell\left(  x\right)  u\left(  x,t\right)
\partial_{t}\left(  \frac{1}{\left(  2\pi\right)  ^{4}}%
{\displaystyle\int_{\mathbb{R}^{3}}}
{\displaystyle\int_{\left\vert \tau\right\vert \geq\lambda}}
\frac{1}{i\tau}~e^{i\left(  x\xi+t\tau\right)  }~\widehat{\varphi f}\left(
\xi,\tau\right)  d\xi d\tau\right)  dxdt\right\vert \text{ ,}\\
& =\left\vert
{\displaystyle\int_{\Omega\times\mathbb{R}}}
\ell\left(  x\right)  \partial_{t}\left(  a_{o}u\left(  x,t\right)  \right)
\left(  \frac{1}{2\pi}%
{\displaystyle\int_{\left\vert \tau\right\vert \geq\lambda}}
\frac{1}{i\tau}~e^{it\tau}\left[
{\displaystyle\int_{\mathbb{R}}}
e^{-i\theta\tau}\left(  \varphi f\right)  \left(  x,\theta\right)
d\theta\right]  d\tau\right)  dxdt\right\vert \text{ .}%
\end{array}
\]
It follows using Cauchy-Schwarz inequality and Parseval identity that%
\[%
\begin{array}
[c]{ll}%
\left\vert \mathcal{R}\left(  f\right)  \right\vert  & \leq%
{\displaystyle\int_{\Omega\times\mathbb{R}}}
\left\vert \ell\left(  x\right)  \partial_{t}\left(  a_{o}u\left(  x,t\right)
\right)  \right\vert \left(  \frac{1}{2\pi}\left[
{\displaystyle\int_{\left\vert \tau\right\vert \geq\lambda}}
\frac{1}{\tau^{2}}d\tau\right]  ^{1/2}\left[
{\displaystyle\int_{\mathbb{R}}}
\left\vert
{\displaystyle\int_{\mathbb{R}}}
e^{-i\theta\tau}\left(  \varphi f\right)  \left(  x,\theta\right)
d\theta\right\vert ^{2}d\tau\right]  ^{1/2}\right)  dxdt\\
& \leq%
{\displaystyle\int_{\Omega\times\mathbb{R}}}
\left\vert \partial_{t}\left(  a_{o}u\left(  x,t\right)  \right)  \right\vert
\left(  \frac{1}{2\pi}\left[
{\displaystyle\int_{\left\vert \tau\right\vert \geq\lambda}}
\frac{1}{\tau^{2}}d\tau\right]  ^{1/2}\left[  2\pi%
{\displaystyle\int_{\mathbb{R}}}
\left\vert \left(  \varphi f\right)  \left(  x,\theta\right)  \right\vert
^{2}d\theta\right]  ^{1/2}\right)  dxdt\\
& \leq%
{\displaystyle\int_{\Omega\times\mathbb{R}}}
\left\vert \partial_{t}\left(  a_{o}u\left(  x,t\right)  \right)  \right\vert
\left(  \frac{1}{\sqrt{2\pi}}\sqrt{\frac{2}{\lambda}}\left\Vert \left(
\varphi f\right)  \left(  x,\cdot\right)  \right\Vert _{L^{2}\left(
\mathbb{R}\right)  }\right)  dxdt\\
& \leq\frac{1}{\sqrt{\pi}}\sqrt{\frac{1}{\lambda}}~%
{\displaystyle\int_{\mathbb{R}}}
\left\Vert \partial_{t}\left(  a_{o}u\right)  \left(  \cdot,t\right)
\right\Vert _{L^{2}\left(  \Omega\right)  }dt\left\Vert \varphi f\right\Vert
_{L^{2}\left(  \Omega\times\mathbb{R}\right)  }\text{ .}%
\end{array}
\]
Since we have the following estimates%
\[%
{\displaystyle\int_{\mathbb{R}}}
e^{-\frac{t^{2}}{c}}\int_{\Omega}\left(  \left\vert u\left(  x,t\right)
\right\vert ^{2}+\left\vert \partial_{t}^{2}u\left(  x,t\right)  \right\vert
^{2}\right)  dxdt\leq c\left(  R_{0}+R_{2}\right)  \text{ ,}%
\]%
\[%
\begin{array}
[c]{ll}%
{\displaystyle\int_{\mathbb{R}}}
\left\Vert \partial_{t}\left(  a_{o}u\right)  \left(  \cdot,t\right)
\right\Vert _{L^{2}\left(  \Omega\right)  }dt & \leq%
{\displaystyle\int_{\mathbb{R}}}
\left[
{\displaystyle\int_{\Omega}}
\left\vert \partial_{t}a_{o}u\left(  x,t\right)  \right\vert ^{2}dx\right]
^{1/2}dt+%
{\displaystyle\int_{\mathbb{R}}}
\left[
{\displaystyle\int_{\Omega}}
\left\vert a_{o}\partial_{t}u\left(  x,t\right)  \right\vert ^{2}dx\right]
^{1/2}dt\\
& \leq%
{\displaystyle\int_{\mathbb{R}}}
\left\vert \frac{2t}{c_{2}}\right\vert e^{-\frac{1}{c_{2}}t^{2}}\left[
{\displaystyle\int_{\Omega}}
\left\vert u\left(  x,t\right)  \right\vert ^{2}dx\right]  ^{1/2}dt+%
{\displaystyle\int_{\mathbb{R}}}
e^{-\frac{1}{c_{2}}t^{2}}\left[
{\displaystyle\int_{\Omega}}
\left\vert \partial_{t}u\left(  x,t\right)  \right\vert ^{2}dx\right]
^{1/2}dt\\
& \leq c\left(  R_{0}+R_{1}\right)  \text{ ,}%
\end{array}
\]
we conclude that%
\[
\left\vert \mathcal{R}\left(  u\right)  \right\vert +\left\vert \mathcal{R}%
\left(  \partial_{t}^{2}u\right)  \right\vert \leq c\sqrt{\frac{1}{\lambda}%
}\left(  R_{0}+R_{1}\right)  \left(  R_{0}+R_{2}\right)  \text{ .}%
\]
That completes the proof of (\ref{A1}).

\bigskip

\textbf{Proof of (\ref{A2}).} We estimate $\int_{\mathbb{R}^{3}}%
\int_{\left\vert \tau\right\vert <\lambda}\left\vert \widehat{\varphi
f}\left(  \xi,\tau\right)  \right\vert d\xi d\tau$ where $f$ solves
$\partial_{t}^{2}f-\Delta f=0$ in $\Omega\times\mathbb{R}$. By Cauchy-Schwarz
inequality and Parseval identity,%
\[%
\begin{array}
[c]{ll}%
{\displaystyle\int_{\mathbb{R}^{3}}}
{\displaystyle\int_{\left\vert \tau\right\vert <\lambda}}
\left\vert \widehat{\varphi f}\left(  \xi,\tau\right)  \right\vert d\xi d\tau
& =%
{\displaystyle\int_{\mathbb{R}^{3}}}
{\displaystyle\int_{\left\vert \tau\right\vert <\lambda}}
\frac{1}{1+\left\vert \xi\right\vert ^{2}}\left\vert \left(  1+\left\vert
\xi\right\vert ^{2}\right)  \widehat{\varphi f}\left(  \xi,\tau\right)
\right\vert d\xi d\tau\\
& \leq%
{\displaystyle\int_{\left\vert \tau\right\vert <\lambda}}
\left[
{\displaystyle\int_{\mathbb{R}^{3}}}
\frac{1}{\left(  1+\left\vert \xi\right\vert ^{2}\right)  ^{2}}d\xi\right]
^{1/2}\left[
{\displaystyle\int_{\mathbb{R}^{3}}}
\left\vert \left(  \widehat{\left(  1-\Delta\right)  \left(  \varphi f\right)
}\right)  \left(  \xi,\tau\right)  \right\vert ^{2}d\xi\right]  ^{1/2}d\tau\\
& \leq\pi^{2}\sqrt{\lambda}\left[
{\displaystyle\int_{\mathbb{R}^{3}}}
{\displaystyle\int_{\left\vert \tau\right\vert <\lambda}}
\left\vert \left(  \widehat{\left(  1-\Delta\right)  \left(  \varphi f\right)
}\right)  \left(  \xi,\tau\right)  \right\vert ^{2}d\xi d\tau\right]
^{1/2}\text{ .}%
\end{array}
\]
On the other hand, remark that $\partial_{x}^{j}\varphi\left(  x,t\right)
=\frac{1}{h^{j/2}}\phi_{j}\left(  x\right)  e^{-\frac{1}{c_{3}h}\left\vert
x-x_{o}\right\vert ^{2}}e^{-\frac{1}{c_{4}}t^{2}}$ for some $\phi_{j}\in
C_{0}^{\infty}\left(  \Omega\right)  $. Since $\Delta\left(  \varphi u\right)
=\varphi\Delta u+\Delta\varphi u+2\nabla\varphi\nabla u$, we obtain when
$f=u$, using Parseval identity and the last remark%
\[%
\begin{array}
[c]{ll}
& \quad%
{\displaystyle\int_{\mathbb{R}^{3}}}
{\displaystyle\int_{\left\vert \tau\right\vert <\lambda}}
\left\vert \widehat{\varphi u}\left(  \xi,\tau\right)  \right\vert d\xi
d\tau\\
& \leq c\sqrt{\lambda}\left(  \left\Vert \varphi\Delta u\right\Vert
_{L^{2}\left(  \Omega\times\mathbb{R}\right)  }+\left\Vert \Delta\varphi
u\right\Vert _{L^{2}\left(  \Omega\times\mathbb{R}\right)  }+\left\Vert
\nabla\varphi\nabla u\right\Vert _{L^{2}\left(  \Omega\times\mathbb{R}\right)
}+\left\Vert \varphi u\right\Vert _{L^{2}\left(  \Omega\times\mathbb{R}%
\right)  }\right) \\
& \leq c\sqrt{\lambda}\left(  R_{2}+\frac{1}{h}R_{0}+\frac{1}{\sqrt{h}}%
R_{1}\right)  \text{ .}%
\end{array}
\]
Since
\[%
\begin{array}
[c]{ll}%
\Delta\left(  \varphi\partial_{t}^{2}u\right)  & =\varphi\partial_{t}%
^{2}\Delta u+\Delta\varphi\Delta u+2\nabla\varphi\nabla\partial_{t}^{2}u\\
& =\partial_{t}^{2}\left(  \varphi\Delta u\right)  -2\partial_{t}\left(
\partial_{t}\varphi\Delta u\right)  +\left(  \partial_{t}^{2}\varphi
+\Delta\varphi\right)  \Delta u+2\nabla\varphi\nabla\partial_{t}^{2}u\\
& =\partial_{t}^{2}\left(  \varphi\Delta u\right)  -2\partial_{t}\left(
\partial_{t}\varphi\Delta u\right)  +\left(  \partial_{t}^{2}\varphi
+\Delta\varphi\right)  \Delta u\\
& \quad+2\partial_{t}^{2}\left(  \nabla\varphi\nabla u\right)  -4\partial
_{t}\left(  \partial_{t}\nabla\varphi\nabla u\right)  +2\partial_{t}^{2}%
\nabla\varphi\nabla u\text{ ,}%
\end{array}
\]
we obtain when $f=\partial_{t}^{2}u$, using the fact that $\left\vert
\tau\right\vert <\lambda$, Parseval identity and the above remark,
\[%
\begin{array}
[c]{ll}
& \quad%
{\displaystyle\int_{\mathbb{R}^{3}}}
{\displaystyle\int_{\left\vert \tau\right\vert <\lambda}}
\left\vert \widehat{\varphi\partial_{t}^{2}u}\left(  \xi,\tau\right)
\right\vert ~d\xi d\tau\\
& \leq c\sqrt{\lambda}\left(  \lambda^{2}\left\Vert \varphi\Delta u\right\Vert
_{L^{2}\left(  \Omega\times\mathbb{R}\right)  }+\lambda\left\Vert \partial
_{t}\varphi\Delta u\right\Vert _{L^{2}\left(  \Omega\times\mathbb{R}\right)
}+\left\Vert \left(  \partial_{t}^{2}\varphi+\Delta\varphi\right)  \Delta
u\right\Vert _{L^{2}\left(  \Omega\times\mathbb{R}\right)  }+\left\Vert
\varphi\Delta u\right\Vert _{L^{2}\left(  \Omega\times\mathbb{R}\right)
}\right) \\
& \quad+c\sqrt{\lambda}\left(  \lambda^{2}\left\Vert \nabla\varphi\nabla
u\right\Vert _{L^{2}\left(  \Omega\times\mathbb{R}\right)  }+\lambda\left\Vert
\partial_{t}\nabla\varphi\nabla u\right\Vert _{L^{2}\left(  \Omega
\times\mathbb{R}\right)  }+\left\Vert \partial_{t}^{2}\nabla\varphi\nabla
u\right\Vert _{L^{2}\left(  \Omega\times\mathbb{R}\right)  }\right) \\
& \leq c\sqrt{\lambda}\left(  \left(  \lambda^{2}+\frac{1}{h}\right)
R_{2}+\frac{\lambda^{2}}{\sqrt{h}}R_{1}\right)  \text{ .}%
\end{array}
\]
We conclude that there exists $c>0$ such that for any $h\in\left(  0,1\right]
$ and $\lambda\geq1$,
\[
\int_{\mathbb{R}^{3}}\int_{\left\vert \tau\right\vert <\lambda}\left\vert
\widehat{\varphi u}\left(  \xi,\tau\right)  \right\vert ~d\xi d\tau
+\int_{\mathbb{R}^{3}}\int_{\left\vert \tau\right\vert <\lambda}\left\vert
\widehat{\varphi\partial_{t}^{2}u}\left(  \xi,\tau\right)  \right\vert d\xi
d\tau\leq c\sqrt{\lambda}\left(  \left(  \lambda^{2}+\frac{1}{h}\right)
R_{2}+\frac{\lambda^{2}}{\sqrt{h}}R_{1}+\frac{1}{h}R_{0}\right)  \text{ .}%
\]
That completes the proof of (\ref{A2}).

\bigskip

\bigskip

\section*{Appendix B}

\bigskip

The goal of this Appendix B is to prove the two following inequalities.

\bigskip

\textbf{Lemma B .-} \textit{For any }$x,y,C,R\in\mathbb{R}$\textit{, any
}$z\in\mathbb{C},\operatorname{Re}z>0$\textit{,}
\[
\left\vert \frac{1}{\sqrt{z}}\sum_{n\in\mathbb{Z}}\left(  -1\right)
^{n}e^{-\frac{1}{z}\left(  2n+C\left(  -1\right)  ^{n}+R\right)  ^{2}}%
e^{inx}e^{i\left(  -1\right)  ^{n}y}\right\vert \leq\frac{\sqrt{\pi}}{2}%
+\frac{2}{\sqrt{\operatorname{Re}z}}\text{ ,}%
\]%
\[
\left\vert \frac{1}{\sqrt{z}}\sum_{n\in\mathbb{Z}}e^{-\frac{1}{z}\left(
2n+C\left(  -1\right)  ^{n}+R\right)  ^{2}}e^{inx}e^{i\left(  -1\right)
^{n}y}\right\vert \leq\frac{\sqrt{\pi}}{2}+\frac{2}{\sqrt{\operatorname{Re}z}%
}\text{ .}%
\]

\bigskip

Proof .- First we recall the Poisson summation formula. Let $u\in C^{2}\left(
\mathbb{R},\mathbb{C}\right)  $ be such that for any $k\in\left\{
0,1,2\right\}  $, the functions $x\longmapsto\left(  1+x^{2}\right)
u^{\left(  k\right)  }\left(  x\right)  $ are bounded on $\mathbb{R}$. Then
for any $x\in\mathbb{R}$,
\[
\sum_{n\in\mathbb{Z}}u\left(  x+n\right)  =\sum_{n\in\mathbb{Z}}\widehat
{u}\left(  2\pi n\right)  e^{2\pi inx}\text{\quad where\quad}\widehat
{u}\left(  2\pi n\right)  =\int_{\mathbb{R}}u\left(  t\right)  e^{-2\pi
int}dt\text{ .}%
\]
Next, by choosing $u\left(  x\right)  =v\left(  x\right)  e^{-2\pi iBx}$ for
some $B\in\mathbb{R}$ and $v\in C^{2}\left(  \mathbb{R},\mathbb{C}\right)  $
such that for any $k\in\left\{  0,1,2\right\}  $, the functions $x\longmapsto
\left(  1+x^{2}\right)  v^{\left(  k\right)  }\left(  x\right)  $ are bounded
on $\mathbb{R}$, we obtain that for any $x,B\in\mathbb{R}$,
\[
\sum_{n\in\mathbb{Z}}\widehat{v}\left(  2\pi\left(  n+B\right)  \right)
e^{2\pi inx}=\sum_{n\in\mathbb{Z}}v\left(  x+n\right)  e^{-2\pi iB\left(
x+n\right)  }\text{\quad where\quad}\widehat{v}\left(  2\pi\left(  n+B\right)
\right)  =\int_{\mathbb{R}}v\left(  t\right)  e^{-2\pi i\left(  n+B\right)
t}dt\text{ .}%
\]
Now, for any $z\in\mathbb{C}$ such that $\operatorname{Re}z>0$, we take
$v\left(  x\right)  =e^{-\frac{z}{2}x^{2}}$ in order that $\widehat{v}\left(
2\pi\left(  n+B\right)  \right)  =\frac{\sqrt{2\pi}}{\sqrt{z}}e^{-\frac{1}%
{2z}\left(  2\pi\left(  n+B\right)  \right)  ^{2}}$. Thus, after simple
changes, the following formula holds for any $x,B\in\mathbb{R}$, any
$z\in\mathbb{C},\operatorname{Re}z>0,a>0$,%
\[
\frac{1}{\sqrt{z}}\sum_{n\in\mathbb{Z}}e^{-\frac{a}{z}\left(  n+B\right)
^{2}}e^{i2nx}=\frac{\sqrt{\pi}}{\sqrt{a}}\sum_{n\in\mathbb{Z}}e^{-\frac{z}%
{a}\left(  x+\pi n\right)  ^{2}}e^{-i2B\left(  x+\pi n\right)  }\text{ .}%
\]
Finally, we deduce that for any $x,y,C,R\in\mathbb{R}$, any $z\in
\mathbb{C},\operatorname{Re}z>0$,
\[%
\begin{array}
[c]{ll}
& \quad\left\vert \frac{1}{\sqrt{z}}\sum_{n\in\mathbb{Z}}\left(  -1\right)
^{n}e^{-\frac{1}{z}\left(  2n+C\left(  -1\right)  ^{n}+R\right)  ^{2}}%
e^{inx}e^{i\left(  -1\right)  ^{n}y}\right\vert \\
& =\left\vert \frac{1}{\sqrt{z}}\sum_{n\in\mathbb{Z}}e^{-\frac{1}{z}\left(
4n+C+R\right)  ^{2}}e^{i2nx}e^{iy}-\frac{1}{\sqrt{z}}\sum_{n\in\mathbb{Z}%
}e^{-\frac{1}{z}\left(  4n+2-C+R\right)  ^{2}}e^{i\left(  2n+1\right)
x}e^{-iy}\right\vert \\
& =\left\vert e^{iy}\frac{1}{\sqrt{z}}\sum_{n\in\mathbb{Z}}e^{-\frac{4^{2}}%
{z}\left(  n+\frac{C+R}{4}\right)  ^{2}}e^{i2nx}-e^{-iy}e^{ix}\frac{1}%
{\sqrt{z}}\sum_{n\in\mathbb{Z}}e^{-\frac{4^{2}}{z}\left(  n+\frac{2-C+R}%
{4}\right)  ^{2}}e^{i2nx}\right\vert \\
& =\left\vert e^{iy}\frac{\sqrt{\pi}}{4}\sum_{n\in\mathbb{Z}}e^{-\frac
{z}{4^{2}}\left(  x+\pi n\right)  ^{2}}e^{-i\frac{C+R}{2}\left(  x+\pi
n\right)  }-e^{-iy}\frac{\sqrt{\pi}}{4}\sum_{n\in\mathbb{Z}}e^{-\frac{z}%
{4^{2}}\left(  x+\pi n\right)  ^{2}}e^{-i\frac{2-C+R}{2}\left(  x+\pi
n\right)  }\right\vert \\
& \leq\frac{\sqrt{\pi}}{2}\sum_{n\in\mathbb{Z}}e^{-\frac{\operatorname{Re}%
z}{4^{2}}\pi^{2}\left(  \frac{x}{\pi}+n\right)  ^{2}}\leq\frac{\sqrt{\pi}}%
{2}+\frac{2}{\sqrt{\operatorname{Re}z}}\text{ ,}%
\end{array}
\]
and similarly
\[%
\begin{array}
[c]{ll}
& \quad\left\vert \frac{1}{\sqrt{z}}\sum_{n\in\mathbb{Z}}e^{-\frac{1}%
{z}\left(  2n+C\left(  -1\right)  ^{n}+R\right)  ^{2}}e^{inx}e^{i\left(
-1\right)  ^{n}y}\right\vert \\
& =\left\vert e^{iy}\frac{1}{\sqrt{z}}\sum_{n\in\mathbb{Z}}e^{-\frac{1}%
{z}\left(  4n+C+R\right)  ^{2}}e^{i2nx}+e^{-iy}e^{ix}\frac{1}{\sqrt{z}}%
\sum_{n\in\mathbb{Z}}e^{-\frac{1}{z}\left(  4n+2-C+R\right)  ^{2}}%
e^{i2nx}\right\vert \\
& \leq\frac{\sqrt{\pi}}{2}+\frac{2}{\sqrt{\operatorname{Re}z}}\text{ .}%
\end{array}
\]

\bigskip

\bigskip

\end{document}